\documentclass[12pt]{article}

\usepackage{makeidx}
\usepackage{latexsym}
\usepackage{amsfonts}
\usepackage{amsmath}
\usepackage{amsthm}
\usepackage{amssymb}
\usepackage{amscd}
\usepackage[dvips]{graphicx}
\usepackage{color} 
\usepackage{eepic} 
\usepackage{setspace}
\usepackage{calrsfs}

\newdimen\Squaresize \Squaresize=20pt
\newdimen\thickness \thickness=1pt         
                                                    
\def\Square#1{\hbox{\vrule width \thickness
   \vbox to \Squaresize{\hrule height \thickness\vss                            
      \hbox to \Squaresize{\hss#1\hss}
   \vss\hrule height\thickness} 
\unskip\vrule width \thickness} 
\kern-\thickness}                                                            
                               
\def\vsquare#1{\vbox{\Square{$#1$}}\kern-\thickness}
\def\blank{\omit\hskip\Squaresize}

\def\fibyoung#1{\let\\=\cr              
\vbox{\smallskip\offinterlineskip
\halign{&\vsquare{##}\cr #1}}\,}

\def\borderlessrect#1#2{\hbox{\hskip \thickness
   \vbox to \Squaresize{\vskip \thickness \vss
      \hbox to #2 {\hss #1\hss}
   \vss\vskip\thickness} 
\unskip\hskip \thickness} 
\kern-\thickness}                                                            
                               
\def\vborderlessrect#1#2{\vbox{\borderlessrect{$#1$}{#2}}\kern-\thickness}

\def\borderless#1{\omit\vborderlessrect{#1}{\Squaresize}}
\def\borderlessrc#1#2{\omit\vborderlessrect{#1}{#2}}

\def\msquare#1{\vbox{\hbox{\vrule width \thickness
   \vbox to \Squaresize{\hrule height \thickness
      \hbox to \Squaresize{\hfil{\sevenrm #1}}
   \vfil\hrule height\thickness}
\unskip\vrule width \thickness}
\kern-\thickness}\kern-\thickness}

\def\twosquare#1#2{\vbox{\hbox{\vrule width \thickness 
   \vbox to \Squaresize{\hrule height \thickness
      \hbox to \Squaresize{\hfil{\rm #1}}\vss
      \hbox to \Squaresize{\hss{#2}\hss}
   \vfil\hrule height\thickness}
\unskip\vrule width \thickness}
\kern-\thickness}\kern-\thickness}

\def\twoblank#1#2{\vbox{\hbox{
   \vbox to \Squaresize{\vskip 2pt
      \hbox to \Squaresize{\hfil{\sevenrm #1}\ }\vss
      \hbox to \Squaresize{\hss{#2}\hss}
   \vfil}\unskip\kern-\thickness}
}\unskip\kern-\thickness}

\def\young#1{
\def\>{\blank}
\def\<{\borderless}
\def\*{\borderlessrc}
\def\p{\omit\msquare}
\def\t{\omit\twosquare}
\def\b{\omit\twoblank}
\let\\=\cr 
\vbox{\smallskip\offinterlineskip
\halign{&\vsquare{##}\cr #1}}}

\newdimen\smsquaresize \smsquaresize=12pt
\newdimen\smthickness \smthickness=.5pt
\font\smcellfont=cmss8 scaled \magstep0

\def\smsquare#1{\hbox{\vrule width \smthickness
   \unskip\vbox to \smsquaresize{\hrule height \smthickness\vss
      \hbox to \smsquaresize{\hss{\smcellfont #1}\hss}
   \vss\hrule height\smthickness} 
\unskip\vrule width \smthickness} 
\kern-\smthickness}

\def\smvsquare#1{\vbox{\smsquare{$#1$}}\kern-\smthickness}
\def\blank{\omit\hskip\smsquaresize}

\def\smyoung#1{\let\\=\cr 
\vbox{\smallskip\offinterlineskip
\halign{&\smvsquare{##}\cr #1}}}
\newdimen\vsmsquaresize \vsmsquaresize=10pt
\newdimen\vsmthickness \vsmthickness=.5pt
\font\vsmcellfont=cmsl8 scaled \magstep0
\font\vsmletterfont=cmr6 scaled \magstep0

\def\vsmsquare#1{\hbox{\vrule width \vsmthickness
   \unskip\vbox to \vsmsquaresize{\hrule height \vsmthickness\vss
      \hbox to \vsmsquaresize{\hss{\vsmcellfont #1}\hss}
   \vss\hrule height\vsmthickness} 
\unskip\vrule width \vsmthickness} 
\kern-\vsmthickness}
\def\vsmvsquare#1{\vbox{\vsmsquare{#1}}\kern-\vsmthickness}
\def\vsmblank{\omit\hskip\vsmsquaresize}
\def\vsmborderless#1{\hbox{\hskip \vsmthickness\unskip
   \vbox to \vsmsquaresize{\vss
      \hbox to \vsmsquaresize{\hss{\vsmletterfont #1}\hss}
   \vss} 
\unskip\hskip \vsmthickness} 
\kern-\vsmthickness}                                                            \def\vsmvborderless#1{\vbox{\vsmborderless{#1}}\kern-\vsmthickness}

\def\vsmyoung#1{
\def\>{\vsmblank}
\def\<{\omit\vsmvborderless}
\let\\=\cr 
\vbox{\smallskip\offinterlineskip
\halign{&\vsmvsquare{##}\cr #1}}}

\newdimen\edgesize \edgesize=20pt
\newdimen\eedgesize \eedgesize=21pt
\newdimen\doublesize \doublesize=41pt
\newdimen\ddoublesize \ddoublesize=43pt

\newdimen\triplesize \triplesize=62pt
\newdimen\tetrasize \tetrasize=85pt
\newdimen\futosa \futosa=1pt         

\def\Hako#1{\hbox{\vrule width \futosa
   \vbox to \eedgesize{\hrule height \futosa\vss                            
      \hbox to \edgesize{\hss#1\hss}
   \vss\hrule height\futosa} 
\unskip\vrule width \futosa} 
\kern-\futosa}                                                            
                               
\def\Seihokei#1{\vbox{\Hako{#1}}\kern-\futosa}

\def\Horizontal#1{\hbox{\vrule width \futosa
   \vbox to \eedgesize{\hrule height \futosa\vss                            
      \hbox to \doublesize{\hss#1\hss}
   \vss\hrule height \futosa} 
\unskip\vrule width \futosa} 
\kern-\futosa}                                                            
                               
\def\Vertical#1{\hbox{\vrule width \futosa
   \vbox to \ddoublesize{\hrule height \futosa\vss                            
      \hbox to \edgesize{\hss#1\hss}
   \vss\hrule height \futosa} 
\unskip\vrule width \futosa} 
\kern-\futosa}                                                            
                               
\def\NS#1{\vbox{\Hako{$#1$}\vskip22pt}\kern-\futosa}
\def\NSS#1{\vbox{\Hako{$#1$}\vskip42pt}\kern-\futosa}
\def\NSSS#1{\vbox{\Hako{$#1$}\vskip64pt}\kern-\futosa}
\def\H#1{\vbox{\Horizontal{$#1$}}\kern-\futosa}
\def\HH#1#2{\vbox{\Horizontal{$#1$}\Horizontal{$#2$}}\kern-\futosa}
\def\VER#1{\vbox{\Vertical{$#1$}}\kern-\futosa}
\def\VS#1{\vbox{\Vertical{$#1$}\vskip20pt}\kern-\futosa}
\def\VSS#1{\vbox{\Vertical{$#1$}\vskip42pt}\kern-\futosa}
\def\VV#1#2{\vbox{\Vertical{$#1$}\Vertical{$#2$}}\kern-\futosa}
\def\NV#1#2{\vbox{\Hako{$#1$}\Vertical{$#2$}}\kern-\futosa}
\def\NVS#1#2{\vbox{\Hako{$#1$}\Vertical{$#2$}\vskip22pt}\kern-\futosa}
\def\YTT#1#2#3{\vbox{\Horizontal{$#1$}\hbox{\vbox{\Vertical{$#2$}}\kern-\futosa\vbox{\Vertical{$#3$}}\kern-\futosa}}\kern-\futosa}

\def\domino#1{
\def\ns{\omit\NS}
\def\nss{\omit\NSS}
\def\nsss{\omit\NSSS}
\def\h{\omit\H}
\def\hh{\omit\HH}
\def\V{\omit\VER}
\def\Vs{\omit\VS}
\def\Vss{\omit\VSS}
\def\vsss{\omit\VSSS}
\def\vv{\omit\VV}
\def\nv{\omit\NV}
\def\nvs{\omit\NVS}
\def\ytt{\omit\YTT}
\let\\=\cr 
\vbox{\smallskip\offinterlineskip
\halign{&\Seihokei{##}\cr #1}}}

\theoremstyle{definition}
\newtheorem{theorem}{Theorem}[section]
\newtheorem{prop}[theorem]{Proposition}
\newtheorem{lemma}[theorem]{Lemma}
\newtheorem{corollary}[theorem]{Corollary}
\newtheorem{definition}[theorem]{Definition}

\newtheorem{remark}[theorem]{Remark}
\newtheorem{conjecture}[theorem]{Conjecture}

\newenvironment{demo}[1]{%
  \trivlist
  \item[\hskip\labelsep
        {\bf #1.}]
}{%
\hfill\qedsymbol
  \endtrivlist
}

\renewcommand{\mathcal}{\mathrsfs}

\newcommand{\p}{\partial}

\renewcommand{\b}{\beta}

\newcommand\Real{\mathbb{R}}

\newcommand\Pos{\mathbb{P}}

\newcommand\Pf{\operatorname{Pf}}

\newcommand\sgn{\operatorname{sgn}}

\newcommand\rdots{\mathinner{\mkern1mu\raise0pt\vbox{\kern7pt\hbox{.}}
     \mkern2mu\raise4pt\hbox{.}\mkern2mu\raise8pt\hbox{.}\mkern1mu}}
\def\TSSCPP#1{\mathcal{T}_{#1}}
\def\TSPP#1{\mathcal{B}_{#1}}
\def\CSPP#1{\mathcal{P}_{#1}}
\def\RCSPP#1{{\mathcal{P}^\text{R}_{#1}}}
\def\CCSPP#1{{\mathcal{P}^\text{C}_{#1}}}
\def\MT#1{\mathcal{M}_{#1}}
\def\BOX#1{{X}_{#1}}
\def\BOXC#1#2{{X}_{#1}({#2})}
\def\BOXU#1#2{{X}_{#1}^{#2}}
\def\BOXCU#1#2#3{{X}_{#1}^{#3}({#2})}
\def\BIJ#1{{\varphi}_{#1}}
\def\CMP#1{{\sigma}_{#1}}

\def\BORDER#1#2#3{{\theta}_{{#1}}{\left({#3}_{\,\geq{#2}}\right)}}
\def\PATH#1#2{{\cal P}\left({#1},{#2}\right)}
\def\NPATH#1#2{{\cal P}_0\left({#1},{#2}\right)}
\def\GF#1{{\operatorname{GF}}\left[{#1}\right]}
\def\SHAPE#1{{{\operatorname{sh}}\!\left({#1}\right)}}
\def\PROFILE#1{{\operatorname{pr}\left({#1}\right)}}

\def\REM#1#2{\operatorname{rem}\!\left({#1},{#2}\right)}
\def\CT#1{{\operatorname{CT}_{#1}}}

\def\V#1#2#3#4{V^{{#1},{#2}}\left({#3};{#4}\right)}

\numberwithin{equation}{section}
\def\rdots{\mathinner{\mkern1mu\raise0pt\vbox{\kern7pt\hbox{.}}
     \mkern2mu\raise4pt\hbox{.}\mkern2mu\raise8pt\hbox{.}\mkern1mu}}

\def\covered{\mathinner{\mkern1mu\raise0pt\vbox{\kern7pt\hbox{$<$}}
     \mkern-4mu\raise2pt\hbox{.}\mkern2mu}}
\def\covers{\mathinner{\mkern1mu\raise0pt\vbox{\kern7pt\hbox{$>$}}
     \mkern-12mu\raise2pt\hbox{.}\mkern8mu}}

\def\defterm#1{{\sl #1}\/}
\def\newterm#1{{\sl #1}\/}

\def\operatorname#1{{\mathrm{#1}\>\!}}

\def\Bbb#1{{\mathbb{#1}}}

\title{
On refined enumerations of
totally symmetric
self-complementary
plane partitions I
}

\author{
Masao Ishikawa\\
\small Faculty of Education, Tottori University\\[-0.8ex]
\small Koyama, Tottori, Japan\\[-0.8ex]
\small \texttt{ishikawa@fed.tottori-u.ac.jp}
}

\date{
\medbreak
}

\begin{document}

\maketitle

\begin{abstract}

In this paper we give Pfaffian expressions and constant term identities
for three conjectures
(i.e. Conjecture~2, Conjecture~3 and Conjecture~7)
by Mills, Robbins and Rumsey in the paper
``Self-complementary totally symmetric plane partitions''
{\em J. Combin. Theory Ser. A} {\bf  42}, 277--292)
concerning the refined enumeration problems of totally symmetric
self-complementary
plane partitions.
We also present some new conjectures and give Pfaffian expressions and constant term identities
for them.
But evaluation problem of these Pfaffians are still difficult.

\end{abstract}

{
\small
\begin{spacing}{0.1}
\tableofcontents 
\end{spacing}
}

\section{
Introduction
}\label{sec:intro}


In the paper \cite{MRR2}
Mills, Robbins and Rumsey presented several conjectures on the enumeration
of the totally symmetric self-complementary plane partitions.
G.E.~Andrews (\cite{A2}) settled the conjecture (\cite[Conjecture~1]{MRR2})
on the cardinality of the totally symmetric self-complementary plane partitions of size $n$
(see also \cite{Ste1}).
D.~Zeilberger gave a constant term identity of this cardinality in \cite{Z1}.
The aim of this paper is to give Pfaffian expressions
of the other conjectures in \cite{MRR2}
which are the enumeration with some weight or enumeration of some subset.
We also generalize Zeilberger's constant term identity,
and show that each enumeration correspond to
each of the classical Littlwood type identities for Schur functions.

In \cite{MRR2}
Mills, Robbins and Rumsey have introduced a class $\TSPP{n}$ of 
triangular shifted plane partitions
\begin{equation*}
\begin{array}{cccc}
b_{11}&b_{12}&\hdots&b_{1,n-1}\\
      &b_{22}&\hdots&b_{2,n-1}\\
      &      &\ddots&\vdots\\
      &      &      &b_{n-1,n-1}\\
\end{array}
\end{equation*}
whose parts are $\leq n$,
weakly decreasing along rows and columns,
and all parts in row $i$ are $\geq n-i$.
For example,
$\TSPP{3}$ consists of the following seven elements.
\begin{equation*}
\begin{array}{cc}
3 & 3\\
  & 3
\end{array}
\quad
\begin{array}{cc}
3 & 3\\
  & 2
\end{array}
\quad
\begin{array}{cc}
3 & 3\\
  & 1
\end{array}
\quad
\begin{array}{cc}
3 & 2\\
  & 2
\end{array}
\quad
\begin{array}{cc}
3 & 2\\
  & 1
\end{array}
\quad
\begin{array}{cc}
2 & 2\\
  & 2
\end{array}
\quad
\begin{array}{cc}
2 & 2\\
  & 1
\end{array}
\end{equation*}
They have established an bijection between the totally symmetric
self-complementary plane partitions of size $n$
and the elements of $\TSPP{n}$
(see Section~\ref{sec:bijections}),
and defined a new statistics $U_{r}(b)$ for
an element $b\in\TSPP{n}$ and $r=1,\dots,n$,
i.e.,
for a $b=(b_{ij})_{1\leq i\leq j\leq n-1}$ in $\TSPP{n}$,
let
\begin{equation}
U_{r}(b)=\sum_{t=1}^{n-r}(b_{t,t+r-1}-b_{t,t+r})
+\sum_{t=n-r+1}^{n-1}\chi\{b_{t,n-1}>n-t\}.
\label{eq:stat_mrr}
\end{equation}
Here $\chi\{\dots\}$ has value $1$ when the statement ``$\dots$'' is true and $0$ otherwise,
and we use the convention that $b_{i,n}=n-i$ for all $i$ and $b_{0,j}=n$ for all $j$.

Mills, Robbins and Rumsey conjectured that $U_{r}$ has the same distribution as the position of the $1$ 
in the top row of an alternating sign matrix,
and presented several conjectures related to the distribution of the statistics $U_{r}$.
The aim of this paper is to obtain the generating functions
for the enumerations concerning these conjectures.
In this introduction,
we briefly review these conjectures by Mills, Robbins and Rumsey,
and present a Pfaffian expression for each problem.
In fact,
we generalize the definition of $\TSPP{n}$ to $\TSPP{n,m}$ in Section~\ref{sec:bijections},
and consider the generating functions in wider classes of plane partitions,
so that these results are special cases of the theorems obtained in the following sections.
For the definition of the numbers $A_{n}$, $A_{n}^{k}$, $A_{n}^{k,l}$, $A^\text{VS}_{n,r}$, $A^{\text{VS},\,r}_{n}$
and the polynomials $A_{n}(t)$, $A_{n}(t,u)$, $A^\text{VS}_{2n+1}(t)$,
the reader should refer to the Section~\ref{sec:preliminaries}.
It seems that these numbers have the standard notation which have appeared concerning
the alternating sign matrices
(see \cite{Ku3,O2,RS1,Str1}).
Let
$
\bar S_{n}=(\bar s_{ij})_{1\leq i,j\leq n}
$
be the skew-symmetric matrix of size $n$ whose $(i,j)$th entry $\bar s_{ij}$ is equal to $(-1)^{j-i-1}$
for $1\leq i<j\leq n$,
and let $O_n$ denote the $n\times n$ zero matrix.
Let $J_n=(\delta_{i,n+1-j})_{1\leq i,j\leq n}$ denote the anti-diagonal matrix
where $\delta_{i,j}$ stands for the Kronecker delta function. 
First of all,
Mills, Robbins and Rumsey presented the following conjecture in the paper \cite{MRR2},
which we call the refined enumeration of TSSCPPs:
\begin{conjecture}
\label{conj:refined}
(\cite[pp.282, Conjecture~2]{MRR2})
Let $n$ be a positive integer.
Let $1\leq k\leq n$ and $1\leq r\leq n$.
Then the number of elements $b$ of $\TSPP{n}$ such that $U_{r}(b)=k-1$ would be $A_{n}^{k}$.
Namely,
 $\sum_{b\in\TSPP{n}}t^{U_{r}(b)}=A_{n}(t)$ would hold.
\end{conjecture}
Let $n$ and $N$ be positive integers,
and
let $B_{n}^{N}(t)=(b_{ij}(t))_{0\leq i\leq n-1,\ 0\leq j\leq n+N-1}$
be the $n\times (n+N)$ matrix whose $(i,j)$th entry is 
\begin{equation*}
b_{ij}(t)=
\begin{cases}
\delta_{0,j}
&\text{ if $i=0$,}\\
\binom{i-1}{j-i}+\binom{i-1}{j-i-1}t
&\text{ otherwise.}
\end{cases}
\end{equation*}
Especially, when $t=1$,
we write $B_{n}^{N}$ for $B_{n}^{N}(1)$
whose $(i,j)$th entry is $\binom{i}{j-i}$.
One of the results we obtain for Conjecture~\ref{conj:refined} is following:
\begin{theorem}
\label{thm:refined}
Let $n$ be a positive integer
and let $N$ be an even integer such that $N\geq n-1$.
Then
\begin{equation}
\sum_{b\in\TSPP{n}}t^{U_{r}(b)}
=\Pf\begin{pmatrix}
O_{n}&J_{n}B_{n}^{N}(t)\\
-{}^t\!B_{n}^{N}(t)J_{n}&\bar S_{n+N}\\
\end{pmatrix}.
\end{equation}
(cf. Corollary~\ref{cor:bijection_TSPPtoCSPP}, Theorem~\ref{th:U_k}, Corollary~\ref{cor:refined} and Corollary~\ref{cor:const_refined}).
\end{theorem}
For example,
if $n=3$ and $N=2$ then the above Pfaffian looks like as follows.
\[
\Pf\left(
  \begin{array}{ccc|ccccc}
   0 & 0 & 0 & 0 & 0 & 1 &1+t& t \\
   0 & 0 & 0 & 0 & 1 & t & 0 & 0 \\
   0 & 0 & 0 & 1 & 0 & 0 & 0 & 0 \\\hline
   0 & 0 &-1 & 0 & 1 &-1 & 1 &-1 \\
   0 &-1 & 0 &-1& 0 & 1 &-1 & 1 \\
  -1 &-t & 0 & 1 &-1 & 0 & 1 &-1 \\
-1-t & 0 & 0 &-1 & 1 &-1 & 0 & 1 \\
  -t & 0 & 0 & 1 &-1 & 1 &-1 & 0 \\
  \end{array}
\right).
\]
In the same paper,
they also presented the following conjecture
which we call the doubly refined enumeration of TSSCPPs:

\begin{conjecture}
\label{conj:double_refined}
(\cite[pp.284, Conjecture~3]{MRR2}, \cite{Str1})
Let $n\geq2$ and $1\leq k,l\leq n$ be integers.
Then the number of elements $b$ of $\TSPP{n}$ such that $U_{1}(b)=k-1$ 
and $U_{2}(b)=n-l$ would be $A_{n}^{k,l}$. 
\end{conjecture}
Let $n$ and $N$ be positive integers.
Let $B_{n}^{N}(t,u)=(b_{ij}(t,u))_{0\leq i\leq n-1,\ 0\leq j\leq n+N-1}$
be the $n\times (n+N)$ matrix whose $(i,j)$th entry is
\begin{equation*}
b_{ij}(t,u)=
\begin{cases}
\delta_{0,j}
&\text{ if $i=0$,}\\
\delta_{0,j-i}+\delta_{0,j-i-1}tu
&\text{ if $i=1$,}\\
\binom{i-2}{j-i}+\binom{i-2}{j-i-1}(t+u)+\binom{i-2}{j-i-2}tu
&\text{ otherwise.}
\end{cases}
\end{equation*}
Note that,
when $u=1$,
$B_{n}^{N}(t,1)$ is equal to $B_{n}^{N}(t)$.
Then one form of the Pfaffian expressions for Conjecture~\ref{conj:double_refined}
which we obtain in this paper
is following:
\begin{theorem}
\label{thm:double_refined}
Let $n$ be a positive integer
and let $N$ be an even integer such that $N\geq n-1$.
If $r$ is an integer such that $2\leq r\leq n$,
then we have
\begin{equation}
\sum_{b\in\TSPP{n}}t^{U_{1}(b)}u^{U_{r}(b)}
=\Pf\begin{pmatrix}
O_{n}&J_{n}B_{n}^{N}(t,u)\\
-{}^t\!B_{n}^{N}(t,u)J_{n}&\bar S_{n+N}\\
\end{pmatrix}.
\end{equation}
(cf. Corollary~\ref{cor:bijection_TSPPtoCSPP}, Theorem~\ref{th:U_k}, Corollary~\ref{cor:doubly_refined} and Corollary~\ref{cor:const_doubly_refined}).
\end{theorem}

The monotone triangles are known to be in one-to-one correspondence
with the alternating sign matrices
(\cite{B,MRR1}).
Here we arrange our definition
following the notation in \cite{MRR2}.
A \defterm{monotone triangle} of size $n$ is,
by definition,
 a triangular array of positive integers
\begin{equation*}
\begin{array}{cccc}
&&&m_{n,n}\\
&&m_{n-1,n-1}&m_{n-1,n}\\
&\rdots&\vdots&\vdots\\
m_{1,1}&\hdots&m_{1,n-1}&m_{1,n}\\
\end{array}
\end{equation*}
subject to the constraints that
\begin{enumerate}
\item[(M1)]
$m_{ij}<m_{i,j+1}$ whenever both sides are defined,
\item[(M2)]
$m_{ij}\geq m_{i+1,j}$ whenever both sides are defined,
\item[(M3)]
$m_{ij}\leq m_{i+1,j+1}$ whenever both sides are defined,
\item[(M4)]
the bottom row $(m_{1,1},m_{1,2},\dots,m_{1,n})$ is $(1,2,\dots, n)$.
\end{enumerate}
Let $\MT{n}$ denote the set of monotone triangles of size $n$.
For example,
$\MT{3}$ consists of the following seven elements.
\begin{equation*}
\begin{array}{ccc}
 & &1\\
 &1&2\\
1&2&3
\end{array}
\quad
\begin{array}{ccc}
 & &2\\
 &1&2\\
1&2&3
\end{array}
\quad
\begin{array}{ccc}
 & &1\\
 &1&3\\
1&2&3
\end{array}
\quad
\begin{array}{ccc}
 & &2\\
 &1&3\\
1&2&3
\end{array}
\quad
\begin{array}{ccc}
 & &3\\
 &1&3\\
1&2&3
\end{array}
\quad
\begin{array}{ccc}
 & &2\\
 &2&3\\
1&2&3
\end{array}
\quad
\begin{array}{ccc}
 & &3\\
 &2&3\\
1&2&3
\end{array}
\end{equation*}
Note that,
if one removes the bottom row of $m\in\MT{n}$
and turn it upside-down,
then he get an array defined in \cite{MRR2}.

For $k=0,1,\dots,n-1$,
let $\MT{n}^{k}$ denote the set of monotone triangles
with all entries $m_{ij}$ in the first $n-k$ columns equal
to their minimum values $j-i+1$.
For $k=0,1,\dots,n-1$,
let $\TSPP{n}^{k}$ be the subset of those $b$ in $\TSPP{n}$ 
such that all $b_{ij}$ in the first $n-1-k$ columns are equal to their
maximal values $n$.
Then they also presented the following conjecture:
\begin{conjecture}
\label{conj:MT}
(\cite[pp.287, Conjecture~7]{MRR2})
For $n\geq2$ and $k=0,1,\dots,n-1$,
the cardinality of $\TSPP{n}^{k}$ is equal to 
the cardinality of $\MT{n}^{k}$.
\end{conjecture}
Let $m$, $n$ and $k$ be integers such that
$1\leq m\leq n$ and $0\leq k\leq n-m$.
We define the $n\times n$ skew-symmetric matrix
$\bar L_{n}^{(m,k)}(\varepsilon)=(\bar l^{(m,k)}_{ij}(\varepsilon))_{1\leq i,j\leq n}$
as follows:
if $k$ is even, then
\[
\bar l^{(m,k)}_{ij}(\varepsilon)
=\begin{cases}
(-1)^{j-i-1}\varepsilon
&\text{ if $1\leq i<j\leq n$ and $i\leq m+k$,}\\
(-1)^{j-i-1}
&\text{ if $m+k< i<j\leq n$,}
\end{cases}
\]
else
\[
\bar l^{(m,k)}_{ij}(\varepsilon)
=\begin{cases}
(-1)^{j-i-1}\varepsilon
&\text{ if $1\leq i<j\leq m+k$,}\\
(-1)^{j-i-1}
&\text{ if $1\leq i<j\leq n$ and $m+k<j$.}
\end{cases}
\]
For example,
\[
\bar L_{6}^{(2,1)}(\varepsilon)
=\begin{pmatrix}
           0& \varepsilon&-\varepsilon&           1&          -1&           1\\
-\varepsilon&           0& \varepsilon&          -1&           1&          -1\\
 \varepsilon&-\varepsilon&           0&           1&          -1&           1\\
          -1&           1&          -1&           0&           1&          -1\\
           1&          -1&           1&          -1&           0&           1\\
          -1&           1&          -1&           1&          -1&           0
\end{pmatrix},
\quad
\bar L_{6}^{(2,2)}(\varepsilon)
=\begin{pmatrix}
           0& \varepsilon&-\varepsilon& \varepsilon&-\varepsilon& \varepsilon\\
-\varepsilon&           0& \varepsilon&-\varepsilon& \varepsilon&-\varepsilon\\
 \varepsilon&-\varepsilon&           0& \varepsilon&-\varepsilon& \varepsilon\\
-\varepsilon& \varepsilon&-\varepsilon&           0& \varepsilon&-\varepsilon\\
 \varepsilon&-\varepsilon& \varepsilon&-\varepsilon&           0&           1\\
-\varepsilon& \varepsilon&-\varepsilon& \varepsilon&          -1&           0
\end{pmatrix}.
\]
Then a Pfaffian expression for Conjecture~\ref{conj:MT}
which we obtain in this paper
is following:
\begin{theorem}
\label{thm:MT}
Let $n$ be a positive integer and let $k=0,1,\dots,n-1$.
Let $N$ be an even integer such that $N\geq k$.
The cardinality of $\TSPP{n}^{k}$ is equal to 
\begin{equation}
\lim\limits_{\varepsilon\to0}\,
\varepsilon^{-\lfloor\frac{k}{2}\rfloor}
\Pf\begin{pmatrix}
O_{n}&B_{n}^{N}J_{n+N}\\
-J_{n+N}{}^t\!B_{n}^{N}&\bar L_{n+N}^{(n,k)}(\varepsilon)\\
\end{pmatrix}
\end{equation}
Here $\lfloor x\rfloor$ stands for the floor function,
i.e. the greatest integer less than or equal to $x$. 
(cf. Theorem~\ref{prop:TSPPtoCSPP_k}, Conjecture~\ref{conj:refined_MT},
Corollary~\ref{cor:MT} and Corollary~\ref{cor:const_mt}).
\end{theorem}

This paper is composed as follows.
All through the paper, we consider more general set
$\TSPP{n,m}$
(see Definition~\ref{def:TSPP})
of shifted plane partitions,
which first appeared in \cite[Theorem~1]{Kr1}.
Thus the above theorems,
which give Pfaffian expressions for the Mills, Robbins and Rumsey conjectures,
are the special cases of the generating functions we obtain in this paper.
The key idea is to construct a bijection
between $\TSPP{n,m}$ and a new set $\CSPP{n,m}$
which is more easy to understand.
In Section~\ref{sec:bijections},
we define this set $\CSPP{n,m}$ (see Definition~\ref{def:CSPP}) of ordinary plane partitions
which is the main object we study throughout the paper.
We establish a bijection between $\CSPP{n,m}$ and a set $\TSSCPP{n,m}$ 
of totally symmetric self-complementary plane partitions
(see Theorem~\ref{thm:bijection_CSPP}),
and also construct a bijection between $\TSPP{n,m}$ and $\TSSCPP{n,m}$
(see Theorem~\ref{thm:bijection_TSPP}).
As a corollary we obtain a bijection 
 between $\CSPP{n,m}$ and $\TSPP{n,m}$
(see Corollary~\ref{cor:bijection_TSPPtoCSPP})
which makes it possible to interpret all the properties of $\TSPP{n}$ studied in \cite{MRR1}
in the words of $\CSPP{n}$
(here we write $\TSPP{n}$ for $\TSPP{n,0}$ and $\CSPP{n}$ for $\CSPP{n,0}$).
Thus we reduce the enumeration problems of the totally symmetric self-complementary
plane partitions to the study of our new object $\CSPP{n,m}$,
which we call ``the restricted column-strict plane partitions''.
This reveals several mysterious properties of this new object $\CSPP{n,m}$
which resembles the classical theory of the tableaux and Schur functions.
Before we proceed to Section~\ref{sec:bijections},
we collect some basic
definitions and several fundamental theorems
in Section~\ref{sec:preliminaries}.
Especially the minor summation formula and its applications
(see Theorem~\ref{msf} and Proposition~\ref{prop:skew})
will be an important tool to obtain the generating functions in Section~\ref{sec:gf}.
But the reader can skip this section now, 
and use it as a reference when he need it.
In Section~\ref{sec:statistics},
we generalize the statistics \thetag{\ref{eq:stat_mrr}}
 to the general set $\TSPP{n,m}$
(see \thetag{\ref{eq:stat}})
and restate it as the statistics of $\CSPP{n,m}$
(see Theorem~\ref{th:U_k}).
We also present a new conjecture (Conjecture~\ref{conj:even_row})
which is not in \cite{MRR1}.
In Section~\ref{sec:mt},
we restate Conjecture~\ref{conj:MT} in the words of $\CSPP{n}$
(Theorem~\ref{prop:TSPPtoCSPP_k})
and also present a new conjecture
(Conjecture~\ref{conj:refined_MT})
which is a refined version of Conjecture~\ref{conj:MT}.
We also give a restatement of \cite[Conjecture~7']{MRR1}
in the words of $\CSPP{n}$
(Theorem~\ref{prop:TSPPtoCSPP_kr}).
In Section~\ref{sec:se},
we translate the strange enumeration (especially $(-1)$-enumeration)
of the totally symmetric self-complementary plane partitions
in the words of $\CSPP{n,m}$
(see Theorem~\ref{thm:strange}).
Then Section~\ref{sec:gf} contains the main results of this paper,
i.e. we obtain several generating functions concerning $\CSPP{n,m}$
using the lattice paths.
We give
Corollary~\ref{cor:refined} for the refined TSSCPP enumeration,
Corollary~\ref{cor:doubly_refined} for the doubly refined TSSCPP enumeration,
and
Corollary~\ref{cor:MT} for Conjecture~\ref{conj:MT}.
In Section~\ref{sec:cti}
we give a constant term identity for each Pfaffian obtained in Section~\ref{sec:gf}
(see Corollary~\ref{cor:const_doubly_refined}, Corollary~\ref{cor:const_refined} and Corollary~\ref{cor:const_mt}).
Each of the constant identities includes a different Littlewood type identity.
This reveals a certain relation between the Littlewood type identity
for the enumerations of the Schur functions and the enumerations of the TSSCPPs.

\section{
Preliminaries
}\label{sec:preliminaries}

First we recall the numbers 
and polynomials
related to the alternating sign matrices
(cf. \cite{B,Ku3,MRR1,O2,Ro,RS1,Str1,Z2}).
In the latter half of this section,
we recall the notation of partitions
and the results on Pfaffians
which will be needed in the following sections,
i.e.
Theorem~\ref{msf} and Proposition~\ref{prop:skew}.
For the details on partitions 
the reader is referred to \cite{M,Sta2},
and for the explanation on Pfaffians
the reader can consult \cite{IW1,IW2,Kn}.
Let $A_{n}$ denote the number defined by
\begin{equation}
A_{n}=\prod_{i=0}^{n-1}\frac{(3i+1)!}{(n+i)!}.
\end{equation}
This number is famous for the alternating sign matrix conjecture
(cf. \cite{B}). 
The number of totally symmetric self-complementary plane partitions was conjectured
to be $A_n$ in \cite[pp.282, Conjecture~1]{MRR2},
and settled in \cite[p.p.127, Theorem~8.3]{Ste1} and \cite{A2}
(see also \cite{A1,AB}),
Another proof was appeared in \cite{Kr1} and several determinant techniques have been developed in it.
Let $n$ be a positive number and let  $1\leq r\leq n$.
Set $A_{n}^{r}$ to be the number
\begin{equation}
A_{n}^{r}=\frac{\binom{n+r-2}{n-1}\binom{2n-r-1}{n-1}}{\binom{2n-2}{n-1}}A_{n-1}
=\frac{\binom{n+r-2}{n-1}\binom{2n-1-r}{n-1}}{\binom{3n-2}{n-1}}A_{n}
.
\end{equation}
Then the number $A_{n}^{r}$ satisfies the recurrence $A_{n}^{1}=A_{n-1}$ and
\[
\frac{A_{n}^{r+1}}{A_{n}^{r}}
=\frac{(n-r)(n+r-1)}{r(2n-r-1)}.
\]
The number has appeared to describe the distribution of the position of the $1$ in the top row of an alternating sign matrix
(see \cite{Ku1,Ku3,O2,Z2}).
We also define the polynomial $A_n(t)=\sum_{r=1}^{n}A_{n}^{r}t^{r-1}$.
For instance, 
the first few terms are
$A_1(t)=1$, $A_{2}(t)=1+t$, $A_{3}(t)=2+3t+2t^2$, $A_{4}(t)=7+14t+14t^2+7t^3$.
Let $n$ be a positive integer and 
let $A_{n}^{k,l}$, $1\leq k,l\leq n$,
 denote the number which satisfies the initial condition
\begin{equation}
A_{n}^{k,1}=A_{n}^{1,k}=\begin{cases}
0&\text{ if $k=1$}\\
A_{n-1}^{k-1}&\text{ if $2\leq k\leq n$}
\end{cases}
\end{equation}
and the recurrence equation
\begin{equation}
A_{n}^{k+1,l+1}-A_{n}^{k,l}
=\frac{A_{n-1}^{k}(A_{n}^{l+1}-A_{n}^{l})+A_{n-1}^{l}(A_{n}^{k+1}-A_{n}^{k})}{A_{n}^{1}}
\end{equation}
for $1\leq k,l\leq n-1$.
This recurrence equation satisfied by $A_{n}^{k,l}$ has been introduced by Stroganov in \cite[Section~5]{Str1}
to describe the double distribution of the positions of the $1$'s 
in the top row and the bottom row of an alternating sign matrix.
For example, if $n=3,4$,
then we have
\begin{equation*}
\left(A_{3}^{k,l}\right)_{1\leq k,l\leq 3}=
\begin{pmatrix}
0&1&1\\
1&1&1\\
1&1&0\\
\end{pmatrix},
\qquad
\left(A_{4}^{k,l}\right)_{1\leq k,l\leq 4}=
\begin{pmatrix}
0&2&3&2\\
2&4&5&3\\
3&5&4&2\\
2&3&2&0\\
\end{pmatrix}.
\end{equation*}
Let $A_{n}(t,u)$ denote the polynomial defined by
$A_{n}(t,u)=\sum_{k,l=1}^{n}A_{n}^{k,l}t^{k-1}u^{n-l}$.
Thus we have $A_{3}(t,u)=1+t+u+tu+t^2u+tu^2+t^2u^2$.
Let $\omega=e^{2i\pi/3}$.
Di~Francesco and Zinn-Justin showed that $A_{n}(t,u)$
can be expressed by the Schur function as
\begin{equation}
A_{n}(t,u)
=\frac{\{\omega^2(\omega+t)(\omega+u)\}^{n-1}}{3^{n(n-1)/2}}
s_{\delta(n-1,n-1)}^{(2n)}\left(
\frac{1+\omega t}{\omega+t},\frac{1+\omega u}{\omega+u},1,\dots,1
\right)
\end{equation}
where $s_{\lambda}^{(n)}(x_1,\dots,x_{n})$ stands for the Schur function
in the $n$ variables $x_1$, $\dots$, $x_n$, corresponding to the partition $\lambda$,
and $\delta(n-1,n-1)=(n-1,n-1,n-2,n-2,\dots,1,1)$
(See \cite[pp.4]{FZ1}, \cite{O2}).

Let $A^\text{VS}_{2n+1}$ be the number defined by
\begin{equation}
A^\text{VS}_{2n+1}
=(-3)^{n^2}\prod_{{1\leq i,j\leq 2n+1}\atop{2|j}}\frac{3(j-i)+1}{j-i+2n+1}
=\frac1{2^n}\prod_{k=1}^{n}\frac{(6k-2)!(2k-1)!}{(4k-1)!(4k-2)!}
\end{equation}
and let ${A^\text{VS, $r$}_{2n+1}}$ be the number given by
\begin{equation}
A^\text{VS, $r$}_{2n+1}=\frac{A^\text{VS}_{2n-1}}{(4n-2)!}
\sum_{k=1}^{r}(-1)^{r+k}
\frac{(2n+k-2)!(4n-k-1)!}{(k-1)!(2n-k)!}.
\end{equation}
This number $A^\text{VS}_{2n+1}$ is equal to the number of vertically symmetric alternating sign matrices of size $2n+1$
(see \cite{Ku3,O2,RS1}).
For example,
the first few terms of $A^\text{VS}_{2n+1}$
are $1$, $3$, $26$, $646$ and $45885$.
We also define the polynomial $A^\text{VS}_{2n+1}(t)$ by
\begin{equation}
A^\text{VS}_{2n+1}(t)=\sum_{r=1}^{2n}A^\text{VS, $r$}_{2n+1}t^{r-1}.
\label{eq:poly_VS}
\end{equation}
For instance,
the first few terms of \thetag{\ref{eq:poly_VS}} are
$A^\text{VS}_{3}(t)=1$,
$A^\text{VS}_{5}(t)=1+t+t^2$,
$A^\text{VS}_{7}(t)=3+6{t}+8{t}^{2}+6{t}^{3}+3{t}^{4}$
and
$A^\text{VS}_{9}(t)=26+78{t}+138{t}^{2}+162{t}^{3}+138{t}^{4}+78{t}^{5}+26{t}^{6}$.

Next we recall and fix the notation of partitions and shifted partitions.
We follow the notation and terminology of Macdonald \cite{M}.
Let $\Pos$ denote the set of positive integers.
A partition is a sequence
$\lambda=(\lambda_1,\lambda_2,\dots)$
of non-negative integers in non-increasing order:
$\lambda_1\geq\lambda_2\geq\dots$
and containing only finitely many non-zero terms.
The non-zero $\lambda_i$ are call the \defterm{parts} of $\lambda$.
The number of parts is the \defterm{length} of $\lambda$,
denoted by $\ell(\lambda)$;
and the sum of parts is the \defterm{weight} of $\lambda$,
denoted by $|\lambda|$.
The \defterm{diagram} of a partition $\lambda$ may be formally defined as the set of lattice points $(i,j)\in\Bbb{P}^2$
such that $1\leq j\leq\lambda_i$.
We identify $\lambda$ with its diagram.
The number of nodes in
the main diagonal of the diagram is called the \defterm{diagonal length} $d(\lambda)=\sharp\{i:\lambda_i\geq i\}$.
The \defterm{conjugate} of a partition $\lambda$ is the partition $\lambda'$
whose diagram is the transpose of the diagram of $\lambda$.
A \defterm{self-conjugate partition} is a partition whose conjugate partition is equal to itself.
For a partition $\lambda$,
let $r(\lambda)$ denote
the number of rows of odd length,
and
let $c(\lambda)$ denote
the number of columns of odd length.

A partition with distinct parts is called a \defterm{strict partition}.
The \defterm{shifted diagram} of a strict partition $\mu$ is the set of lattice points $(i,j)\in\Bbb{P}^2$
such that $i\leq j\leq\mu_i+i$.
We identify a strict partition with its shifted diagram.
If $\lambda$ is a self-conjugate partition,
we can associate a strict partition $(\mu_i)_{1\leq i\leq d(\lambda)}$ to $\lambda$
where $\mu_i=\lambda_{i}-i+1$,
and this defines a bijection of the set of self-conjugate partitions onto the set of strict partitions.
We say that a partition is \defterm{even} if all of its parts $\lambda_i$ are even.

A \newterm{$q$-binomial coefficient} is,
by definition,
$\left[{n\atop r}\right]_{q}=\frac{(q)_{n}}{(q)_{r}(q)_{n-r}},$ 
where $(q)_{k}=\prod_{i=1}^{k}(1-q^i)$.
A binomial coefficient is written as $\left({n\atop r}\right)=\left[{n\atop r}\right]_{1}$.
For any finite set $S$ and a non-negative integer $r$,
let $\binom{S}{r}$ denote the set of all $r$-element subsets of $S$.
For a subset $I=\{i_1,\dots,i_r \}\in\binom{[n]}{r}$,
let $\overline I$ denote the set-theoretic complement
of $I$ in $[n]$.
Let $m$, $n$ and $r$ be integers such that $r\leq m,n$ and let $T$ be an $m$ by $n$ matrix.
For any index sets $I=\{i_1,\dots,i_r \}\in\binom{[m]}{r}$ and $J=\{j_1,\dots,j_r\}\in\binom{[n]}{r}$,
let $A^{I}_{J}$ denote the submatrix obtained by selecting the rows indexed by $I$ and the columns indexed by $J$.
If $r=m$ and $I=[m]$,
we simply write $A_{J}$ for $A^{[m]}_{J}$.
Similarly,
if $r=n$ and $J=[n]$,
we write $A^{I}$ for $A_{[n]}^{I}$.
The following identity \thetag{\ref{eq:msf}} follows from the proof of \cite[Theorem~4.2]{IW2})
which we call the minor summation formula here
(see also \cite{Ha}).
The formula is equivalent to the even case \thetag{\ref{eq:msf_even}} 
and the odd case \thetag{\ref{eq:msf_odd}} which appeared in \cite{IW1},
but the merit of using \thetag{\ref{eq:msf}} is that it does not depend on whether $n$
is even or odd.
\begin{theorem}
\label{msf}
Let $m$ and $n$ be positive integers such that $m\le n$ and $n-m$ is even.
Let $T=(t_{ij})_{1\le i\le m, 1\le j\le n}$ be an $m$ by $n$ rectangular matrix.
Let $B$ be a skew-symmetric matrix of size $n$.
Then
\begin{align}
\sum_{I\in\binom{[n]}{m}}
(-1)^{s(\overline I,I)}
\Pf(B^{\overline I}_{\overline I}) \det(T_I)
&=\Pf\begin{pmatrix}
   O_m   &   T J_n \\
-J_n\, {}^t\!T &J_n\,{}^t\!B J_n \\
\end{pmatrix}
=\Pf\begin{pmatrix}
   O_m  & J_m T \\
-{}^tT J_m    & B \\
\end{pmatrix}.
\label{eq:msf}
\end{align}
Here $\overline I=[n]\setminus I$,
and $s(\overline I,I)$ denote the shuffle number
to merge $\overline I$ with $I$ into $[n]$.
\end{theorem}

From here we define several skew-symmetric matrices
which play an important role in the applications.
Let $n$ be a positive integer.
Let
$
S_{n}=(s_{ij})_{1\leq i,j\leq n}
$
be the skew-symmetric matrix of size $n$ whose $(i,j)$th entry $s_{ij}$ is $1$
for $1\leq i<j\leq n$,
and let $\bar S_{n}$ be as defined in Section~\ref{sec:intro}.
Let $\REM{a}{b}$ denote the remainder of $a$ divided by $b$,
and let $t$ be an indeterminate.
Let the $n\times n$ skew-symmetric matrices
$R_{n}(t)=(r_{ij}(t))_{1\leq i,j\leq n}$,
$C_{n}(t)=(c_{ij}(t))_{1\leq i,j\leq n}$,
$\bar R_{n}(t)=(\bar r_{ij}(t))_{1\leq i,j\leq n}$
and
$\bar C_{n}(t)=(\bar c_{ij}(t))_{1\leq i,j\leq n}$
be defined as
$r_{ij}(t)=t^{\REM{i-1}{2}+\REM{j}{2}}$,
$c_{ij}(t)=t^{j-i-1}$,
$\bar r_{ij}(t)=(-1)^{j-i-1}\,t^{j-i-1}$
and
$\bar c_{ij}(t)=(-1)^{j-i-1}\,t^{\REM{n+1-i}{2}+\REM{n-j}{2}}$
 for $1\leq i<j\leq n$,
respectively.
Further we write
$R_{n}$,
$C_{n}$,
$\bar R_{n}$
and
$\bar C_{n}$
 for 
$R_{n}(0)$,
$C_{n}(0)$,
$\bar R_{n}(0)$
and
$\bar C_{n}(0)$,
respectively.
For example,
\begin{equation*}
R_{4}(t)
=\begin{pmatrix}
  0 &  1  &   t &  1  \\
 -1 &  0  & t^2 & t \\
-t  &-t^2 &   0 &  1  \\
-1  & -t  &  -1 &  0 
\end{pmatrix}
\quad\text{ and }\quad
C_{4}(t)
=\begin{pmatrix}
  0 &  1  &   t & t^2 \\
 -1 &  0  &  1  &  t  \\
-t  & -1  &   0 &  1  \\
-t^2& -t  &  -1 &  0 
\end{pmatrix}.
\end{equation*}
Let $m$, $n$ and $k$ be integers such that
$1\leq m\leq n$ and $0\leq k\leq n-m$,
and let $\varepsilon$ be an indeterminate.
Let $\bar L_{n}^{(m,k)}(\varepsilon)$ be as in Section~\ref{sec:intro},
and
let $L_{n}^{(m,k)}(\varepsilon)=(l^{(m,k)}_{ij}(\varepsilon))_{1\leq i,j\leq n}$
denote the $n\times n$ skew-symmetric matrix
whose $(i,j)$th entry is
\[
l^{(m,k)}_{ij}(\varepsilon)
=\begin{cases}
1
&\text{ if $1\leq i<j\leq m+k$,}\\
\varepsilon
&\text{ if $1\leq i<j\leq n$ and $m+k<j$.}
\end{cases}
\]
The following lemma 
(cf. \cite[Section~4, Lemma~7]{IW1})
is very useful to compute the subpfaffians of a given skew-matrix,
and the basic idea to prove Proposition~\ref{prop:skew}.
\begin{lemma}
\label{lemma:product}
Let $x_i$ and $y_j$ be indeterminates,
and let $n$ is a non-negative integer.
Then
\begin{equation*}
\label{prod_pf}
\Pf[x_iy_j]_{1\leq i<j\leq 2n}=\prod_{i=1}^{n}x_{2i-1}\prod_{i=1}^{n}y_{2i}.
\ \Box
\end{equation*}
\end{lemma}
Let $\lambda=(\lambda_1,\dots,\lambda_m)$ be a partition such that $\ell(\lambda)\leq m$.
Let $I_m(\lambda)$ denote the $m$-element set
$\{\lambda_{m}+1,\lambda_{m-1}+2,\dots,\lambda_{1}+m\}$.
For example,
if $m=4$ and $\lambda=(4,3,1)$,
then $I_4(\lambda)=\{1,3,6,8\}$.
We use this notation to fix a row/column index set of a given matrix.
The following proposition is useful in combination with
Theorem~\ref{msf}:
\begin{prop}
\label{prop:skew}
Let $m$ and $n$ be positive integers such that $m\leq n$.
Let $\lambda=(\lambda_1,\dots,\lambda_{m})$ be a partition such that 
$\ell(\lambda)\leq m$ and $\lambda_{1}=\ell(\lambda')\leq n-m$.
Here $\overline I=[n]\setminus I$ denote the complement in the set $[n]$.
\begin{enumerate}
\item[(i)]
If $m$ or $n-m$ is even then we have
$
(-1)^{s(\overline{I_{m}(\lambda)},I_{m}(\lambda))}
=(-1)^{s(I_{m}(\lambda),\overline{I_{m}(\lambda)})}
=(-1)^{|\lambda|}
$.
\item[(ii)]
Assume $m$ is even and $0\leq k\leq n-m$.
Then we have
$
\Pf\left(C_{n}(t)^{I_{m}(\lambda)}_{I_{m}(\lambda)}\right)=t^{c(\lambda)}
$
and
$
\Pf\left(R_{n}(t)^{I_{m}(\lambda)}_{I_{m}(\lambda)}\right)=t^{r(\lambda)}
$.
In particular,
we have
$
\Pf\left({S_{n}}^{I_{m}(\lambda)}_{I_{m}(\lambda)}\right)=1
$
for any $\lambda$,
and
$\Pf\left({R_{n}}^{I_{m}(\lambda)}_{I_{m}(\lambda)}\right)$
(resp. $\Pf\left({C_{n}}^{I_{m}(\lambda)}_{I_{m}(\lambda)}\right)$)
equals $1$ if all rows (resp. columns)
of $\lambda$ have even length,
and $0$ otherwise.
Further,
$
\lim\limits_{\varepsilon\to0}\Pf\left({L_{n}^{(m,k)}(\varepsilon)}^{I_{m}(\lambda)}_{I_{m}(\lambda)}\right)
$
equals $1$ if $\lambda_1\leq k$,
and $0$ otherwise.
\item[(iii)]
Assume $n-m$ is even and $0\leq k\leq n-m$.
Then, we have
$
\Pf\left(\bar C_{n}(t)^{\overline{I_{m}(\lambda)}}_{\overline{I_{m}(\lambda)}}\right)
=(-1)^{|\lambda|}\,t^{c(\lambda)}
$
and
$
\Pf\left(\bar R_{n}(t)^{\overline{I_{m}(\lambda)}}_{\overline{I_{m}(\lambda)}}\right)
=(-1)^{|\lambda|}\,t^{r(\lambda)}
$.
In particular,
we have
$
\Pf\left(\bar {S_{n}}^{\overline{I_{m}(\lambda)}}_{\overline{I_{m}(\lambda)}}\right)
=(-1)^{|\lambda|}
$
for any $\lambda$,
and
$\Pf\left({R_{n}}^{\overline{I_{m}(\lambda)}}_{\overline{I_{m}(\lambda)}}\right)$
(resp. $\Pf\left({C_{n}}^{\overline{I_{m}(\lambda)}}_{\overline{I_{m}(\lambda)}}\right)$)
equals $(-1)^{|\lambda|}$ if all rows (resp. columns)
of $\lambda$ have even length,
and $0$ otherwise.
Further,\\
$
\lim\limits_{\varepsilon\to0}
\varepsilon^{-\lfloor\frac{k}{2}\rfloor}
\Pf\left({\bar L_{n}^{(m,k)}(\varepsilon)}^{\overline{I_{m}(\lambda)}}_{\overline{I_{m}(\lambda)}}\right)
$
equals $(-1)^{|\lambda|}$ if $\lambda_1\leq k$,
and $0$ otherwise.
\end{enumerate}
\end{prop}
We can also use Lemma~\ref{lemma:product} to prove this proposition.
But the important part is that we can combine Proposition~\ref{prop:skew}
with Theorem~\ref{msf} to compute several sums of determinants.
For example,
if we take $B=\bar R_{n}(t)$ in \thetag{\ref{eq:msf}},
then we obtain
\begin{equation*}
\sum_{{\lambda}\atop{\ell(\lambda)\leq m,\ \ell(\lambda')\leq n-m}}
t^{r(\lambda)}\det(T_{I_{m}(\lambda)})
=\Pf\begin{pmatrix}
   O_m  & J_m T \\
-{}^tT J_m    & \bar R_{n}(t) \\
\end{pmatrix}
\end{equation*}
for any positive integers $n\geq m$.
If we put $t=1$ in the left-hand side,
then we obtain the sum of $\det(T_{I_{m}(\lambda)})$ 
where $\lambda$ runs over all partitions,
meanwhile,
if we put $t=0$,
then we obtain the sum of $\det(T_{I_{m}(\lambda)})$ 
where $\lambda$ runs over all even partitions.

In \cite{IW2} we defined the notion of copfaffian matrices.
Let $n$ be an even integer,
and let $A$ be a skew symmetric matrix of size $n$.
Thus $A^{\overline{\{i,j\}}}_{\overline{\{i,j\}}}$
denote the $(n-2)$ by $(n-2)$ skew symmetric sub-matrix
obtained by removing both the $i$th and $j$th rows and both the $i$th and $j$th columns of $A$
for $1\le i<j\le n$.
Let us define $\gamma(i,j)$ by
\begin{equation}
\gamma(i,j)=(-1)^{j-i-1} \Pf\left(A^{\overline{\{i,j\}}}_{\overline{\{i,j\}}}\right)
\end{equation}
for $1\le i<j\le n$.
We define the values of $\gamma(i,j)$ for $1\leq j\leq i\leq n$ so that $\gamma(j,i)=-\gamma(i,j)$ always holds.
Let $n$ be an even integer.
Given a skew symmetric matrix $A$ of size $n$,
let us call $\gamma(i,j)$ a \defterm{copfaffian} corresponding to $a_{ij}$
(or \defterm{$(i,j)$-copfaffian}),
and let $\hat A$ denote the skew symmetric matrix whose $(i,j)$th entry is $\gamma(i,j)$,
which we call the \defterm{copfaffian matrix} of $A$.
From Lemma~\ref{prod_pf},
we can see that,
when $n$ is even,
$\bar S_{n}$, $\bar R_{n}(t)$ and $\bar C_{n}(t)$ 
are the copfaffian matrices of $S_n$, $R_{n}(t)$ and $C_{n}(t)$
respectively,
and vise versa.

In the following,
we deduce another form \thetag{\ref{eq:msf_even}}, \thetag{\ref{eq:msf_odd}}
of the minor summation formula \thetag{\ref{eq:msf}}
which appeared in \cite{IW1}.
Let $n=2n'$ be an even integer and let $A$ be a non-singular $n\times n$ skew symmetric matrix.
In \cite[Theorem~2.6]{IW2},
we have obtained the fact that,
for any $I\in\binom{[n]}{2r}$,
\begin{equation}
\Pf\left(\hat A^{I}_{I}\right)
=(-1)^{s(\overline I,I)}\Pf(A)^{r-1}
\Pf\left(A^{\overline I}_{\overline I}\right).
\label{eq:plucker}
\end{equation}
In particular,
\thetag{\ref{eq:plucker}} implies $\Pf(\hat A)=\Pf(A)^{n'-1}$ and $\hat{\hat A}=\Pf(A)^{n'-2}A$.
If we take $B=\frac1{\Pf(A)}\hat A$ in \thetag{\ref{eq:msf}},
\thetag{\ref{eq:plucker}} also shows that
$(-1)^{s(I,\overline I)}\Pf(B^{\overline I}_{\overline I})
=\frac1{\Pf(A)}\Pf(A^{I}_{I}).$
From this,
one can deduce that
\begin{equation}
\sum_{I\in\binom{[n]}{m}}
\Pf(A^{I}_{I})\det(T_{I})
=\Pf(A)\Pf\begin{pmatrix}
   O_m  & J_m T \\
-{}^tT J_m    & \frac1{\Pf(A)}\hat A \\
\end{pmatrix}
=\Pf(Q_{ij})_{1\leq i,j\leq n},
\label{eq:msf_even}
\end{equation}
where $Q=(Q_{ij})=TA\,{}^t\! T$, and its entries are given by
\begin{equation}
\label{pf_msf}
Q_{ij}=\sum_{1\le k<l\le N} a_{kl} \det(T^{ij}_{kl}),
\qquad(1\leq i,j\leq m).
\end{equation}
Here we write $T^{ij}_{kl}$ for $T^{\{ij\}}_{\{kl\}}$.

Assume $m$ and $n=2n'-1$ is an odd integer and 
let $A_{0}=(\alpha_{ij})_{0\leq i,j\leq n}$ be a non-singular $(n+1)\times(n+1)$ skew symmetric matrix,
and let $\hat A_{0}$ denote its copfaffian matrix.
Let $T$ be an $m\times n$ matrix.
Write the $(n+1)\times(n+1)$ matrix
$B_{0}=\frac1{\Pf(A_{0})}\hat A_{0}
=\left(\begin{array}{c|c}
1&{b}\\\hline
-{}^{t}\!{b}&B
\end{array}\right)
$
with $n\times n$ matrix $B$,
and the $(m+1)\times(n+1)$ matrix 
$T^{0}=\left(\begin{array}{c|c}
1&O_{1,n}\\\hline
O_{m,1}&T
\end{array}\right)
$.
Then apply \thetag{\ref{eq:msf_even}} and we obtain
\begin{equation}
\sum_{I\in\binom{[n]}{m}}
\Pf({A^{0}}^{\{0\}\uplus I}_{\{0\}\uplus I})\det(T_{I})
=\Pf(A^{0})\Pf\begin{pmatrix}
   O_m  & J_m T \\
-{}^tT J_m    & B \\
\end{pmatrix}
=\Pf(Q_{ij})_{0\leq i,j\leq n},
\label{eq:msf_odd}
\end{equation}
where $Q_{0j}=\sum_{1\leq k\leq n}a_{0k}t_{jk}$
(see also \cite[Theorem~1]{IW1}).

\section{
Definitions and bijections
}\label{sec:bijections}

In this section we study three classes of (shifted) plane partitions
which are denoted by $\TSSCPP{n,m}$, $\TSPP{n,m}$ and $\CSPP{n,m}$,
and we establish bijections between them.
The set $\TSPP{n,m}$ is a generalization of the set $\TSPP{n}$
defined in \cite{MRR2},
and the set $\CSPP{n,m}$ is newly defined in this paper.
Thus the study of the totally symmetric plane partitions
reduce to the study of $\CSPP{n,m}$,
which we call the set of ``restricted column-strict plane partitions''
(see Definition~\ref{def:CSPP}).
In the later sections,
we intensively study the set $\CSPP{n,m}$
which enable us to reveal the several interesting properties
of this set of plane partitions. 
First of all we have to recall the basic definitions and notation 
concerning plane partitions.
For the general theory of plane partitions 
the reader may consult \cite{B,M,Sta1,Sta2,Ste1}.

A \defterm{plane partition} is an array
$\pi=(\pi_{ij})_{i,j\geq1}$
of nonnegative integers such that $\pi$ has finite support
(i.e. finitely many nonzero entries)
and is weakly decreasing in rows and columns.
If $\sum_{i,j\geq1}\pi_{ij}=n$,
then we write $|\pi|=n$ and say that
$\pi$ is a plane partition of $n$,
or $\pi$ has \defterm{weight} $n$.
A \defterm{part} of a plane partition $\pi=(\pi_{ij})_{i,j\geq1}$
is a positive entry $\pi_{ij}>0$.
The \defterm{shape} of $\pi$ is the ordinary partition $\lambda$
for which $\pi$ has $\lambda_i$ nonzero parts in the $i$th row.
We denote the shape of $\pi$ by $\SHAPE{\pi}$.
We also say that $\pi$ has $r$ \defterm{rows} if $r=\ell(\lambda)$.
Similarly,
$\pi$ has $s$ \defterm{columns} if $s=\ell(\lambda')$.
A plane partition is said to be \defterm{column-strict}
if it is strictly decreasing in columns.
For example,
\begin{equation*}
\begin{array}{ccccccccc}
5&5&4&3&2&1\\
4&4&2&2&1\\
2&2&1&1\\
1\\
\end{array}
\end{equation*}
is a column-strict plane partition and has shape $(6,5,4,1)$,
$4$ rows, $6$ columns and weight $40$.
Consider the elements of $\Pos^3$,
regarded as the lattice points of $\Real^3$ in the positive orthant.
The \defterm{Ferrers graph} $F(\pi)$ of $\pi$ is the set of all lattice points $(i,j,k)\in\Pos^3$
such that $k\leq \pi_{ij}$.
A subset $F$ of $\Pos^3$ is a Ferrers graph 
if and only if it satisfies
\[
x_1\leq x_2,\ y_1\leq y_2,\ z_1\leq z_2\text{ and }(x_2,y_2,z_2)\in F\Rightarrow(x_1,y_1,z_1)\in F.
\]
Hereafter we identify a plane partition and its Ferrers graph,
and write $\pi$ for $F(\pi)$.
The symmetric group $S_3$ is acting on $\Pos^3$ as permutations of the coordinate axises.
A plane partition is said to be \defterm{totally symmetric}
if its Ferrers graph is mapped to itself under all 6 permutations in $S_3$.

In this section we mainly consider three classes,
i.e. $\TSSCPP{n,m}$, $\TSPP{n,m}$ and $\CSPP{n,m}$,
of (shifted) plane partitions
and construct bijections between them.
First of all,
the following set $\CSPP{n,m}$ of plane partitions
plays a crucial role throughout this paper 
in the study of the totally symmetric plane partitions.
\begin{definition}
\label{def:CSPP}
Let $m$ and $n\geq1$ be nonnegative integers.
Let $\CSPP{n,m}$ denote the set
of plane partitions $c=(c_{ij})_{1\leq i, j}$
subject to the constraints that
\begin{enumerate}
\item[(C1)]
$c$ has at most $n$ columns;
\item[(C2)]
$c$ is column-strict 
and each part in the $j$th column does not exceed $n+m-j$.
\end{enumerate}
We call an element of $\CSPP{n,m}$ a \defterm{restricted column-strict plane partition} (abbreviated to RCSPP).
When $m=0$,
we write $\CSPP{n}$ for $\CSPP{n,0}$.
If a part in the $j$th column of $c$ is equal to $n+m-j$
(that can happen only in the first row, i.e. $c_{1j}=n+m-j$),
we call the part a \defterm{saturated part}.
Further we define two subclasses of $\CSPP{n,m}$.
Let $\RCSPP{n,m}$ denote the set of plane partitions $c$ in $\CSPP{n,m}$
where each row has even length.
and let $\CCSPP{n,m}$ denote the set of plane partitions $c$ in $\CSPP{n,m}$
with each column of even length.
We also write  $\RCSPP{n}$ (resp. $\CCSPP{n}$) for $\RCSPP{n,0}$ (resp. $\CCSPP{n,0}$).
\end{definition}
For instance,
$\CSPP{1,2}$ consists of the following 4 plane partitions:
\begin{equation*}
\begin{array}{cc}
\emptyset\\
&
\end{array}
\qquad
\begin{array}{cc}
1\\
&
\end{array}
\qquad
\begin{array}{cc}
\pmb{2}\\
&
\end{array}
\qquad
\begin{array}{cc}
\pmb{2}\\
1
\end{array}
\end{equation*}
In the above four RCSPPs
the boldfaced letters $\pmb{2}$ stand for saturated parts
since they are in the first column and equals $n+m-1=2$.

Next we define the following set $\TSPP{n,m}$ of shifted
plane partitions
which is a generalization of $\TSPP{n}$ defined in \cite[pp.281]{MRR2}.
Let $\mu$ be a strict partition.
A \defterm{shifted plane partition} $\tau$ of \defterm{shifted shape} $\mu$ 
is an arbitrary filling of the cells of $\mu$ with nonnegative integers 
such that
each entry is weakly decreasing in rows and columns.
In this paper we allow parts to be zero 
for shifted plane partitions of a fixed shifted shape $\mu$.
\begin{definition}
\label{def:TSPP}
(See \cite[Theorem~1]{Kr1}).
Let $m$ and $n\geq1$ be nonnegative integers.
Let $\TSPP{n,m}$ denote the set
of shifted plane partitions $b=(b_{ij})_{1\leq i\leq j}$
subject to the constraints that
\begin{enumerate}
\item[(B1)]
the shifted shape of $b$ is  $(n+m-1,n+m-2,\dots,2,1)$;
\item[(B2)]
$\max\{n-i,0\}\leq b_{ij}\leq n$ for $1\leq i\leq j\leq n+m-1$.
\end{enumerate}
When $m=0$,
we write $\TSPP{n}$ for $\TSPP{n,0}$.
In this paper we call an element of $\TSPP{n,m}$ a \defterm{triangular shifted plane partition}
(abbreviated to TSPP).
\end{definition}
When $n=1$ and $m=2$,
$\TSPP{1,2}$ consists of the following 4 elements:
\begin{equation*}
\smyoung{
1&1\\
\blank&1\\}
\qquad
\smyoung{
1&1\\
\blank&0\\
}
\qquad
\smyoung{
1&0\\
\blank&0\\
}
\qquad
\smyoung{
0&0\\
\blank&0\\
}
\end{equation*}

The final object $\TSSCPP{n,m}$ we need to define in this section
is a subclass of totally symmetric self-complementary plane partitions.
Before we proceed to the definition of $\TSSCPP{n,m}$,
we need to define additional terminology and symbols concerning plane partitions.
Let $\BOX{r,s,t}=[r]\times[s]\times[t]$ denote the $r\times s\times t$ box.
Assume $r$, $s$ and $t$ are all even.
We divide this box into the eight regions 
$\BOXU{r,s,t}{+++}$, $\BOXU{r,s,t}{++-}$, $\BOXU{r,s,t}{+-+}$, $\BOXU{r,s,t}{+--}$, 
$\BOXU{r,s,t}{-++}$, $\BOXU{r,s,t}{-+-}$, $\BOXU{r,s,t}{--+}$ and $\BOXU{r,s,t}{---}$
depending on each of $x-r/2$, $y-s/2$ and $z-t/2$ is plus ($>0$) or minus ($\leq0$).
For example $\BOXU{r,s,t}{-+-}=[1,r/2]\times[s/2+1,s]\times[1,t/2]$.
Further we use the notation $\BOXU{r,s,t}{+}=\BOXU{r,s,t}{+++}\uplus\BOXU{r,s,t}{++-}\uplus\BOXU{r,s,t}{+-+}\uplus\BOXU{r,s,t}{-++}$
and
$\BOXU{r,s,t}{-}=\BOXU{r,s,t}{+--}\uplus\BOXU{r,s,t}{-+-}\uplus\BOXU{r,s,t}{--+}\uplus\BOXU{r,s,t}{---}$.
More generally we write $\BOXC{r,s,t}{a,b,c}=[a-r/2+1,a+r/2]\times[b-s/2+1,b+s/2]\times[c-t/2+1,c+t/2]$
for the $r\times s\times t$ box centered at $(a,b,c)$.
We also use the notation $\BOXCU{r,s,t}{a,b,c}{\pm\pm\pm}$ as the same meaning as above
where each stands for one of the eight regions of $\BOXC{r,s,t}{a,b,c}$.
For example $\BOXCU{r,s,t}{a,b,c}{+-+}=[a+1,a+r/2]\times[b-s/2+1,b]\times[c+1,c+t/2]$.
The symbols $\BOXCU{r,s,t}{a,b,c}{\pm}$ should be defined similarly.
The involution $\CMP{r,s,t}:(x,y,z)\mapsto(r+1-x,s+1-y,t+1-z)$ is called the \defterm{complementation}.
When $r=s=t$ and $a=b=c$,
we use the abbreviation 
$\BOX{r}$ for $\BOX{r,r,r}$,
$\BOXU{r}{\pm\pm\pm}$ for $\BOXU{r,r,r}{\pm\pm\pm}$,
$\BOXU{r}{\pm}$ for $\BOXU{r,r,r}{\pm}$
and $\CMP{r}$ for $\CMP{r,r,r}$.
The symbols $\BOXC{r}{a}$, $\BOXCU{r}{a}{\pm\pm\pm}$ and $\BOXCU{r}{a}{\pm}$ should be interpreted similarly.

\begin{figure}[b]
\begin{center}
\setlength{\unitlength}{1pt} 
\setlength{\unitlength}{0.5mm}
\begin{picture}(200,200)
\put(  0,100){\line(0,1){30}}
\put( 10,100){\line(0,1){30}}
\put( 20,100){\line(0,1){30}}
\put( 30,100){\line(0,1){30}}
\put( 35,105){\line(0,1){30}}
\put( 40,110){\line(0,1){30}}
\put( 45,115){\line(0,1){60}}
\put( 55,115){\line(0,1){30}}
\put( 65,115){\line(0,1){30}}
\put( 75,115){\line(0,1){30}}
\put( 80,120){\line(0,1){30}}
\put( 85,125){\line(0,1){30}}
\put( 90,130){\line(0,1){30}}
\put( 15,145){\line(0,1){30}}
\put( 25,145){\line(0,1){30}}
\put( 35,145){\line(0,1){30}}
\put( 50,150){\line(0,1){30}}
\put( 55,155){\line(0,1){30}}
\put( 60,160){\line(0,1){30}}
\put(  0,100){\line(1,0){30}}
\put(  0,110){\line(1,0){30}}
\put(  0,120){\line(1,0){30}}
\put(  0,130){\line(1,0){30}}
\put(  5,135){\line(1,0){30}}
\put( 10,140){\line(1,0){30}}
\put( 15,145){\line(1,0){60}}
\put( 15,155){\line(1,0){30}}
\put( 15,165){\line(1,0){30}}
\put( 15,175){\line(1,0){30}}
\put( 20,180){\line(1,0){30}}
\put( 25,185){\line(1,0){30}}
\put( 30,190){\line(1,0){30}}
\put( 45,115){\line(1,0){30}}
\put( 45,125){\line(1,0){30}}
\put( 45,135){\line(1,0){30}}
\put( 50,150){\line(1,0){30}}
\put( 55,155){\line(1,0){30}}
\put( 60,160){\line(1,0){30}}
\put(  0,130){\line(1,1){15}}
\put( 10,130){\line(1,1){15}}
\put( 20,130){\line(1,1){15}}
\put( 30,130){\line(1,1){30}}
\put( 30,120){\line(1,1){15}}
\put( 30,110){\line(1,1){15}}
\put( 30,100){\line(1,1){15}}
\put( 15,175){\line(1,1){15}}
\put( 25,175){\line(1,1){15}}
\put( 35,175){\line(1,1){15}}
\put( 45,175){\line(1,1){15}}
\put( 45,165){\line(1,1){15}}
\put( 45,155){\line(1,1){15}}
\put( 55,145){\line(1,1){15}}
\put( 65,145){\line(1,1){15}}
\put( 75,145){\line(1,1){15}}
\put( 75,135){\line(1,1){15}}
\put( 75,125){\line(1,1){15}}
\put( 75,115){\line(1,1){15}}
\put(100,100){\line(0,1){30}}
\put(110,100){\line(0,1){30}}
\put(120,100){\line(0,1){30}}
\put(130,100){\line(0,1){20}}
\put(135,105){\line(0,1){30}}
\put(140,110){\line(0,1){30}}
\put(145,125){\line(0,1){40}}
\put(155,115){\line(0,1){30}}
\put(165,115){\line(0,1){30}}
\put(175,115){\line(0,1){20}}
\put(180,120){\line(0,1){30}}
\put(185,125){\line(0,1){30}}
\put(190,130){\line(0,1){30}}
\put(115,155){\line(0,1){20}}
\put(125,145){\line(0,1){30}}
\put(135,145){\line(0,1){30}}
\put(150,150){\line(0,1){30}}
\put(155,155){\line(0,1){30}}
\put(160,170){\line(0,1){20}}
\put(110,140){\line(0,1){10}}
\put(120,140){\line(0,1){10}}
\put(125,125){\line(0,1){10}}
\put(140,170){\line(0,1){10}}
\put(150,110){\line(0,1){10}}
\put(165,155){\line(0,1){10}}
\put(170,140){\line(0,1){10}}
\put(170,160){\line(0,1){10}}
\put(100,100){\line(1,0){30}}
\put(100,110){\line(1,0){30}}
\put(100,120){\line(1,0){30}}
\put(100,130){\line(1,0){20}}
\put(105,135){\line(1,0){30}}
\put(110,140){\line(1,0){30}}
\put(125,145){\line(1,0){40}}
\put(115,155){\line(1,0){30}}
\put(115,165){\line(1,0){30}}
\put(115,175){\line(1,0){20}}
\put(120,180){\line(1,0){30}}
\put(125,185){\line(1,0){30}}
\put(130,190){\line(1,0){30}}
\put(155,115){\line(1,0){20}}
\put(145,125){\line(1,0){30}}
\put(145,135){\line(1,0){30}}
\put(150,150){\line(1,0){30}}
\put(155,155){\line(1,0){30}}
\put(170,160){\line(1,0){20}}
\put(110,150){\line(1,0){10}}
\put(125,125){\line(1,0){10}}
\put(140,110){\line(1,0){10}}
\put(140,120){\line(1,0){10}}
\put(140,170){\line(1,0){10}}
\put(155,165){\line(1,0){10}}
\put(160,170){\line(1,0){10}}
\put(170,140){\line(1,0){10}}
\put(100,130){\line(1,1){10}}
\put(110,130){\line(1,1){15}}
\put(120,130){\line(1,1){15}}
\put(135,135){\line(1,1){20}}
\put(130,120){\line(1,1){15}}
\put(130,110){\line(1,1){15}}
\put(130,100){\line(1,1){10}}
\put(115,175){\line(1,1){15}}
\put(125,175){\line(1,1){15}}
\put(135,175){\line(1,1){15}}
\put(150,180){\line(1,1){10}}
\put(145,165){\line(1,1){15}}
\put(145,155){\line(1,1){15}}
\put(155,145){\line(1,1){15}}
\put(165,145){\line(1,1){15}}
\put(180,150){\line(1,1){10}}
\put(175,135){\line(1,1){15}}
\put(175,125){\line(1,1){15}}
\put(175,115){\line(1,1){15}}
\put(110,150){\line(1,1){ 5}}
\put(120,120){\line(1,1){ 5}}
\put(120,150){\line(1,1){ 5}}
\put(135,165){\line(1,1){ 5}}
\put(150,110){\line(1,1){ 5}}
\put(150,120){\line(1,1){ 5}}
\put(165,135){\line(1,1){ 5}}
\put(165,165){\line(1,1){ 5}}
\put(  0,  0){\line(0,1){30}}
\put( 10,  0){\line(0,1){30}}
\put( 20,  0){\line(0,1){20}}
\put( 30,  0){\line(0,1){10}}
\put( 15, 25){\line(0,1){10}}
\put(  5, 35){\line(0,1){10}}
\put( 15, 35){\line(0,1){10}}
\put( 10, 50){\line(0,1){10}}
\put( 15, 35){\line(0,1){10}}
\put( 20, 40){\line(0,1){20}}
\put( 25, 15){\line(0,1){20}}
\put( 35,  5){\line(0,1){30}}
\put( 40, 20){\line(0,1){20}}
\put( 45, 25){\line(0,1){40}}
\put( 55, 25){\line(0,1){20}}
\put( 65, 15){\line(0,1){30}}
\put( 75, 15){\line(0,1){10}}
\put( 80, 20){\line(0,1){20}}
\put( 85, 25){\line(0,1){30}}
\put( 90, 30){\line(0,1){30}}
\put( 15, 65){\line(0,1){10}}
\put( 25, 45){\line(0,1){30}}
\put( 35, 45){\line(0,1){20}}
\put( 50, 50){\line(0,1){20}}
\put( 55, 55){\line(0,1){30}}
\put( 60, 80){\line(0,1){10}}
\put( 30, 70){\line(0,1){10}}
\put( 40, 70){\line(0,1){10}}
\put( 45, 75){\line(0,1){10}}
\put( 45,  5){\line(0,1){10}}
\put( 50, 10){\line(0,1){10}}
\put( 60, 10){\line(0,1){10}}
\put( 65, 55){\line(0,1){20}}
\put( 70, 70){\line(0,1){10}}
\put( 75, 55){\line(0,1){10}}
\put( 80, 60){\line(0,1){10}}
\put( 70, 30){\line(0,1){20}}
\put( 75, 45){\line(0,1){10}}
\put(  0,  0){\line(1,0){30}}
\put(  0, 10){\line(1,0){30}}
\put(  0, 20){\line(1,0){20}}
\put(  0, 30){\line(1,0){10}}
\put( 25, 15){\line(1,0){20}}
\put( 35,  5){\line(1,0){10}}
\put( 35,  5){\line(1,0){10}}
\put( 15, 25){\line(1,0){20}}
\put( 40, 20){\line(1,0){20}}
\put( 50, 10){\line(1,0){10}}
\put(  5, 35){\line(1,0){30}}
\put( 20, 40){\line(1,0){20}}
\put( 25, 45){\line(1,0){40}}
\put( 25, 55){\line(1,0){20}}
\put(  5, 45){\line(1,0){10}}
\put( 10, 50){\line(1,0){10}}
\put( 10, 60){\line(1,0){10}}
\put( 15, 65){\line(1,0){30}}
\put( 15, 75){\line(1,0){10}}
\put( 20, 80){\line(1,0){20}}
\put( 30, 70){\line(1,0){20}}
\put( 60, 80){\line(1,0){10}}
\put( 45, 75){\line(1,0){20}}
\put( 55, 65){\line(1,0){20}}
\put( 70, 70){\line(1,0){10}}
\put( 25, 85){\line(1,0){30}}
\put( 30, 90){\line(1,0){30}}
\put( 65, 15){\line(1,0){10}}
\put( 45, 25){\line(1,0){30}}
\put( 45, 35){\line(1,0){20}}
\put( 50, 50){\line(1,0){20}}
\put( 55, 55){\line(1,0){30}}
\put( 80, 60){\line(1,0){10}}
\put( 70, 30){\line(1,0){10}}
\put( 70, 40){\line(1,0){10}}
\put( 75, 45){\line(1,0){10}}
\put(  0, 30){\line(1,1){ 5}}
\put( 10, 30){\line(1,1){15}}
\put( 10, 20){\line(1,1){ 5}}
\put( 20, 20){\line(1,1){ 5}}
\put( 20, 10){\line(1,1){ 5}}
\put( 25, 35){\line(1,1){10}}
\put( 35, 35){\line(1,1){20}}
\put( 35, 25){\line(1,1){10}}
\put( 30, 10){\line(1,1){15}}
\put( 30,  0){\line(1,1){ 5}}
\put( 15, 75){\line(1,1){15}}
\put( 25, 75){\line(1,1){15}}
\put( 40, 80){\line(1,1){10}}
\put( 55, 85){\line(1,1){ 5}}
\put( 45, 65){\line(1,1){15}}
\put( 45, 55){\line(1,1){10}}
\put( 55, 45){\line(1,1){10}}
\put( 65, 45){\line(1,1){15}}
\put( 85, 55){\line(1,1){ 5}}
\put( 80, 40){\line(1,1){10}}
\put( 75, 25){\line(1,1){15}}
\put( 75, 15){\line(1,1){15}}
\put(  5, 45){\line(1,1){ 5}}
\put( 10, 60){\line(1,1){ 5}}
\put( 20, 60){\line(1,1){10}}
\put( 35, 65){\line(1,1){10}}
\put( 15, 45){\line(1,1){10}}
\put( 45,  5){\line(1,1){ 5}}
\put( 60, 10){\line(1,1){ 5}}
\put( 60, 20){\line(1,1){10}}
\put( 65, 35){\line(1,1){10}}
\put( 45, 15){\line(1,1){10}}
\put( 75, 65){\line(1,1){ 5}}
\put( 65, 65){\line(1,1){ 5}}
\put( 65, 75){\line(1,1){ 5}}
\put(100,  0){\line(0,1){30}}
\put(110,  0){\line(0,1){30}}
\put(120,  0){\line(0,1){10}}
\put(130,  0){\line(0,1){10}}
\put(130, 70){\line(0,1){10}}
\put(135, 75){\line(0,1){10}}
\put(145, 75){\line(0,1){10}}
\put(135,  5){\line(0,1){30}}
\put(140, 20){\line(0,1){20}}
\put(115, 15){\line(0,1){40}}
\put(125, 15){\line(0,1){20}}
\put(145, 25){\line(0,1){40}}
\put(155, 25){\line(0,1){20}}
\put(165, 15){\line(0,1){30}}
\put(175, 15){\line(0,1){10}}
\put(180, 20){\line(0,1){10}}
\put(185, 25){\line(0,1){30}}
\put(190, 30){\line(0,1){30}}
\put(115, 65){\line(0,1){10}}
\put(125, 45){\line(0,1){30}}
\put(135, 45){\line(0,1){20}}
\put(150, 50){\line(0,1){20}}
\put(155, 55){\line(0,1){30}}
\put(160, 80){\line(0,1){10}}
\put(105, 35){\line(0,1){20}}
\put(120, 40){\line(0,1){20}}
\put(125, 25){\line(0,1){10}}
\put(145,  5){\line(0,1){10}}
\put(155,  5){\line(0,1){10}}
\put(160, 10){\line(0,1){10}}
\put(165, 55){\line(0,1){20}}
\put(170, 30){\line(0,1){20}}
\put(175, 35){\line(0,1){20}}
\put(175, 55){\line(0,1){20}}
\put(180, 60){\line(0,1){20}}
\put(100,  0){\line(1,0){30}}
\put(100, 10){\line(1,0){30}}
\put(100, 20){\line(1,0){10}}
\put(100, 30){\line(1,0){10}}
\put(105, 35){\line(1,0){30}}
\put(105, 45){\line(1,0){10}}
\put(105, 55){\line(1,0){10}}
\put(110, 60){\line(1,0){10}}
\put(120, 40){\line(1,0){20}}
\put(125, 45){\line(1,0){40}}
\put(125, 55){\line(1,0){20}}
\put(115, 65){\line(1,0){30}}
\put(115, 75){\line(1,0){10}}
\put(120, 80){\line(1,0){10}}
\put(125, 85){\line(1,0){30}}
\put(130, 90){\line(1,0){30}}
\put(165, 15){\line(1,0){10}}
\put(145, 25){\line(1,0){30}}
\put(145, 35){\line(1,0){20}}
\put(150, 50){\line(1,0){20}}
\put(155, 55){\line(1,0){30}}
\put(180, 60){\line(1,0){10}}
\put(115, 15){\line(1,0){40}}
\put(115, 25){\line(1,0){20}}
\put(135,  5){\line(1,0){20}}
\put(140, 20){\line(1,0){20}}
\put(130, 70){\line(1,0){20}}
\put(135, 75){\line(1,0){40}}
\put(155, 65){\line(1,0){20}}
\put(160, 80){\line(1,0){20}}
\put(175, 45){\line(1,0){10}}
\put(175, 35){\line(1,0){10}}
\put(170, 30){\line(1,0){10}}
\put(100, 30){\line(1,1){ 5}}
\put(105, 55){\line(1,1){10}}
\put(110, 30){\line(1,1){15}}
\put(115, 45){\line(1,1){10}}
\put(115, 55){\line(1,1){20}}
\put(125, 35){\line(1,1){10}}
\put(135, 35){\line(1,1){20}}
\put(135, 25){\line(1,1){10}}
\put(130, 10){\line(1,1){15}}
\put(130,  0){\line(1,1){ 5}}
\put(115, 75){\line(1,1){15}}
\put(125, 75){\line(1,1){15}}
\put(145, 85){\line(1,1){ 5}}
\put(155, 85){\line(1,1){ 5}}
\put(145, 65){\line(1,1){15}}
\put(145, 55){\line(1,1){10}}
\put(155, 45){\line(1,1){10}}
\put(165, 45){\line(1,1){15}}
\put(185, 55){\line(1,1){ 5}}
\put(185, 45){\line(1,1){ 5}}
\put(175, 25){\line(1,1){15}}
\put(175, 15){\line(1,1){15}}
\put(110, 10){\line(1,1){ 5}}
\put(110, 20){\line(1,1){ 5}}
\put(120, 10){\line(1,1){ 5}}
\put(135, 65){\line(1,1){10}}
\put(145, 15){\line(1,1){10}}
\put(155, 15){\line(1,1){10}}
\put(155,  5){\line(1,1){10}}
\put(165, 25){\line(1,1){10}}
\put(165, 35){\line(1,1){10}}
\put(165, 75){\line(1,1){ 5}}
\put(175, 75){\line(1,1){ 5}}
\put(175, 65){\line(1,1){ 5}}
\end{picture}
\caption{TSSCPP ($n=1$, $m=2$)}\label{figure:tsscpp}
\end{center}
\end{figure}
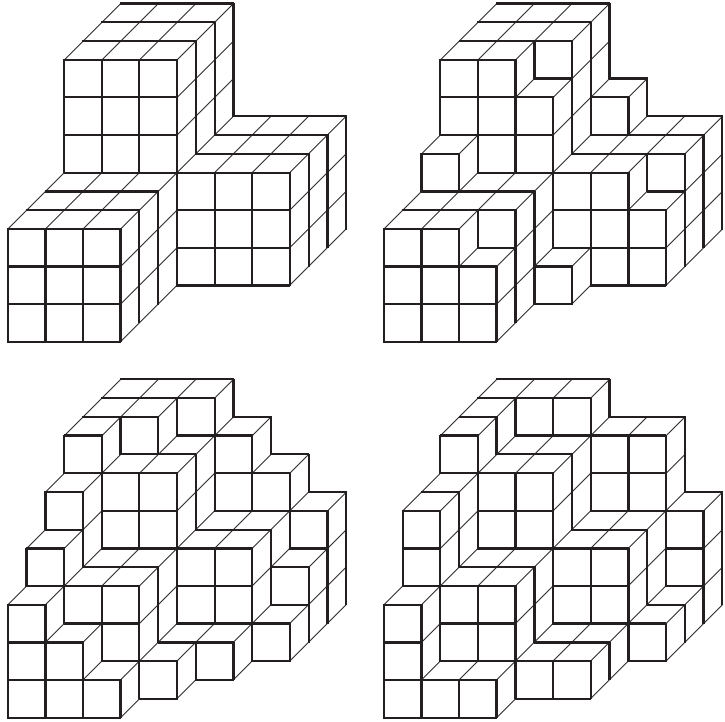

A plane partition
$\pi\subseteq\BOX{r,s,t}$ is $(r,s,t)$-self-complementary if we have,
for all $p\in\BOX{r,s,t}$,
$p\in\pi$ if and only if $\CMP{r,s,t}(p)\not\in\pi$.
Let $\TSSCPP{n}$ denote the set of all plane partitions
which is contained in the cube $\BOX{2n}$,
$(2n,2n,2n)$-self-complementary and totally symmetric.
\begin{definition}
\label{def:TSSCPP}
An element of $\TSSCPP{n}$ is called 
a totally symmetric self-complementary plane partition
(abbreviated to TSSCPP) of size $n$.
For nonnegative integers $m$ and $n\geq1$,
let $\TSSCPP{n,m}$ denote the set of TSSCPPs $\pi\in\TSSCPP{n+m}$ of size $(n+m)$
which satisfy
\begin{enumerate}
\item[(T)]
each $p\in\pi\cap \BOXC{2m}{n}$ must be contained in $\BOXU{2(n+m)}{-}$.
\end{enumerate}
Note that,
when $\pi\in\TSSCPP{n,m}$,
$\pi\cap \BOXC{2m}{n}$ is uniquely determined by the condition (T),
i.e. $\pi\cap \BOXC{2m}{n}=\BOXCU{2m}{n}{-}$.
\end{definition}
For instance
$\TSSCPP{1,2}$ is composed of four elements which are designated in Figure~\ref{figure:tsscpp}.
%

Mills, Robbins and Rumsey
have constructed a bijection between $\TSSCPP{n}$ and $\TSPP{n}$
(see \cite[Theorem~1]{MRR2}).
The set $\TSPP{n,m}$ is first considered in \cite[Theorem~1 (2)]{Kr1}.
The main results of this section is the bijections between
the three classes $\TSSCPP{n,m}$, $\TSPP{n,m}$ and $\CSPP{n,m}$.
Thus, in the later sections,
 we intensively study the properties of $\CSPP{n,m}$.
The first theorem establish  a bijection between $\TSSCPP{n,m}$ and $\TSPP{n,m}$,
and this bijection is a generalization of the bijection between $\TSSCPP{n}$ and $\TSPP{n}$
constructed in \cite[Theorem~1]{MRR2}.
\begin{theorem}
\label{thm:bijection_TSPP}
Let $m$ and $n\geq1$ be nonnegative integers and $(a_{ij})$
be a TSSCPP in $\TSSCPP{n,m}$.
Associate to the array $(a_{ij})$
the triangular array $a_{i+1,j+1}-(n+2m)$ with $1\leq i\leq j\leq n+m-1$.
Then the triangular array is in $\TSPP{n,m}$ 
and this mapping from $\TSSCPP{n,m}$ to $\TSPP{n,m}$ is a bijection.
\end{theorem}
In the set-theoretic approach this bijection $\TSSCPP{n,m}\rightarrow\TSPP{n,m},\ a\mapsto b$ can be restated as follows:
\begin{align}
\forall(x,y,z)\in\Bbb{P}^3&\text{ such that $1\leq x\leq n+m-1$:}
\nonumber\\
&(x,y,z)\in b\Leftrightarrow(x+1,y+1,z+n+2m)\in a.
\label{eq:set_TSPP}
\end{align}
The second theorem establish a bijection between $\TSSCPP{n,m}$ and $\CSPP{n,m}$,
\begin{theorem}
\label{thm:bijection_CSPP}
Let $m$ and $n\geq1$ be nonnegative integers and $a=(a_{ij})_{1\leq i,j\leq2(n+m)}$
be a TSSCPP in $\TSSCPP{n,m}$.
Associate to the array $a$
the array $a_{i+n+m,j}-(n+m)-i+1$ with $1\leq i,j\leq n+m$.
Then,
by ignoring the zeros and negative entries,
this array is a plane partition in $\CSPP{n,m}$ 
and this mapping is a bijection of $\TSSCPP{n,m}$ onto $\CSPP{n,m}$.
\end{theorem}
This bijection $\TSSCPP{n,m}\rightarrow\CSPP{n,m},\ a\mapsto c$ can be restated as follows:
\begin{align}
\forall(x,y,z)\in\Bbb{P}^3&\text{ such that $1\leq y+z\leq n+m$:}
\nonumber\\
&(x,y,z)\in c\Leftrightarrow(x+n+m,y,x+z+n+m-1)\in a.
\label{eq:set_CSPP}
\end{align}
Finally
as a corollary of these theorems 
we establish a bijection between $\TSPP{n,m}$ and $\CSPP{n,m}$,
Let $c=(c_{ij})_{1\leq i\leq n+m,1\leq j\leq n}$ be a RCSPP in $\CSPP{n,m}$
and let $k$ be a positive integer.
Let $c_{\,\geq k}$ denote the plane partition formed by the parts $\geq k$.
Let
\begin{equation}
\BORDER{i}{k}{c}=\sharp\{l:c_{i,l}\geq k\}
\label{eq:border}
\end{equation}
denote the length of the $i$th row of $c_{\,\geq k}$,
i.e.
the rightmost column
 containing a letter $\geq k$ in the $i$th row of $c$.
\begin{corollary}
\label{cor:bijection_TSPPtoCSPP}
Let $m$ and $n\geq1$ be nonnegative integers
and $c=(c_{ij})_{1\leq i\leq n+m,1\leq j\leq n}$ be a RCSPP in $\CSPP{n,m}$.
Associate to the array $c=(c_{ij})_{1\leq i\leq n+m,1\leq j\leq n}$ 
the array $b=(b_{ij})_{1\leq i\leq j\leq n+m-1}$ defined by
\begin{equation}
n-b_{ij}=\BORDER{n+m-j}{1-i+j}{c}
\label{eq:CSPPtoTSPP}
\end{equation}
with $1\leq i\leq j\leq n+m-1$.
Then $b$ is in $\TSPP{n,m}$,
 and this mapping $\BIJ{n,m}$,
which associate to a RCSPP $c$ the TSPP $b=\BIJ{n.m}(c)$,
is a bijection of $\CSPP{n,m}$ onto $\TSPP{n,m}$.
\end{corollary}
In the set-theoretic interpretation 
the mapping $\CSPP{n,m}\rightarrow\TSPP{n,m}$ $c\mapsto b$ is formulated by
\begin{align}
\forall(x,y,z)\in\Bbb{P}^3&\text{ such that $1\leq x\leq y\leq n+m-1$:}
\nonumber\\
&(x,y,z)\in b\Leftrightarrow(n+m-y,n+1-z,1+y-x)\not\in c,
\label{eq:set_CSPPtoTSPP}
\end{align}
and
the mapping $\TSPP{n,m}\rightarrow\CSPP{n,m}$ $b\mapsto c$ is formulated by
\begin{align}
\forall(x,y,z)\in\Bbb{P}^3&\text{ such that $1\leq x+z\leq n+m$ and $1\leq y\leq n$:}
\nonumber\\
&(x,y,z)\in c\Leftrightarrow(n+m+1-x-z,n+m-x,n+1-y)\not\in b.
\label{eq:set_TSPPtoCSPP}
\end{align}

\begin{demo}{Proof of Theorem~\ref{thm:bijection_TSPP}}
Since $\TSSCPP{n,m}$ is a subset of $\TSSCPP{n+m}$,
we can imitate the proof of Theorem~1 in \cite[pp.280]{MRR2}.
Suppose that $(a_{ij})$ is in $\TSSCPP{n,m}$.
Associate to the array $(a_{ij})$ the subarray consisting of those parts $a_{ij}$ with $1\leq i\leq j\leq n+m$.
Then the subarray $(a_{ij})_{1\leq i\leq j\leq n+m}$ satisfies the conditions that
\begin{enumerate}
\item[(A1)]
all rows and columns are weakly decreasing;
\item[(A2)]
$\max\left\{2(n+m)+1-i,n+2m\right\}\leq a_{ij}\leq 2(n+m)$ for $1\leq i\leq n+m$.
\end{enumerate}
Exactly the same argument as in the proof of Theorem~1 in \cite[pp.280]{MRR2}
works to show that this mapping defines a bijection from $\TSSCPP{n,m}$ to the set
of triangular arrays $(a_{ij})_{1\leq i\leq j\leq n+m}$ satisfying the constraints (C1) and (C2).
The reader can consult \cite{MRR2} to fill the details.
Since $a_{ij}\geq n+2m$ and the top row of $(a_{ij})_{1\leq i\leq j\leq n+m}$ consists of $(n+m)$ $2(n+m)$'s,
no informations is lost by omitting the top row and subtracting $(n+2m)$ from the remaining parts.
When we do this,
we obtain $\TSPP{n,m}$.
\end{demo}

\begin{demo}{Proof of Theorem~\ref{thm:bijection_CSPP}}
Suppose that $(a_{ij})$ is in $\TSSCPP{n,m}$.
Associate to the array $a=(a_{ij})_{1\leq i,j\leq2(n+m)}$ the array consisting of $a_{i+n+m,j}-(n+m)$ with $1\leq i,j\leq n+m$.
If $a_{i+n+m,j}-(n+m)$ is negative, then we ignore the part and regard it as zero.
Then the array $\gamma=(\gamma_{ij})=(a_{i+n+m,j}-(n+m))_{1\leq i,j\leq n+m}$ satisfies the conditions that
\begin{enumerate}
\item[(D1)]\label{D1}
$\gamma$ is a plane partition with at most $n$ columns;
\item[(D2)]\label{D2}
each column of $\gamma$ is a self-conjugate partition;
\item[(D3)]\label{D3}
each part in the $j$th column of $\gamma$ is $\leq n+m-j$.
\end{enumerate}
We shall show that this defines a bijection $(a_{ij})\mapsto\gamma$
where $\gamma$ satisfies the constraints \thetag{D1}, \thetag{D2} and \thetag{D3}.
In fact $\gamma$ is trivially a plane partition,
and \thetag{D2} is easy to see from the symmetry.
We shall show that it satisfies \thetag{D3}.
Suppose that $\gamma_{ij}>n+m-j$,
i.e. $a_{i+n+m,j}>2(n+m)-j$.
Then we would have $(i+n+m,j,2(n+m)-j+1)\in a$
which implies that $(n+m+1-i,j,2(n+m)-j+1)\not\in a$ by self-complementarity and total symmetry.
Since $n+m+1-i<i+n+m$,
this contradict the fact that $a$ is a plane partition.
Next we shall show that $\gamma$ has at most $n$ columns.
In fact if $\gamma$ had more than $n$ columns
then $(1,n+1,1)\in\gamma$ ,
i.e. $(n+m+1,n+1,n+m+1)\in a$.
Since $(n+m+1,n+1,n+m+1)\in\BOXCU{2m}{n+m}{+-+}$,
this contradicts the condition (T).
Thus $\gamma$ satisfies \thetag{D1}, \thetag{D2} and \thetag{D3},
and the mapping from $a$ to $\gamma$ is well-defined.

To see that the mapping is one-to-one we show that we can recover $(a_{ij})$
from $\gamma$.
Let $a=(a_{ij})\in\TSSCPP{n,m}$.
Since $a$ is self-complementary,
all cells in $\BOXU{2(n+m)}{---}$ are in $a$,
no cell in $\BOXU{2(n+m)}{+++}$ is in $a$.
Since $a$ is totally symmetric,
the information in $\BOXU{2(n+m)}{+--}$ and $\BOXU{2(n+m)}{--+}$ is recovered from that of $\BOXU{2(n+m)}{-+-}$,
and the information in $\BOXU{2(n+m)}{++-}$ and $\BOXU{2(n+m)}{-++}$ is recovered from that of $\BOXU{2(n+m)}{+-+}$,
whereas the information in $\BOXU{2(n+m)}{-+-}$ is recovered from that of $\BOXU{2(n+m)}{+-+}$ from self-complementarity.
Thus the information in $\BOXU{2(n+m)}{+-+}$ is exactly what we need,
which is completely determined by $\gamma$.

Finally we show that our mapping is onto.
Assume that we have recovered $a$ from $\gamma$ in the above way using total symmetry and self-complementarity.
We have to check $a$ is a plane partition.
For example,
suppose $p=(x+n+m,y,z)\in\BOXU{2(n+m)}{+--}$ and
$p'=(x'+n+m,y',z'+n+m)\in\BOXU{2(n+m)}{+-+}$ with $x\leq x'$ and $y\leq y'$,
and $p'\in a$,
i.e. $(x',y',z')\in\gamma\Rightarrow(z',y',x')\in\gamma\Rightarrow x'+y'\leq n+m$ by (D2) and (D3).
Thus we have $(n+m+1-x)+(n+m+1-y)\geq n+m+2$,
which implies $(n+m+1-z,n+m+1-x,n+m+1-y)\not\in\gamma\Rightarrow(2(n+m)+1-z,n+m+1-x,2(n+m)+1-y)\not\in a\Rightarrow p=(x+n+m,y,z)\in a$.
The proofs in the other cases are similar.
This shows $a$ is a plane partition.
Next we show that $a$ satisfies the condition (T).
Since $\gamma$ has at most $n$ columns which implies $(1,n+1,1)\not\in\gamma$,
i.e. $(n+m+1,n+1,n+m+1)\not\in a$.
From self-complementarity and total symmetry we have $(n+m,n+m,n+2m)\in a$,
which implies $\BOXCU{2m}{n+m}{--+}\subseteq a$.
Similarly we have $\BOXCU{2m}{n+m}{-+-}\subseteq a$ and $\BOXCU{2m}{n+m}{+--}\subseteq a$,
and we conclude that $\BOXCU{2m}{n+m}{-}\subseteq a$.

Finally we use the bijection between self-conjugate partitions and strict partitions,
i.e. convert each column of $\gamma$ into a strict partition.
Thus the map $(a_{ij})_{1\leq i,j\leq2(n+m)}\mapsto(\max\{a_{i+n+m,j}-(n+m)-i+1,0\})_{1\leq i,j\leq n+m}$
defines a bijection of $\TSSCPP{n,m}$ onto $\CSPP{n,m}$.
This completes the proof.
\end{demo}

\begin{demo}{Proof of Corollary~\ref{cor:bijection_TSPPtoCSPP}}
Combining the bijections in Theorem~\ref{thm:bijection_TSPP} and Theorem~\ref{thm:bijection_CSPP},
we have a bijection between $\TSPP{n,m}$ and $\CSPP{n,m}$,
which is directly computed using \thetag{\ref{eq:set_TSPP}} and \thetag{\ref{eq:set_CSPP}}.
Thus \thetag{\ref{eq:set_CSPPtoTSPP}} and \thetag{\ref{eq:set_TSPPtoCSPP}} gives a bijection between $\TSPP{n,m}$ and $\CSPP{n,m}$.
Using these identities,
we have
$
n-b_{ij}\geq k
\Leftrightarrow
b_{ij}\leq n-k
\Leftrightarrow
(i,j,n+1-k)\not\in b
\Leftrightarrow
(n+m-j,k,1+j-i)\in c
\Leftrightarrow
c_{n+m-j,k}\geq1+j-i
$
for $1\leq i\leq j\leq n+m-1$,
which implies
\begin{equation*}
n-b_{ij}=\sharp\{k:c_{n+m-j,k}\geq1-i+j\}
\end{equation*}
for $1\leq i\leq j\leq n+m-1$.
This proves \thetag{\ref{eq:CSPPtoTSPP}}.
\end{demo}

\section{
The statistics
}\label{sec:statistics}

In \cite[pp.282]{MRR2} Mills, Robbins and Rumsey have defined the statistics $U_{r}$
(see \thetag{\ref{eq:stat_mrr}})
for $\TSPP{n}$ which have been conjectured 
to have the same distribution as the position of the $1$ in the top row of an alternating sign matrix.
(See \cite{MRR1} for detailed explanation about the alternating sign matrices).
In this section we generalize this $U_{r}$ as the statistics for the generalized set $\TSPP{n,m}$
(see \thetag{\ref{eq:stat}}),
and the main goal is to translate $U_{r}$ into the statistics for $\CSPP{n,m}$
(see Theorem~\ref{th:U_k})
using the bijection in Corollary~\ref{cor:bijection_TSPPtoCSPP}.
We also define new statistics $V^\text{R}$ and $V^\text{C}$
which is not in \cite{MRR2},
and give a new conjecture (see Conjecture~\ref{conj:even_row}).

Throughout this paper,
for $b=(b_{ij})_{1\leq i\leq j\leq n+m-1}\in\TSPP{n,m}$,
we set $b_{i,n+m}=n-i$ for all $i$ and $b_{0,j}=n$ for all $j$ by convention.
As an extension of these statistics
we define the following statistics.

For a $b=(b_{ij})_{1\leq i\leq j\leq n+m-1}$ in $\TSPP{n,m}$ and integers $r=1,\dots,n+m$,
let
\begin{equation}
U_{r}(b)=\sum_{t=1}^{n+m-r}(b_{t,t+r-1}-b_{t,t+r})
+\sum_{t=n+m-r+1}^{n+m-1}\chi\{b_{t,n+m-1}>n-t\}.
\label{eq:stat}
\end{equation}
This $U_{r}(b)$ agrees with \thetag{\ref{eq:stat_mrr}}
when $m=0$.
It is easy to check that each of these functions $U_{r}$ can vary between $0$ and $n+m-1$ as $b$ varies over $\TSPP{n,m}$.
We put $\overline U_{r}(b)=n+m-1-U_{r}(b)$.
Furthermore
we identify each element in $\CSPP{n,m}$ and each element in $\TSPP{n,m}$
by the bijection $\BIJ{n,m}$ defined in Corollary~\ref{cor:bijection_TSPPtoCSPP},
and we define $U_{r}(c)=U_{r}(\BIJ{n,m}(c))$ and $\overline U_{r}(c)=\overline U_{r}(\BIJ{n,m}(c))$
for $c\in\CSPP{n,m}$.
The following theorem enable us to compute $\overline U_{r}(c)$ directly.
\begin{theorem}
\label{th:U_k}
Let $m$ and $n\geq1$ be nonnegative integers and let $c\in\CSPP{n,m}$.
Then
$\overline U_{r}(c)$ is the number of parts equal to $r$
plus the number of saturated parts less than $r$,
i.e.
\begin{align}
\overline U_{r}(c)=\sharp\{(i,j):c_{ij}=r\}+\sharp\{1\leq k<r:c_{1,n+m-k}=k\}.
\label{eq:statUk}
\end{align}
Especially $\overline U_{1}(c)$ is the number of $1$'s in $c$ 
and $\overline U_{n+m}(c)$ is the number of saturated parts in $c$.
It is also easy to see that $\overline U_{n+m-1}(c)=\overline U_{n+m}(c)$
since, if a part of $c\in\CSPP{n,m}$ is equal to $n+m-1$, then it is saturated.
\end{theorem}
We also define two new statistics.
For $c\in\CSPP{n,m}$,
let $V^\text{R}(c)$ denote the number of rows of $c$ of odd length,
and let $V^\text{C}(c)$ denote the number of columns of $c$ of odd length.
For example,
$\CSPP{3}$ consists of the following 7 elements:
\begin{equation*}
\begin{array}{cc}
\emptyset\\
&
\end{array}
,\qquad
\begin{array}{cc}
1\\
&
\end{array}
,\quad
\begin{array}{cc}
1&\pmb{1}\\
&
\end{array}
,\qquad
\begin{array}{cc}
\pmb{2}\\
&
\end{array}
,\quad
\begin{array}{cc}
\pmb{2}&\pmb{1}\\
&
\end{array}
,\qquad
\begin{array}{cc}
\pmb{2}\\
1
\end{array}
,\quad
\begin{array}{cc}
\pmb{2}&\pmb{1}\\
1
\end{array}.
\end{equation*}
The distribution statistics for $\overline U_{k}(c)$, $k=1,2,3$,
$V^\text{R}(c)$ and $V^\text{C}(c)$ in $\CSPP{3}$
are as in Table~\ref{table:statistics}
where $c\in\CSPP{3}$ are in this order.
\begin{table}[htb]
\begin{center}
\begin{tabular}{|c||c|c|c|c|c|c|c|}
\hline
$\overline U_{1}(c)$ & \quad$0$\quad & \quad$1$\quad &\quad $2$\quad & \quad$0$\quad & \quad$1$\quad & \quad$1$\quad & \quad$2$\quad \\
\hline
$\overline U_{2}(c)$ & \quad$0$\quad & \quad$0$\quad & \quad$1$\quad & \quad$2$\quad & \quad$1$\quad & \quad$1$\quad & \quad$2$\quad \\
\hline
$\overline U_{3}(c)$ & \quad$0$\quad & \quad$0$\quad & \quad$1$\quad & \quad$2$\quad & \quad$1$\quad & \quad$1$\quad & \quad$2$\quad \\
\hline
$V^\text{R}(c)$      & \quad$0$\quad & \quad$1$\quad & \quad$0$\quad & \quad$1$\quad & \quad$0$\quad & \quad$2$\quad & \quad$1$\quad \\
\hline
$V^\text{C}(c)$      & \quad$0$\quad & \quad$1$\quad & \quad$2$\quad & \quad$1$\quad & \quad$2$\quad & \quad$0$\quad & \quad$1$\quad \\
\hline
\end{tabular}
\end{center}
\caption{The distribution statistics table in $\CSPP{3}$}
\label{table:statistics}
\end{table}
As the reader may easily see by comparison, 
the distribution of all the $\overline U_{k}(c)$, $k=1,2,3$, in $\CSPP{3}$
is independent of $k$,
and the function $V^\text{C}(c)$ also have the same distribution.
The first fact was proved in \cite{MRR2} when $m=0$,
and we also see it for general $\CSPP{n,m}$ in Section~\ref{sec:gf}.
We also give a proof that $V^\text{C}(c)$ have the same distribution
as $\overline U_{k}(c)$ in this paper.
From this fact we also can see that the number of $c$ in $\CSPP{n}$ with all columns
of even length
is equal to $A_{n-1}$.
In fact.
as the reader can see from the above example,
the two plane partitions $\emptyset$ and
$
\begin{array}{cc}
2\\
1
\end{array}
$ 
have all columns of even length,
and its generating function $\sum_{c\in\CCSPP{3}} t^{\overline U_{k}(t)}=1+t$,
$k=1,2,3$,
is equal to $A_{2}(t)$.
Further one may check the case of $c\in\CSPP{n}$
with all rows of even length.
In the above example,
the three elements $\emptyset$,
$
\begin{array}{cc}
1&1
\end{array}
$
and
$
\begin{array}{cc}
2&1
\end{array}
$
have all rows of even length,
and its generating function is $\sum_{c\in\RCSPP{3}} t^{\overline U_{k}(t)}=1+t+t^2$,
$k=1,2,3$,
which coincidently equals $A^\text{VS}_{5}(t)$.
It is also easily checked $\sum_{c\in\RCSPP{n}} t^{\overline U_{k}(t)}$ becomes $3(1+t+t^2)$ if $n=4$,
$3(3+6t+8t^2+6t^3+3t^4)$ if $n=5$,
$26(3+6t+8t^2+6t^3+3t^4)$ if $n=6$,
and so on.
From this numerical experiment
one may expect the following conjecture could hold:
\begin{conjecture}
\label{conj:even_row}
Let $n\geq 1$ be a positive integer,
and let $1\leq r\leq n$.
Then 
\[
\sum_{c\in\RCSPP{n}}t^{\overline U_{r}(c)}=
\begin{cases}
A^\text{VS}_{2m+1}\cdot A^\text{VS}_{2m+1}(t)
&\text{ if $n=2m$,}\\
A^\text{VS}_{2m+1}\cdot A^\text{VS}_{2m+3}(t)
&\text{ if $n=2m+1$.}
\end{cases}
\]
would hold.
Especially,
if we put $t=1$,
the number of $c$ in $\RCSPP{n}$ would be
\[
\begin{cases}
\left(A^\text{VS}_{2m+1}\right)^2
&\text{ if $n=2m$,}\\
A^\text{VS}_{2m+1}\cdot A^\text{VS}_{2m+3}
&\text{ if $n=2m+1$.}
\end{cases}
\]
\end{conjecture}
Di~Francesco defined another weight in \cite{F1} for NILPs
(non-intersecting lattice paths)
which seems to have the same distribution as the position of the $1$ in the top row of an alternating sign matrix.
\begin{demo}{Proof of Theorem~\ref{th:U_k}}
Using the identity \thetag{\ref{eq:CSPPtoTSPP}} and $b_{n+m-r,n+m}=r-m$,
we have
\begin{align*}
\sum_{t=1}^{n+m-r}(b_{t,t+r-1}-b_{t,t+r})
&=\sum_{t=1}^{n+m-r-1}\sharp\{k:c_{t,k}\geq r+1\}\\
&-\sum_{t=1}^{n+m-r}\sharp\{k:c_{t,k}\geq r\}
+n+m-r.
\end{align*}
Using $c_{t,k}\leq c_{1,k}-t+1\leq n+m-k-t+1$,
we have $\sharp\{k:c_{n+m-r,k}\geq r+1\}=0$.
Thus the right-hand-side is equal to
$n+m-r-\sum_{t=1}^{n+m-r}\sharp\{k:c_{t,k}=r\}$.
Note that $c_{n+m-r+1,1}\leq r-1$,
which implies $r$ does not appear in the $(n+m-r+1)$st row,
and $\sum_{t=1}^{n+m-r}\sharp\{k:c_{t,k}=r\}$ is the number of $r$ appearing in $c$.
Also by \thetag{\ref{eq:CSPPtoTSPP}},
it is easy to see that $b_{t,n+m-1}>n-t$ if and only if $c_{1,t}<n+m-t$.
Thus we have
\[
\sum_{t=n+m-r+1}^{n+m-1}\{b_{t,n+m-1}>n-t\}=r-1-\sum_{t=1}^{r-1}\chi\{c_{1,n+m-t}=t\}
\]
and this shows
$
\overline U_{r}(c)=\sum_{t=1}^{n+m-r}\sharp\{k:c_{t,k}=r\}+\sum_{t=1}^{r-1}\chi\{c_{1,n+m-t}=t\}.
$
\end{demo}

\section{
The monotone triangles and TSSCPPs
}\label{sec:mt}


In \cite[pp.287]{MRR2}
a subset $\TSPP{n}^{k}$ of $\TSPP{n}$ is defined.
In this section we 
show that we can generalize this definition naturally to the subset $\TSPP{n,m}^{k}$ of $\TSPP{n,m}$,
and also show that there is a very nice interpretation in the words of $\CSPP{n,m}$
through the bijection defined in Section~\ref{sec:bijections}
(see Theorem~\ref{prop:TSPPtoCSPP_k}).
As a result we can present a new conjecture (see Conjecture~\ref{conj:refined_MT})
which is a refined version of Conjecture~\ref{conj:MT}
(\cite[Conjecture~7]{MRR2}).
We also try to restate \cite[Conjecture~7']{MRR2} by means of $\CSPP{n,m}$
(see Theorem~\ref{prop:TSPPtoCSPP_kr}).

For $k=0,\dots,n+m-1$,
let $\TSPP{n,m}^{k}$ denote the subset of those $b=(b_{ij})_{1\leq i\leq j\leq n+m-1}$ in $\TSPP{n,m}$
such that all $b_{ij}$ in the first $n+m-1-k$ columns are equal to
their maximum values $n$.
Also, for $k=0,\dots,n+m-1$,
let $\CSPP{n,m}^{k}$ denote the subset of those $c=(c_{ij})$ in $\CSPP{n,m}$
which has at most $k$ rows.
For example,
if $n=3$ and $m=0$,
$\TSPP{3}$ consists of the following seven elements:
\begin{equation*}
\smyoung{
3&3\\
\blank&3\\}
\qquad
\smyoung{
3&3\\
\blank&2\\
}
\qquad
\smyoung{
3&3\\
\blank&1\\
}
\qquad
\smyoung{
3&2\\
\blank&2\\
}
\qquad
\smyoung{
3&2\\
\blank&1\\
}
\qquad
\smyoung{
2&2\\
\blank&2\\
}
\qquad
\smyoung{
2&2\\
\blank&1\\
}
\end{equation*}
whereas
$\CSPP{3}$ consists of the following seven plane partitions.
\begin{equation*}
\begin{array}{cc}
\emptyset\\
&
\end{array}
\qquad
\begin{array}{cc}
1\\
&
\end{array}
\qquad
\begin{array}{cc}
1&\pmb{1}\\
&
\end{array}
\qquad
\begin{array}{cc}
\pmb{2}\\
&
\end{array}
\qquad
\begin{array}{cc}
\pmb{2}&\pmb{1}\\
&
\end{array}
\qquad
\begin{array}{cc}
\pmb{2}\\
1&
\end{array}
\qquad
\begin{array}{cc}
\pmb{2}&\pmb{1}\\
1&
\end{array}
\end{equation*}
There are only one element,
i.e. $\emptyset$,
of $\CSPP{3}$ with no row,
five elements of $\CSPP{3}$ with with at most one row,
and seven elements of $\CSPP{3}$ with at most two rows.

\begin{theorem}
\label{prop:TSPPtoCSPP_k}
Let $m$ and $n\geq1$ be nonnegative integers.
Let $0\leq k\leq n+m-1$.
By the bijection $\BIJ{n,m}$ defined in Corollary~\ref{cor:bijection_TSPPtoCSPP},
the subset $\TSPP{n,m}^{k}$ of $\TSPP{n,m}$ is in one-to-one correspondence with the subset 
$\CSPP{n,m}^{k}$
of $\CSPP{n,m}$.
Especially,
we have 
\begin{equation*}
\sum_{b\in\TSPP{n,m}^k}t^{U_r(b)}
=\sum_{c\in\CSPP{n,m}^k}t^{\overline U_r(c)}.
\end{equation*}
\end{theorem}
We will prove Theorem~\ref{prop:TSPPtoCSPP_kr}
which is a refined version of this proposition
and immediately implies this theorem.

By Theorem~\ref{prop:TSPPtoCSPP_k},
we can reduce Conjecture~\ref{conj:MT} (\cite[pp.287, Conjecture~7]{MRR2})
to the enumeration problem of $\CSPP{n}^k$.
Furthermore,
we can present the following new conjecture
(Conjecture~\ref{conj:refined_MT}) which is a refined version of Conjecture~\ref{conj:MT}.
We give the weight to an element $m=(m_{ij})_{1\leq i\leq j\leq n}$ of $\MT{n}^{k}$
by the vale $m_{n,n}-1$.
Set $M_{n}^{k}(t)$ to be the polynomial
\begin{equation}
\sum_{m=(m_{ij})\in\MT{n}^{k}}t^{m_{n,n}-1}.
\end{equation}
For instance,
from the above example,
one easily sees
$M_{3}^{0}(t)=1$,
$M_{3}^{1}(t)=t^2+2t+2$,
and $M_{3}^{2}(t)=2t^2+3t+2$.
The reader who is familiar with the alternating sign matrices
may notice that this weight corresponds to the position of the 1 
in the top row of an alternating sign matrix.
\begin{conjecture}
\label{conj:refined_MT}
Let $n\geq1$ and $1\leq r\leq n$.
Then, for $k=0,1,\dots,n-1$,
we would have
\begin{equation}
\sum_{c\in\CSPP{n}^k}t^{\overline U_r(c)}=M_{n}^{k}(t).
\end{equation}
Later we will see that the left-hand side does not depend on $r$.
\end{conjecture}
More generally 
the following set theoretic partitions of $\MT{n}^{k}$ and $\TSPP{n}^{k}$
are defined in \cite{MRR2}.
Suppose $n\geq1$, $0\leq k\leq n-1$ and that $x$ and $y$ are nonnegative integers.
Let $\MT{n}^{k,x,y}$ denote the subset of all those $m=(m_{ij})_{1\leq i\leq j\leq n}$ in $\MT{n}^{k}$
such that there are precisely $x+1$ parts in column $n-k+1$ which are equal to their minimum possible value $n-k-i+2$
and $y+1$ parts in column $n$ equal to their maximum values $n$.
From the above example
one sees that
the following two elements of $\MT{3}^{2}$
has precisely $2$ parts in the second column 
which are equal to their minimum possible values
and exactly two parts in the last column which equal to their maximum values $3$
so that $\MT{3}^{2,1,1}$ is composed of them:
\begin{equation*}
\begin{array}{ccc}
 & &1\\
 &1&3\\
1&2&3
\end{array}
\quad
\begin{array}{ccc}
 & &2\\
 &1&3\\
1&2&3
\end{array}
\end{equation*}
A similar argument shows that
$\sharp\MT{3}^{1,1,1}=\sharp\MT{3}^{1,2,0}=\sharp\MT{3}^{0,0,0}
=\sharp\MT{3}^{2,0,1}=\sharp\MT{3}^{1,0,1}=\sharp\MT{3}^{1,1,0}
=\sharp\MT{3}^{2,0,2}=\sharp\MT{3}^{1,0,2}=\sharp\MT{3}^{2,1,2}=1$,
$\sharp\MT{3}^{2,1,0}=\sharp\MT{3}^{2,1,1}=2$
and the others has the cardinality zero.

Similarly,
suppose $n\geq1$, $m\geq0$, $0\leq k\leq n+m-1$ and that $x$ and $y$ are nonnegative integers.
Let $\TSPP{n,m}^{k,x,y}$ be the subset of all those $b$ in $\TSPP{n,m}^{k}$
such that there are precisely $x$ entries in column $n+m-k$
which are equal to their maximum values $n$ and
there are exactly $y$ parts $b_{i,n+m-1}$ in column $n+m-1$
equal to their minimum value $\max\{n-i,0\}$.
The following conjecture is due to Mills, Robbins and Rumsey:
\begin{conjecture}
\label{conj:MT2}
(\cite[pp.291, Conjecture~7']{MRR2})
The cardinality of $\TSPP{n}^{krs}$ is equal to 
the cardinality of $\MT{n}^{krs}$.
\end{conjecture}

We can restate this conjecture by the following theorem.
Let $\CSPP{n,m}^{k,x,y}$ be the subset of all those $c$ in $\CSPP{n,m}^{k}$
such that the $k$th row of $c$ has exactly $n+m-k-x$ parts
 and
there are exactly $y$ saturated parts in the first row of $c$.
For example,
if $n=3$, $m=0$, $k=2$ and $x=y=1$,
then the following two elements of $\CSPP{3}^{2}$ has no parts in the second row
and has precisely one saturated part.
\begin{equation*}
\begin{array}{cc}
1&\pmb{1}
\end{array}
\qquad
\begin{array}{cc}
\pmb{2}
\end{array}
\end{equation*}
\begin{theorem}
\label{prop:TSPPtoCSPP_kr}
Let $m$, $n\geq1$, $x$ and $y$ be nonnegative integers.
Let $0\leq k\leq n+m-1$.
By the bijection $\BIJ{n,m}$ defined in Corollary~\ref{cor:bijection_TSPPtoCSPP},
the subset $\TSPP{n,m}^{k,x,y}$ of $\TSPP{n,m}$ is in one-to-one correspondence with the subset 
$\CSPP{n,m}^{k,m,n}$
of $\CSPP{n,m}$.
\end{theorem}

\section{
Refined strange enumeration
}\label{sec:se}

In \cite{E},
T.~Eisenk\"olbl gave the $(-1)$-enumeration of $\TSSCPP{n}$.
In \cite{Ste2},
J.~Stembridge proposed more general ``strange enumeration''
of CSPPs (cyclically symmetric plane partitions).
The aim of this section is to restate the ``strange enumeration''
of $\TSSCPP{n,m}$
in the words of $\CSPP{n,m}$
(see Theorem~\ref{thm:strange})
to obtain the generating functions in Section~\ref{sec:gf}.

In \cite[pp.25/26]{Ku2}
natural $(-1)$-enumeration for the six symmetry classes of plane partitions
which involve complementation
are proposed (also see \cite{E}).
Here we restrict our attention to the TSSCPP case,
i.e. $\TSSCPP{n,m}$.
In this case the symmetry group $G$ is generated by
$S_3$ acting on the coordinates and the complementation $\CMP{n+m}$.
Then each orbit of a TSSCPP under the group action of $G$ is always half-filled.
Thus there is a natural move between TSSCPPs which replaces half of an orbit of cubes by the opposite half.
Any two TSSCPPs differs by either an odd or an even number of moves,
and we can define a relative sign between them.
This sign becomes absolute if we assign the weight $1$ to the $\BOXU{n+m}{-}\in\TSSCPP{n,m}$.
We call this type of signed enumeration the \defterm{$(-1)$-enumeration} of TSSCPPs.

In the strange enumeration \cite{Ste2},
Stembridge has proposed a way of signed enumerations of CSPPs assigning each orbit $\pm1$ 
depending on which family it belongs
where each family corresponds to the faces of Coxeter complex of type $A_2$
modulo the action of $C_3$.
Here we show that we can define a similar strange enumerations for TSSCPPs.
Let $\langle(i,j,k)\rangle$ denote the $G$-orbit of a cube $(i,j,k)$.
Let
\[
\BOX{2(n+m)}/G=P_1\cup P_2\cup P_3
\]
denote the partition of the $G$-orbits of $\BOX{2(n+m)}$ into the families
where $P_1=\{\langle(i,i,i)\rangle\}$,
$P_2=\{\langle(i,j,j)\rangle: i\neq j\}$
and $P_3=\{\langle(i,j,k)\rangle: \text{ $i$, $j$ and $k$ are all distinct}\}$.
For $a\in\TSSCPP{n,m}$,
let $m_i(a)$ denote the the number of moves in $P_i$
which is needed to obtain $a$ from $\BOXU{2(n+m)}{-}$.
Choose an arbitrary weight function $w:\TSSCPP{n,m}\rightarrow R$
which assigns values in some commutative ring $R$ to each TSSCPP $a\in\TSSCPP{n,m}$.
We may consider the signed enumeration
\[
G_{n,m}(w)_{s,t,u}=\sum_{a\in\TSSCPP{n,m}}s^{m_1(a)}t^{m_2(a)}u^{m_3(a)}w(a)
\]
where 
$s$ is the sign of moves in the orbits $\langle(i,i,i)\rangle$,
$t$ is the sign of moves in the orbits $\langle(i,j,j)\rangle$ with $i\neq j$
and
$u$ is the sign of moves in the orbits $\langle(i,j,k)\rangle$ with $i$, $j$ and $k$ all distinct.
Here $s$, $t$ and $u$ are $\pm1$,
$G_{n,m}(w)_{1,1,1}$ is the ordinary generating function,
and $G_{n,m}(w)_{-1,-1,-1}$ is the $(-1)$-enumeration.
But note that $s$ has no meaning since we need no move in $P_1$ to obtain a TSSCPP
from $\BOXU{2(n+m)}{-}$,
which implies $G_{n,m}(w)_{s,t,u}$ does not depend on $s$,
i.e. $G_{n,m}(w)_{-1,t,u}=G_{n,m}(w)_{1,t,u}$.
Thus we may consider two other strange enumerations,
i.e.
$G_{n,m}(w)_{1,-1,1}$ and $G_{n,m}(w)_{1,1,-1}$.
\begin{theorem}
\label{thm:strange}
Let $m$ and $n\geq1$ be non-negative integers.
Let $a$ be a TSSCPP in $\TSSCPP{n,m}$.
Suppose $a$ is mapped to a RCSPP $c$ in $\CSPP{n,m}$ by the bijection 
defined in Theorem~\ref{thm:bijection_CSPP}.
Then $m_1(a)=0$,
$m_2(a)$ is the sum of parts in the first row of $c$
and $m_3(a)$ is the sum of parts of $c$ which is not in the first row.
Especially the total number $m_1(a)+m_2(a)+m_3(a)$ of moves equals $|c|$.
\end{theorem}
In \cite{E} Eisenk\"olbl considered $(-1)$-enumeration of several classes of plane partitions
with complementary symmetry
and obtain the result that the $(-1)$-enumeration of totally symmetric self-complementary plane partitions
contained in $\BOX{2n}$ equals
$
\begin{cases}
A_{n+2}^\text{VS}
&\text{ if $n$ is odd,}\\
0
&\text{ otherwise,}
\end{cases}
$
(see \cite[Theorem 5]{E}).
By Theorem~\ref{thm:strange},
we can consider the refined $(-1)$-enumeration 
$\sum_{c\in\CSPP{n,m}}(-1)^{|c|}\,t^{\overline U_r(c)}$
($1\leq r\leq n$)
and the doubly refined $(-1)$-enumeration 
$\sum_{c\in\CSPP{n,m}}(-1)^{|c|}\,t^{\overline U_1(c)}u^{\overline U_r(c)}$
($2\leq r\leq n$)
of  totally symmetric self-complementary plane partitions.
For $c$ in $\CSPP{n,m}$,
let $\PROFILE{c}$ denote the sum of parts in the first row of $c$.
Moreover we can propose the problem to consider
the refined strange enumeration
$\sum_{c\in\CSPP{n,m}}(-1)^{\PROFILE{c}}\,t^{\overline U_r(c)}$
($1\leq r\leq n$)
and
$\sum_{c\in\CSPP{n,m}}(-1)^{\PROFILE{c}}(-1)^{|c|}\,t^{\overline U_r(c)}$
($1\leq r\leq n$),
and the doubly refined strange enumeration
$\sum_{c\in\CSPP{n,m}}(-1)^{\PROFILE{c}}\,t^{\overline U_1(c)}u^{\overline U_r(c)}$
($2\leq r\leq n$)
and
$\sum_{c\in\CSPP{n,m}}(-1)^{\PROFILE{c}}(-1)^{|c|}\,t^{\overline U_1(c)}u^{\overline U_r(c)}$
($2\leq r\leq n$).
Even the strange enumeration $\sum_{c\in\CSPP{n,m}}(-1)^{\PROFILE{c}}$
and $\sum_{c\in\CSPP{n,m}}(-1)^{\PROFILE{c}}(-1)^{|c|}$
(i.e. $t=u=1$)
are considered by nobody.

\section{
The generating functions
}\label{sec:gf}

Using the bijection \thetag{\ref{eq:CSPPtoTSPP}} in Corollary~\ref{cor:bijection_TSPPtoCSPP},
we see 
\begin{equation}
\sharp\CSPP{n,m}
=\prod_{k=0}^{n-1}
\frac{(3k+3m+1)!\prod_{i=0}^{m}(k+2i)!}
{(2k+m)!(2k+3m+1)!\prod_{i=1}^{m}(k+2i-1)!}
\label{eq:card_Pnm}
\end{equation}
from Krattenthaler's result (see \cite[Theorem~2]{Kr1}).
In this section we give the generating function of $\CSPP{n,m}$,
$\CCSPP{n,m}$ and $\RCSPP{n,m}$
with the weight $\overline U_{k}$ defined in  \thetag{\ref{eq:statUk}},
 $V^\text{C}$ and/or $V^\text{R}$,
which give several refinements of of \thetag{\ref{eq:card_Pnm}}
and Pfaffian expressions for Conjecture~\ref{conj:refined},
Conjecture~\ref{conj:double_refined},
Conjecture~\ref{conj:MT} and Conjecture~\ref{conj:even_row}.
First we give a general theorem (see Theorem~\ref{thm:gf})
and derive the Pfaffian expression for each conjecture
as a corollary of the general theorem.
The reader should see 
Corollary~\ref{cor:refined} for Conjecture~\ref{conj:refined} and Conjecture~\ref{conj:even_row},
Corollary~\ref{cor:doubly_refined} for Conjecture~\ref{conj:double_refined},
Corollary~\ref{cor:MT} for Conjecture~\ref{conj:MT}
and Corollary~\ref{cor:(-1)enumeration} for $(-1)$-enumeration of $\CSPP{n,m}$.
We also give two new conjectures,
i.e. Conjecture~\ref{conj:row} and \thetag{\ref{eq:conj-1}} \thetag{\ref{eq:conj-1-2}} in this section.

\bigbreak

First of all we fix the notation.
The reader can consult Macdonald's
book \cite{M} for information on the symmetric functions. 
For a positive integer $l$,
let $\pmb{x}=(x_{1},x_{2},\dots,x_{l})$ be an $l$-tuple of variables.
We write 
the $r$th elementary symmetric function in $x_{1},\dots,x_{l}$
as $e^{(l)}_{r}(\pmb{x})=e^{(l)}_{r}(x_{1},\dots,x_{l})$,
i.e.
\[
\sum_{r=0}^{\infty}e^{(l)}_{r}(\pmb{x})y^r=\prod_{i=1}^{l}(1+x_{i}y).
\]
Let $m$ and $n\geq1$ be non-negative integers.
Let $\pmb{x}=(x_1,\dots,x_{n+m})$ and $\pmb{t}=(t_1,\dots,t_{n+m})$ be commutative variables.
We write $T_i$ for $\prod_{k=i}^{n+m}t_k$.
For each plane partition $c$ in $\CSPP{n,m}$,
we assign the weight
\[
\pmb{t}^{\overline U(c)}\pmb{x}^c
=
t_1^{\overline U_1(c)}t_2^{\overline U_2(c)}\cdots t_{n+m}^{\overline U_{n+m}(c)}
x_1^{\mu_1}x_2^{\mu_2}\cdots x_{n+m}^{\mu_{n+m}},
\]
where $\mu_i$ is the number of $i$'s in $c$ for $i=1,\dots,n+m$.
Notice that $i$ appears at most $n+m-i$ times in $c$.
As an application of this Gessel-Viennot formula \cite{GV},
we obtain the following fundamental lemma
to enumerate the elements of $\CSPP{n,m}$:
\begin{lemma}
\label{lem:lp}
Let $m$ and $n\geq1$ be non-negative integers,
and put $N=n+m$.
Let $\lambda$ be a partition with $\ell(\lambda)\leq n$.
Then
the generating function of all plane partitions $c\in\CSPP{n,m}$
of shape $\lambda'$ with the weight $\pmb{t}^{\overline U(c)}\pmb{x}^c$ is given by
\begin{equation}
\sum_{{c\in\CSPP{n,m}}\atop{\SHAPE{c}=\lambda'}}\pmb{t}^{\overline U(c)}\pmb{x}^c
=\det\left(e^{(N-i)}_{\lambda_{j}-j+i}(t_{1}x_{1},\dots,t_{N-i-1}x_{N-i-1},T_{N-i}x_{N-i})\right)_{1\leq i,j\leq n},
\label{eq:GenFunc}
\end{equation}
where $T_{i}=\prod_{k=i}^{N}t_k$.
\end{lemma}
As a consequence of Lemma~\ref{lem:lp} and Theorem~\ref{msf},
we obtain the following fundamental theorem
from which we can derive all the results as corollaries.
\begin{theorem}
\label{thm:gf}
Let $m$ and $n\geq1$ be non-negative integers.
Let $N$ be an even integer such that $N\geq n+m-1$.
Let $A$ be an 
$(n+N)\times(n+N)$ skew-symmetric matrix.
Let 
$B_{n,m}^{N}(\pmb{t},\pmb{x})
=(b^{(m)}_{ij}(\pmb{t},\pmb{x}))_{0\leq i\leq n-1,\ 0\leq j\leq n+N-1}$
be the rectangular matrix whose $(i,j)$th entry is
\begin{equation}
b^{(m)}_{ij}(\pmb{t},\pmb{x})
=e^{(i+m)}_{j-i}(t_{1}x_{1},\dots,t_{i+m-1}x_{i+m-1},T_{i+m}x_{i+m}),
\label{eq:entry_b}
\end{equation}
where $T_{i}=\prod_{k=i}^{n+m}t_k$.
Then
the generating function for all plane partitions $c\in\CSPP{n,m}$
 with the weight 
$
(-1)^{|\SHAPE{c}|}
\Pf\left(A^{\overline{I_{n}(\SHAPE{c}')}}_{\overline{I_{n}(\SHAPE{c}')}}\right)
\pmb{t}^{\overline U(c)}\pmb{x}^c
$ is given by
\begin{equation}
\sum_{{\lambda}\atop{\ell(\lambda)\leq n}}
\sum_{{c\in\CSPP{n,m}}\atop{\SHAPE{c}=\lambda'}}
(-1)^{|\lambda|}
\Pf\left(A^{\overline{I_{n}(\lambda)}}_{\overline{I_{n}(\lambda)}}\right)\pmb{t}^{\overline U(c)}\pmb{x}^c
=\Pf\begin{pmatrix}
O_{n}&J_{n}B_{n,m}^{N}(\pmb{t},\pmb{x})\\
-{}^t\!B_{n,m}^{N}(\pmb{t},\pmb{x})J_{n}& A
\end{pmatrix}.
\label{eq:General}
\end{equation}
Here $\overline I=[n+N]\setminus I$ stands for the complement in the set $[n+N]$.
\end{theorem}
To specialize \thetag{\ref{eq:General}},
we use the following notation.
Let $n$ and $N$ be positive integers,
and let $m$ be a nonnegative integer.
Let $B_{n,m}^{N}(t,u)=(b^{(m)}_{ij}(t,u))_{0\leq i\leq n-1,\ 0\leq j\leq n+N-1}$
be the $n\times(n+N)$ matrix whose $(i,j)$th entry is
\begin{equation}
b^{(m)}_{ij}(t,u)=
\begin{cases}
\delta_{0,j}
&\text{ if $i+m=0$,}\\
\binom{i+m-1}{j-i}+\binom{i+m-1}{j-i-1}tu
&\text{ if $i+m=1$,}\\
\binom{i+m-2}{j-i}+\binom{i+m-2}{j-i-1}(t+u)+\binom{i+m-2}{j-i-2}tu
&\text{ otherwise.}
\end{cases}
\label{eq:entry_doubly_refined}
\end{equation}
For example,
\[
B_{3,0}^{2}(t,u)
=\begin{pmatrix}
1&0&0&0&0\\
0&1&tu&0&0\\
0&0&1&t+u&tu
\end{pmatrix}.
\]
We define the $n\times (n+N)$ matrices $B_{n,m}^{N}(t)=B_{n,m}^{N}(t,1)$ and $B_{n,m}^{N}=B_{n,m}^{N}(1)$.
Then the $(i,j)$th entry of $B_{n,m}^{N}(t)$ is
\begin{equation}
b^{(m)}_{ij}(t)=
\begin{cases}
\delta_{0,j}
&\text{ if $i+m=0$,}\\
\binom{i+m-1}{j-i}+\binom{i+m-1}{j-i-1}t
&\text{ otherwise.}
\end{cases}
\label{eq:entry_refined}
\end{equation}
and the $(i,j)$th entry of $B_{n,m}^{N}$ is $\binom{i+m}{j-i}$
where the row index runs $0\leq i\leq n-1$
and the column index runs $0\leq j\leq n+N-1$.
When $m=0$,
these $B_{n,m}^{N}(t,u)$, $B_{n,m}^{N}(t)$ and $B_{n,m}^{N}$
agree with $B_{n}^{N}(t,u)$, $B_{n}^{N}(t)$ and $B_{n}^{N}$ introduced in Section~\ref{sec:intro}.
The following corollary (i) gives a Pfaffian expression for
the doubly refined TSSCPP conjecture (Conjecture~\ref{conj:double_refined}).
\begin{corollary}
\label{cor:doubly_refined}
Let $m$ and $n\geq1$ be non-negative integers,
and let $N$ be an even integer such that $N\geq n+m-1$.
\begin{enumerate}
\item[(i)]
If $r$ is a positive integer such that $2\leq r\leq n+m$,
then
the generating function for all plane partitions $c\in\CSPP{n,m}$
 with the weight $t^{\overline U_1(c)}u^{\overline U_r(c)}$ is
\begin{align}
\sum_{c\in\CSPP{n,m}}t^{\overline U_1(c)}u^{\overline U_r(c)}
&=\Pf\begin{pmatrix}
O_{n}&J_{n}B_{n,m}^{N}(t,u)\\
-{}^t\!B_{n,m}^{N}(t,u)J_{n}&\bar S_{n+N}
\end{pmatrix}.
\label{eq:gen_double_refined}
\end{align}
\item[(ii)]
If $r$ is a positive integer such that $2\leq r\leq n+m$,
then
the generating function for all plane partitions $c\in\CCSPP{n,m}$
 with the weight $t^{\overline U_1(c)}u^{\overline U_r(c)}$ is
\begin{align}
\sum_{c\in\CCSPP{n,m}}t^{\overline U_1(c)}u^{\overline U_r(c)}
&=\Pf\begin{pmatrix}
O_{n}&J_{n}B_{n,m}^{N}(t,u)\\
-{}^t\!B_{n,m}^{N}(t,u)J_{n}&\bar R_{n+N}
\end{pmatrix}.
\label{eq:gen_double_refined2}
\end{align}
\item[(iii)]
If $r$ is a positive integer such that $1\leq r\leq n+m$,
then
the generating function for all plane partitions $c\in\CSPP{n,m}$
 with the weight $t^{\overline U_r(c)}u^{V^\text{C}(c)}$ is
\begin{align}
\sum_{c\in\CSPP{n,m}}t^{\overline U_r(c)}u^{V^\text{C}(c)}
&=\Pf\begin{pmatrix}
O_{n}&J_{n}B_{n,m}^{N}(t)\\
-{}^t\!B_{n,m}^{N}(t)J_{n}&\bar C_{n+N}(u)
\end{pmatrix}.
\label{eq:gen_double_refined3}
\end{align}
\end{enumerate}
\end{corollary}
From Corollary~\ref{cor:doubly_refined}(i)(ii)
we obtain the following corollary:
\begin{corollary}
\label{cor:doubly_refined_prop}
Let $m$ and $n\geq1$ be non-negative integers.
Let $r$ and $s$ be integers such that $2\leq r,s\leq n$,
and
let $k$ be an integer such that $1\leq k\leq n$.
Then
we have
\begin{align*}
\sum_{c\in\CSPP{n,m}} t^{\overline U_{1}(c)}u^{\overline U_{r}(c)}
=\sum_{c\in\CCSPP{n,m+1}} t^{\overline U_{1}(c)}u^{\overline U_{s}(c)}
=\sum_{c\in\CSPP{n,m}} t^{\overline U_{k}(c)}u^{V^\text{C}(c)}
\end{align*}
\end{corollary}
If we put $u=1$ in Corollary~\ref{cor:doubly_refined}(i),
then we immediately obtain
the following (i),
which gives a progress to Conjecture~\ref{conj:refined}.
Also the following (iii) gives a Pfaffian expression for Conjecture~\ref{conj:even_row}.
\begin{corollary}
\label{cor:refined}
Let $m$ and $n\geq1$ be non-negative integers,
and let $N$ be an even integer such that $N\geq n+m-1$.
\begin{enumerate}
\item[(i)]
If $r$ is a positive integer such that $1\leq r\leq n+m$,
then the generating function for all plane partitions $c\in\CSPP{n,m}$
 with the weight $t^{\overline U_r(c)}$ is given by
\begin{align}
\sum_{c\in\CSPP{n,m}}t^{\overline U_r(c)}
&=\Pf\begin{pmatrix}
O_{n}&J_{n}B_{n,m}^{N}(t)\\
-{}^t\!B_{n,m}^{N}(t)J_{n}&\bar S_{n+N}
\end{pmatrix}.
\label{eq:gen_refined}
\end{align}
\item[(ii)]
If $r$ is a positive integer such that $1\leq r\leq n+m$,
then the generating function for all plane partitions $c\in\CCSPP{n,m}$
with all columns of even length
 with the weight $t^{\overline U_r(c)}$ is given by
\begin{align}
\sum_{c\in\CCSPP{n,m}}t^{\overline U_r(c)}
&=\Pf\begin{pmatrix}
O_{n}&J_{n}B_{n,m}^{N}(t)\\
-{}^t\!B_{n,m}^{N}(t)J_{n}&\bar R_{n+N}
\end{pmatrix}.
\label{eq:column_even}
\end{align}
\item[(iii)]
If $r$ is a positive integer such that $1\leq r\leq n+m$,
then the generating function for $c\in\RCSPP{n,m}$ with all rows of even length
 with the weight $t^{\overline U_r(c)}$ is given by
\begin{align}
\sum_{c\in\RCSPP{n,m}}t^{\overline U_r(c)}
&=\Pf\begin{pmatrix}
O_{n}&J_{n} B_{n,m}^{N}(t)\\
-{}^t\!B_{n,m}^{N}(t) J_{n}&\bar C_{n+N}
\end{pmatrix}.
\label{eq:row_even}
\end{align}
\item[(iv)]
The generating function for all plane partitions $c\in\CSPP{n,m}$
 with the weight $t^{V^\text{C}(c)}$ is given by
\begin{align}
\sum_{c\in\CSPP{n,m}}t^{V^\text{C}(c)}
&=\Pf\begin{pmatrix}
O_{n}&J_{n}B_{n,m}^{N}\\
-{}^t\!B_{n,m}^{N}J_{n}&\bar R_{n+N}(t)
\end{pmatrix}.
\label{eq:odd_column}
\end{align}
\item[(v)]
The generating function for all plane partitions $c\in\CSPP{n,m}$
 with the weight $t^{V^\text{R}(c)}$ is given by
\begin{align}
\sum_{c\in\CSPP{n,m}}t^{V^\text{R}(c)}
&=\Pf\begin{pmatrix}
O_{n}&J_{n}B_{n,m}^{N}\\
-{}^t\!B_{n,m}^{N}J_{n}&\bar C_{n+N}(t)
\end{pmatrix}.
\label{eq:odd_row}
\end{align}
\end{enumerate}
\end{corollary}
Note that the right-hand sides of 
\thetag{\ref{eq:gen_refined}}, \thetag{\ref{eq:column_even}} and \thetag{\ref{eq:row_even}}
do not depend on $r$.
The following corollary can be proven from 
\thetag{\ref{eq:gen_refined}},
\thetag{\ref{eq:column_even}}
and 
\thetag{\ref{eq:odd_column}},
or is a direct consequence of Corollary~\ref{cor:doubly_refined_prop}.
\begin{corollary}
\label{cor:column_even}
Let $m\geq0$ and $n\geq1$ be non-negative integers.
Let $r$ and $s$ be integers such that $1\leq r,s\leq n$,
Then
we have
\begin{align*}
\sum_{c\in\CSPP{n,m}} t^{\overline U_{r}(c)}
=\sum_{c\in\CCSPP{n,m+1}}t^{\overline U_{s}(c)}
=\sum_{c\in\CSPP{n,m}} t^{V^\text{C}(c)}.
\end{align*}
In particular we have $\sharp\CSPP{n,m}=\sharp\CCSPP{n,m+1}$.
\end{corollary}

Corollary~\ref{cor:column_even} reveal the relation between \thetag{\ref{eq:gen_refined}},
\thetag{\ref{eq:column_even}} and \thetag{\ref{eq:odd_column}}.
These are all related to the refined TSSCPP conjecture (Conjecture~\ref{conj:refined}) when $m=0$.

The identity
\thetag{\ref{eq:row_even}}
gives a Pfaffian expression for Conjecture~\ref{conj:even_row} when $m=0$.
It is not so easy to guess the explicit form of the polynomial
$
\sum_{c\in\RCSPP{n,m}} t^{\overline U_{r}(c)}
$
even for small $m\geq2$.
But,
if we put $t=1$,
then $\sharp\RCSPP{n,m}$ can be easily guessed as follows for small $m$.
\begin{conjecture}
\label{conj:row}
Let $n\geq1$, $r=0,1$ and $m$ be nonnegative integers.
Let $f(n,m)$ denote
\[
\frac{
\left(6n+6\left\lfloor\frac{m}2\right\rfloor+4\right)!
\left(6n+6\left\lceil\frac{m}2\right\rceil+4\right)!
(2n+1)!
\left(2n+2\left\lceil\frac{m}2\right\rceil\right)!
(2n+2m+1)!
\left(n+\left\lfloor\frac{m}2+1\right\rfloor\right)!
}{
(4n+m+1)!(4n+m+3)!(4n+3m+2)!(4n+3m+4)!
\left(2n+2\left\lceil\frac{m}2\right\rceil+1\right)!
\left(n+\left\lfloor\frac{m}2\right\rfloor\right)!
}.
\]
Then the number of elements $c$ in $\RCSPP{2n+r,m}$ would be
\begin{equation}
2^{-n}\frac{g(n,m+r)}{g(0,m+r)}\prod_{k=0}^{n-1}f(k,m+r)
\end{equation}
where
\[
g(n,m)=\begin{cases}
h_{m}(n)
&\text{ if $\REM{m}{4}=0$ or $1$,}\\
(4n+2m+1)h_{m}(n)
&\text{ if $\REM{m}{4}=2$ or $3$,}
\end{cases}
\]
and $h_{m}(n)$ is a polynomial of degree $2\left\lfloor\frac{m}4\right\rfloor$ in the variable $n$.
\end{conjecture}
For small $m$, $h_{0}(n)=h_{1}(n)=h_{2}(n)=h_{3}(n)=1$,
$h_{4}(n)=26n^2+117n+132$,
$h_{5}(n)=94n^2+517n+715$,
$h_{6}(n)=526n^2+3419n+5610$,
$h_{7}(n)=2062n^2+15465n+29393$,
$h_{8}(n)=18788n^4+319396n^3+2042275n^2+5821157n+6240360$,
$h_{9}(n)=8564n^4+162716n^3+1163679n^2+3712391n+4457400$,
and so on.
The author checked this conjecture for $0\leq m\leq 20$.
For example,
if $m=6$ and $r=1$,
then the number of $c$ in $\RCSPP{2n+1,6}$ would be equal to
\begin{align*}
&2^{-n}\prod_{k=0}^{n-1}
\frac{(6k+22)!(6k+28)!(2k+1)!(2k+8)!(2k+15)!(k+4)!}
{(4k+8)!(4k+10)!(4k+23)!(4k+25)!(2k+9)!(k+3)!}
\\&\qquad\qquad\times
\frac{(4n+15)(2062n^2+15465n+29393)}{15\cdot29393}
\end{align*}
and the first few terms are 
$\sharp\RCSPP{3,6}=3432$, 
$\sharp\RCSPP{5,6}=65934024$ and $\sharp\RCSPP{7,6}=9034911255456$.

Concerning \thetag{\ref{eq:odd_row}},
let $p^{V^\text{R}}_{n,m}(t)=\sum_{c\in\CSPP{n,m}}t^{V^\text{R}(c)}$.
At this point we can guess the explicit form of the polynomial
$p^{V^\text{R}}_{n,m}(t)$ even for $m=0$.
Let us denote $p^{V^\text{R}}_{n}(t)=p^{V^\text{R}}_{n,0}(t)$.
Then we observe
$p^{V^\text{R}}_{1}(t)=1$,
$p^{V^\text{R}}_{2}(t)=1+t$,
$p^{V^\text{R}}_{3}(t)=t^2+3t+3$,
$p^{V^\text{R}}_{4}(t)=3(t+1)p^{V^\text{R}}_{3}(t)$,
$p^{V^\text{R}}_{5}(t)=3(3t^4+18t^3+44t^2+52t+26)$,
$p^{V^\text{R}}_{6}(t)=26(t+1)p^{V^\text{R}}_{5}(t)$
and so on.

\bigbreak

Theorem~\ref{prop:TSPPtoCSPP_k} also give us a Pfaffian expression
of Conjecture\ref{conj:MT} and Conjecture~\ref{conj:refined_MT}
with the help of Theorem~\ref{thm:gf}.

\begin{corollary}
\label{cor:MT}
Let $m$ and $n\geq1$ be non-negative integers,
and let $N$ be an even integer such that $N\geq n+m-1$.
If $r$ is a positive integer such that $1\leq r\leq n+m$,
then the generating function for all plane partitions $c\in\CSPP{n,m}$
 with the weight $t^{\overline U_r(c)}$ is given by
\begin{align}
\sum_{c\in\CSPP{n,m}^{k}}t^{\overline U_r(c)}
&=\lim\limits_{\varepsilon\to0}\,\varepsilon^{-\lfloor\frac{k}{2}\rfloor}
\Pf\begin{pmatrix}
O_{n}&J_{n}B_{n,m}^{N}(t)\\
-{}^t\!B_{n,m}^{N}(t)J_{n}&\bar L_{n+N}^{(n,k)}(\varepsilon)
\end{pmatrix}.
\label{eq:gen_mt_refined}
\end{align}
Especially,
when $t=1$,
the number of elements of $\CSPP{n,m}^{k}$ is equal to
\begin{align}
\lim\limits_{\varepsilon\to0}\,\varepsilon^{-\lfloor\frac{k}{2}\rfloor}
\Pf\begin{pmatrix}
O_{n}&J_{n}B_{n,m}^{N}\\
-{}^t\!B_{n,m}^{N}J_{n}&\bar L_{n+N}^{(n,k)}(\varepsilon)
\end{pmatrix}.
\label{eq:gen_mt}
\end{align}
\end{corollary}
Note that this corollary shows that $\sum_{c\in\CSPP{n,m}^{k}}t^{\overline U_r(c)}$
does not depend on $r$.
For example,
if $n=3$, $m=0$ and $k=1$,
then the Pfaffian in the right-hand side of \thetag{\ref{eq:gen_mt_refined}} is
\[
\Pf\left( \begin {array}{ccc|ccccccc} 
0&0&0&0&0&1&1+t&t&0&0\\
0&0&0&0&1&t&0&0&0&0\\
0&0&0&1&0&0&0&0&0&0\\\hline
0&0&-1&0&\varepsilon&-\varepsilon&\varepsilon&-1&1&-1\\
0&-1&0&-\varepsilon&0&\varepsilon&-\varepsilon&1&-1&1\\
-1&-t&0&\varepsilon&-\varepsilon&0&\varepsilon&-1&1&-1\\
-1-t&0&0&-\varepsilon&\varepsilon&-\varepsilon&0&1&-1&1\\
-t&0&0&1&-1&1&-1&0&1&-1\\
0&0&0&-1&1&-1&1&-1&0&1\\
0&0&0&1&-1&1&-1&1&-1&0
\end {array}\right)
\]
which equals $({t}^{2}+2t+2)+({t}^{2}+t)\varepsilon$.
This tends to ${t}^{2}+2t+2$ when $\varepsilon\to0$.
The reader may notice that
there is another expression for 
$\sum_{c\in\CSPP{n,m}^{k}}t^{\overline U_r(c)}$
which directly follows from Theorem~\ref{prop:TSPPtoCSPP_k}
and Lemma~\ref{lem:lp}.
\begin{remark}
Let $n$ be a positive integer and let $k=0,1,\dots,n-1$.
Let $N$ be an even integer such that $N\geq k$.
Let $B_{n,m}^{(k),N}(t)=\left(b_{ij}^{(m,k)}(t)\right)_{0\leq i\leq n-1,\ 0\leq j\leq n+N-1}$
be the $n\times (n+K)$ rectangular matrix
whose $(i,j)$the entry is
\[
b_{ij}^{(m,k)}(t)=\begin{cases}
b_{ij}^{(m)}(t)
&\text{ if $0\leq j\leq n+k-1$,}\\
0
&\text{ if $j\geq n+k$}
\end{cases}
\]
where $b_{ij}^{(m,k)}(t)$ is as in \thetag{\ref{eq:entry_doubly_refined}}.
Then we have
\begin{align}
\sum_{c\in\CSPP{n,m}^{k}}t^{\overline U_r(c)}
&=
\Pf\begin{pmatrix}
O_{n}&J_{n}B_{n,m}^{(k),N}(t)\\
-{}^t\!B_{n,m}^{(k),N}(t)&\bar S_{n+N}
\end{pmatrix}.
\label{eq:gen_mt_refined2}
\end{align}
\end{remark}
This expression \thetag{\ref{eq:gen_mt_refined2}}
looks simpler than \thetag{\ref{eq:gen_mt_refined}} apparently,
but the other expression \thetag{\ref{eq:gen_mt_refined}}
will be more useful to derive the constant term identity \thetag{\ref{eq:const_mt_refined}}
in Corollary~\ref{cor:const_mt}.

\bigbreak

Next we consider the $(-1)$-enumeration of TSSCPPs.
In \cite{E},
it is shown that 
$
\sum_{c\in\CSPP{n}}(-1)^{|c|}
=\begin{cases}
A^\text{VS}_{n+2}
&\text{ if $n$ is odd,}\\
0
&\text{ otherwise.}
\end{cases}
$.
For $q$-binomial coefficient
it is well-known that
\[
\left[{n\atop r}\right]_{-1}
=\begin{cases}
0
&\text{ if $n$ is even and $r$ is odd,}\\
\binom{\lfloor n/2\rfloor}{\lfloor r/2\rfloor}
&\text{ otherwise.}
\end{cases}
\]
For the general $\CSPP{n,m}$,
we can express the $(-1)$-enumeration of $\CSPP{n,m}$ by a Pfaffian
as follows.
\begin{corollary}
\label{cor:(-1)enumeration}
Let $m\geq0$ and $n\geq1$ be non-negative integers.
Let $N$ be an even integer such that $N\geq n+m-1$.
Then
\begin{equation}
\sum_{c\in\CSPP{n,m}}(-1)^{|c|}
=\Pf\begin{pmatrix}
O_{n}&J_{n}M_{n,m}^{N}\\
-{}^t\!M_{n,m}^{N}J_{n}&\bar S_{n+N}
\end{pmatrix},
\end{equation}
where $M_{n,m}^{N}=(M^{(m)}_{i,j})_{0\leq i\leq n-1,0\leq j\leq n+N-1}$ is 
the $n\times(n+N)$ matrix defined by
\begin{equation}
M^{(m)}_{ij}
=(-1)^{\binom{j-i+1}2}\left[{{m+i}\atop{j-i}}\right]_{-1}.
\end{equation}
\end{corollary}
From this Pfaffian expression we can see that
\[
\sum_{c\in\CSPP{n,1}}(-1)^{|c|}=\sum_{c\in\CSPP{n+1}}(-1)^{|c|}
\]
(this identity is also trivial from the definition) and
\[
\sum_{c\in\CSPP{n,3}}(-1)^{|c|}=\sum_{c\in\CSPP{n+3}}(-1)^{|c|}.
\]
For $m=2,4$ we can observe the following conjecture:
\begin{equation}
\sum_{c\in\CSPP{n,2}}(-1)^{|c|}
=\begin{cases}
3^{k}\prod_{i=0}^{k-1}
\frac{(6i+4)!(3i+5)!(2i+1)!(2i+3)!(i+1)!}{(4i+3)!(4i+6)!(3i+3)!(2i)!(i+2)!}
&\text{ if $n=2k$ is even,}\\
0
&\text{ if $n$ is odd,}
\end{cases}.
\label{eq:conj-1}
\end{equation}
and
\begin{equation}
\sum_{c\in\CSPP{n,4}}(-1)^{|c|}
=\sum_{c\in\CSPP{n+2,2}}(-1)^{|c|}
\label{eq:conj-1-2}
\end{equation}
would hold.
For example,
the first few terms of $\sum_{c\in\CSPP{n,2}}(-1)^{|c|}$ are
$1$, $4$, $50$, $1862$, $202860$, and so on.
For $m\geq5$, it seems hard to guess the explicit form of 
$\sum_{c\in\CSPP{n,m}}(-1)^{|c|}$.

From Theorem~\ref{thm:strange},
$
\sum_{c\in\CSPP{n,m}}(-1)^{|c|}\,t^{\overline U_{r}(c)}
$
and
$
\sum_{c\in\CSPP{n,m}}(-1)^{|c|}\,t^{V^\text{C}(c)}
$
can be regarded
as refined $(-1)$-enumeration of TSSCPPs.
Meanwhile,
in this case,
we should note that
$
\sum_{c\in\CSPP{n,m}}(-1)^{|c|}\,t^{\overline U_{r}(c)}
$
does depend on $r$.
Recall that
$\overline U_{1}(c)$ is the number of $1$'s in $c$,
and
$\overline U_{n+m}(c)$ is the number of saturated parts in $c$.
These two cases can be easily described as in
the following corollary (i)(ii)
which is obtained as a corollary of Theorem~\ref{thm:gf}.

\begin{corollary}
\label{cor:refined-1}
Let $m$ and $n\geq1$ be non-negative integers,
and let $N$ be an even integer such that $N\geq n+m-1$.
\begin{enumerate}
\item[(i)]
We have
\begin{align}
\sum_{c\in\CSPP{n,m}} (-1)^{|c|}\, t^{\overline U_1(c)}
&=\Pf\begin{pmatrix}
O_{n}&J_{n}M_{n,m}^{\overline U_{1},N}(t)\\
-{}^t\!M_{n,m}^{\overline U_{1},N}(t)J_{n}&\bar S_{n+N}
\end{pmatrix},
\label{eq:gen_refined-1}
\end{align}
where 
$M_{n,m}^{\overline U_{1},N}(t)
=\left(M_{ij}^{(m),\overline U_{1}}(t)\right)_{0\leq i\leq n-1,0\leq j\leq n+N-1}$
is defined by
\[
M_{ij}^{(m),\overline U_{1}}(t)
=(-1)^{\binom{j-i+1}2}
\left\{
(-1)^{j-i}\left[{{m+i-1}\atop{j-i}}\right]_{-1}
+\left[{{m+i-1}\atop{j-i-1}}\right]_{-1}t
\right\}.
\]
\item[(ii)]
We have
\begin{align}
\sum_{c\in\CSPP{n,m}} (-1)^{|c|}\, t^{\overline U_{n+m}(c)}
&=\Pf\begin{pmatrix}
O_{n}&J_{n}M_{n,m}^{\overline U_{n+m},N}(t)\\
-{}^t\!M_{n,m}^{\overline U_{n+m},N}(t)J_{n}&\bar S_{n+N}
\end{pmatrix},
\label{eq:refined_saturated-1}
\end{align}
where 
$M_{n,m}^{\overline U_{n+m},N}(t)
=\left(M_{ij}^{(m),\overline U_{n+m}}(t)\right)_{0\leq i\leq n-1,0\leq j\leq n+N-1}$
is defined by
\[
M_{ij}^{(m),\overline U_{n+m}}(t)
=(-1)^{\binom{j-i+1}2}
\left\{
\left[{{m+i-1}\atop{j-i}}\right]_{-1}
+(-1)^{m+2i-j}\left[{{m+i-1}\atop{j-i-1}}\right]_{-1}t
\right\}.
\]
\item[(iii)]
We have
\begin{align}
\sum_{c\in\CSPP{n,m}} (-1)^{|c|}\, t^{V^\text{C}(c)}
&=\Pf\begin{pmatrix}
O_{n}&J_{n}M_{n,m}^{N}\\
-{}^t\!M_{n,m}^{N}J_{n}&\bar R_{n+N}(t)
\end{pmatrix}.
\label{eq:odd_column-1}
\end{align}
\end{enumerate}
\end{corollary}
If one puts $m=0$ in
\thetag{\ref{eq:gen_refined-1}}, 
\thetag{\ref{eq:refined_saturated-1}} and \thetag{\ref{eq:odd_column-1}}
and checks the first few terms of 
$\sum_{c\in\CSPP{n}} (-1)^{|c|}\, t^{\overline U_{1}(c)}$,
$\sum_{c\in\CSPP{n}} (-1)^{|c|}\, t^{\overline U_{n}(c)}$
and
$\sum_{c\in\CSPP{n}} (-1)^{|c|}\, t^{V^\text{C}(c)}$,
then he will see that neither of these polynomials equals
$A^\text{VS}_{n+2}(t)$ when $n$ is odd.
For example,
the first few terms of $\sum_{c\in\CSPP{n}} (-1)^{|c|}\, t^{V^\text{C}(c)}$
are 
$1$, 
$t-1$,
$t$,
$(t-1)  ( {t}^{2}-t+1 )$,
$t ( {t}^{2}+t+1 )$,
$( t-1)  ( {t}^{2}+1 ) ( 3\,{t}^{2}-4\,t+3 )  $,
$2t ( 2{t}^{4}+3{t}^{3}+3{t}^{2}+3t+2 ) $
and
$2( t-1)( 13{t}^{6}-20{t}^{5}+37{t}^{4}-35{t}^{3}+37{t}^{2}-20t+13 ) $.
It will be an interesting problem to find a new weight
whose distribution gives the polynomial $A^\text{VS}_{n+2}(t)$.

\bigbreak

We are now in the position to give proofs.
First we recall notation and definitions used for the lattice path method due to Gessel and Viennot \cite{GV}.
Let $D=(V,E)$ be an acyclic digraph without multiple edges.
If $u$ and $v$ are any pair of vertices,
let $\PATH{u}{v}$ denote the set of all directed $D$-paths from $u$ to $v$.
For a fixed positive integer $n$,
an \defterm{$n$-vertex} is an $n$-tuple of vertices of $D$.
If $\pmb{u}=(u_1,\dots,u_n)$ and $\pmb{v}=(v_1,\dots,v_n)$ are $n$-vertices,
an \defterm{$n$-path} from $\pmb{u}$ to $\pmb{v}$ is an $n$-tuple $\pmb{P}=(P_1,\dots,P_n)$
such that $P_i\in\PATH{u_i}{v_i}$, $i=1,\dots,n$.
The $n$-path $\pmb{P}=(P_1,\dots,P_n)$ is said to be \defterm{non-intersecting}
if any two different paths $P_i$ and $P_j$ have no vertex in common.
We will write $\PATH{\pmb{u}}{\pmb{v}}$ for the set of all $n$-paths from $\pmb{u}$ to $\pmb{v}$,
and write $\NPATH{\pmb{u}}{\pmb{v}}$ for the subset of $\PATH{\pmb{u}}{\pmb{v}}$ 
consisting of non-intersecting $n$-paths.
If $\pmb{u}=(u_1,\dots,u_m)$ and $\pmb{v}=(v_1,\dots,v_n)$ are linearly ordered sets of vertices of $D$,
then $\pmb{u}$ is said to be \defterm{$D$-compatible} with $\pmb{v}$ if
every path $P\in{\cal P}(u_i,v_l)$ intersects with every path $Q\in{\cal P}(u_j,v_k)$ whenever $i<j$ and $k<l$.
Let $S_n$ denote the symmetric group on $\{1,2,\dots,n\}$.
Then for $\pi\in S_n$,
by $\pmb{v}^\pi$ we mean the $n$ vertex $(v_{\pi(1)},\dots,v_{\pi(n)})$.

We assign a commutative indeterminate $x_e$ to each edge $e$ of $D$
and call it the weight of the edge.
Set the weight of a path $P$ to be the product of the weights of its edges and denote it by $w(P)$.
The weight $w(\pmb{P})$ of an $n$-path $\pmb{P}$ is defined to be the product
of the weights of its components.
Given any family ${\cal F}$ of edge multi-sets,
we will write $\GF{{\cal F}}$ for the generating function with respect to the weight function $w$.
Thus,
if $\pmb{u}=(u_1,\dots,u_n)$ and $\pmb{v}=(v_1,\dots,v_n)$ are $n$-vertices,
 we define the generating functions $F(\pmb{u},\pmb{v})=\GF{\PATH{\pmb{u}}{\pmb{v}}}=\sum_{\pmb{P}\in\PATH{\pmb{u}}{\pmb{v}}}w(\pmb{P})$ 
and $F_0(\pmb{u},\pmb{v})=\GF{\NPATH{\pmb{u}}{\pmb{v}}}=\sum_{\pmb{P}\in\NPATH{\pmb{u}}{\pmb{v}}}w(\pmb{P})$.
In particular, if $u$ and $v$ are any pair of vertices,
we write
\begin{equation*}
h(u,v)=\GF{\PATH{u}{v}}=\sum_{P\in\PATH{u}{v}}w(P).
\end{equation*}
The following lemma is called the Gessel-Viennot formula for counting lattice paths
in terms of determinants.
(See \cite{GV,IW2,Ste1}.)
\begin{lemma}
(Lidstr\"om-Gessel-Viennot)

Let $\pmb{u}=(u_1,\dots,u_n)$ and $\pmb{v}=(v_1,\dots,v_n)$ be two $n$-vertices in an acyclic digraph $D$.
Then
\begin{equation}
\sum_{\pi\in S_n}\sgn\pi\ F_0(\pmb{u}^{\pi},\pmb{v})=\det[h(u_i,v_j)]_{1\le i,j\le n}.
\label{eq:LGV1}
\end{equation}
In particular, if $\pmb{u}$ is $D$-compatible with $\pmb{v}$,
then
\begin{equation}
F_0(\pmb{u},\pmb{v})=\det[h(u_i,v_j)]_{1\le i,j\le n}.
\label{eq:LGV2}
\end{equation}
\end{lemma}

\begin{demo}{Proof of Lemma~\ref{lem:lp}}
We give a lattice path realization of each $c\in\CSPP{n,m}$.
Let $V=\{(x,y)\in\Bbb{N}^2:0\leq y\leq x\}$ be the vertex set,
and direct an edge from $u$ to $v$ whenever $v-u=(1,-1)$ or $(0,-1)$.
\begin{enumerate}
\item[(i)]
We assign the weight
\[
\begin{cases}
\prod_{k=j}^{N}t_{k}\cdot x_{j}
&\text{ if $j=i$,}\\
t_{j}x_{j}
&\text{ if $j<i$,}
\end{cases}
\]
to the horizontal edge from $u=(i,j)$ to $v=(i+1,j-1)$.
\item[(ii)]
We assign the weight $1$ to the vertical edge from $u=(i,j)$ to $v=(i,j-1)$.
\end{enumerate}
Let $u_j=(N-j,N-j)$
and
$v_j=(\lambda_j+N-j,0)$
for $j=1,\dots,n$,
and let $\pmb{u}=(u_1,\dots,u_n)$ and $\pmb{v}=(v_1,\dots,v_n)$.
We claim that the RCSPPs $c\in\CSPP{n,m}$ of shape $\lambda'$
can be identified as $n$-tuples of nonintersecting $D$-paths 
in $\PATH{\pmb{u}}{\pmb{v}}$.
To see this,
consider the correspondence between the following plane partition
and the paths illustrated in Figure~\ref{figure:Paths}:
\[
\young{
8&{\bf 8}&{\bf 7}&5&{\bf 5}&3&{\bf 3}\\
7&7&6&3&3&2\\
5&5&5&2&2\\
3&2&2&1&1\\
2&1&1\\
1\\}
\]
%
%
%
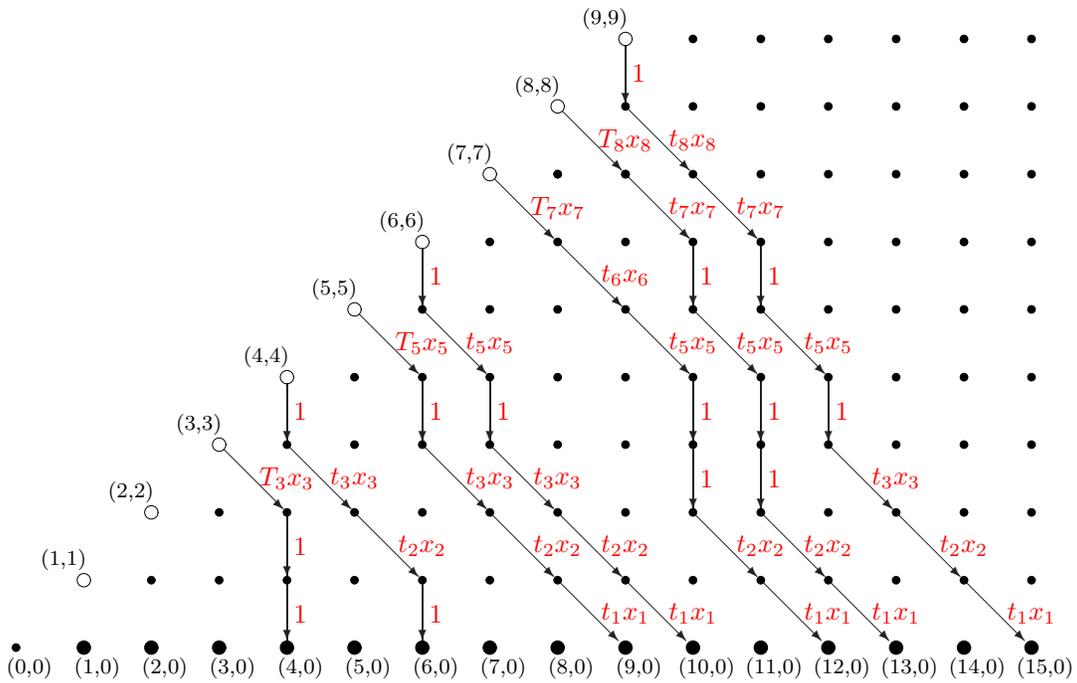
\begin{figure}[b]
\begin{center}
\setlength{\unitlength}{0.9mm}
\begin{picture}(150,100)
%
%
\put(30.7, 34.3){\vector(1,-1){8.8}}
\put( 40, 25){\vector(0,-1){9.5}}
\put( 40, 15){\vector(0,-1){9.0}}
\put( 40, 44){\vector(0,-1){8.5}}
\put( 40, 35){\vector(1,-1){9.5}}
\put( 50, 25){\vector(1,-1){9.5}}
\put( 60, 15){\vector(0,-1){9.0}}
\put(50.7, 54.3){\vector(1,-1){8.8}}
\put( 60, 45){\vector(0,-1){9.5}}
\put( 60, 35){\vector(1,-1){9.5}}
\put( 70, 25){\vector(1,-1){9.5}}
\put( 80, 15){\vector(1,-1){9.0}}
\put( 60, 64){\vector(0,-1){8.5}}
\put( 60, 55){\vector(1,-1){9.5}}
\put( 70, 45){\vector(0,-1){9.5}}
\put( 70, 35){\vector(1,-1){9.5}}
\put( 80, 25){\vector(1,-1){9.5}}
\put( 90, 15){\vector(1,-1){9.0}}
\put(70.7, 74.3){\vector(1,-1){8.8}}
\put( 80, 65){\vector(1,-1){9.5}}
\put( 90, 55){\vector(1,-1){9.5}}
\put(100, 45){\vector(0,-1){9.5}}
\put(100, 35){\vector(0,-1){9.5}}
\put(100, 25){\vector(1,-1){9.5}}
\put(110, 15){\vector(1,-1){9.0}}
\put(80.7, 84.3){\vector(1,-1){8.8}}
\put( 90, 75){\vector(1,-1){9.5}}
\put(100, 65){\vector(0,-1){9.5}}
\put(100, 55){\vector(1,-1){9.5}}
\put(110, 45){\vector(0,-1){9.5}}
\put(110, 35){\vector(0,-1){9.5}}
\put(110, 25){\vector(1,-1){9.5}}
\put(120, 15){\vector(1,-1){9.0}}
\put( 90, 94){\vector(0,-1){8.5}}
\put( 90, 85){\vector(1,-1){9.5}}
\put(100, 75){\vector(1,-1){9.5}}
\put(110, 65){\vector(0,-1){9.5}}
\put(110, 55){\vector(1,-1){9.5}}
\put(120, 45){\vector(0,-1){9.5}}
\put(120, 35){\vector(1,-1){9.5}}
\put(130, 25){\vector(1,-1){9.5}}
\put(140, 15){\vector(1,-1){9.0}}
%
\put(  0,  5){\circle*{1.2}}
\put( 10,  5){\circle*{2}}
\put( 20,  5){\circle*{2}}
\put( 30,  5){\circle*{2}}
\put( 40,  5){\circle*{2}}
\put( 50,  5){\circle*{2}}
\put( 60,  5){\circle*{2}}
\put( 70,  5){\circle*{2}}
\put( 80,  5){\circle*{2}}
\put( 90,  5){\circle*{2}}
\put(100,  5){\circle*{2}}
\put(110,  5){\circle*{2}}
\put(120,  5){\circle*{2}}
\put(130,  5){\circle*{2}}
\put(140,  5){\circle*{2}}
\put(150,  5){\circle*{2}}
\put( 10, 15){\circle{2}}
\put( 20, 25){\circle{2}}
\put( 30, 35){\circle{2}}
\put( 40, 45){\circle{2}}
\put( 50, 55){\circle{2}}
\put( 60, 65){\circle{2}}
\put( 70, 75){\circle{2}}
\put( 80, 85){\circle{2}}
\put( 90, 95){\circle{2}}
\put( 20, 15){\circle*{1.2}}
\put( 30, 15){\circle*{1.2}}
\put( 40, 15){\circle*{1.2}}
\put( 50, 15){\circle*{1.2}}
\put( 60, 15){\circle*{1.2}}
\put( 70, 15){\circle*{1.2}}
\put( 80, 15){\circle*{1.2}}
\put( 90, 15){\circle*{1.2}}
\put(100, 15){\circle*{1.2}}
\put(110, 15){\circle*{1.2}}
\put(120, 15){\circle*{1.2}}
\put(130, 15){\circle*{1.2}}
\put(140, 15){\circle*{1.2}}
\put(150, 15){\circle*{1.2}}
\put( 30, 25){\circle*{1.2}}
\put( 40, 25){\circle*{1.2}}
\put( 50, 25){\circle*{1.2}}
\put( 60, 25){\circle*{1.2}}
\put( 70, 25){\circle*{1.2}}
\put( 80, 25){\circle*{1.2}}
\put( 90, 25){\circle*{1.2}}
\put(100, 25){\circle*{1.2}}
\put(110, 25){\circle*{1.2}}
\put(120, 25){\circle*{1.2}}
\put(130, 25){\circle*{1.2}}
\put(140, 25){\circle*{1.2}}
\put(150, 25){\circle*{1.2}}
\put( 40, 35){\circle*{1.2}}
\put( 50, 35){\circle*{1.2}}
\put( 60, 35){\circle*{1.2}}
\put( 70, 35){\circle*{1.2}}
\put( 80, 35){\circle*{1.2}}
\put( 90, 35){\circle*{1.2}}
\put(100, 35){\circle*{1.2}}
\put(110, 35){\circle*{1.2}}
\put(120, 35){\circle*{1.2}}
\put(130, 35){\circle*{1.2}}
\put(140, 35){\circle*{1.2}}
\put(150, 35){\circle*{1.2}}
\put( 50, 45){\circle*{1.2}}
\put( 60, 45){\circle*{1.2}}
\put( 70, 45){\circle*{1.2}}
\put( 80, 45){\circle*{1.2}}
\put( 90, 45){\circle*{1.2}}
\put(100, 45){\circle*{1.2}}
\put(110, 45){\circle*{1.2}}
\put(120, 45){\circle*{1.2}}
\put(130, 45){\circle*{1.2}}
\put(140, 45){\circle*{1.2}}
\put(150, 45){\circle*{1.2}}
\put( 60, 55){\circle*{1.2}}
\put( 70, 55){\circle*{1.2}}
\put( 80, 55){\circle*{1.2}}
\put( 90, 55){\circle*{1.2}}
\put(100, 55){\circle*{1.2}}
\put(110, 55){\circle*{1.2}}
\put(120, 55){\circle*{1.2}}
\put(130, 55){\circle*{1.2}}
\put(140, 55){\circle*{1.2}}
\put(150, 55){\circle*{1.2}}
\put( 70, 65){\circle*{1.2}}
\put( 80, 65){\circle*{1.2}}
\put( 90, 65){\circle*{1.2}}
\put(100, 65){\circle*{1.2}}
\put(110, 65){\circle*{1.2}}
\put(120, 65){\circle*{1.2}}
\put(130, 65){\circle*{1.2}}
\put(140, 65){\circle*{1.2}}
\put(150, 65){\circle*{1.2}}
\put( 80, 75){\circle*{1.2}}
\put( 90, 75){\circle*{1.2}}
\put( 80, 75){\circle*{1.2}}
\put( 90, 75){\circle*{1.2}}
\put(100, 75){\circle*{1.2}}
\put(110, 75){\circle*{1.2}}
\put(120, 75){\circle*{1.2}}
\put(130, 75){\circle*{1.2}}
\put(140, 75){\circle*{1.2}}
\put(150, 75){\circle*{1.2}}
\put( 90, 85){\circle*{1.2}}
\put(100, 85){\circle*{1.2}}
\put(110, 85){\circle*{1.2}}
\put(120, 85){\circle*{1.2}}
\put(130, 85){\circle*{1.2}}
\put(140, 85){\circle*{1.2}}
\put(150, 85){\circle*{1.2}}
\put(100, 95){\circle*{1.2}}
\put(110, 95){\circle*{1.2}}
\put(120, 95){\circle*{1.2}}
\put(130, 95){\circle*{1.2}}
\put(140, 95){\circle*{1.2}}
\put(150, 95){\circle*{1.2}}
\put(  2,  2){\makebox(0,0){$\scriptstyle(0,0)$}}
\put( 12,  2){\makebox(0,0){$\scriptstyle(1,0)$}}
\put( 22,  2){\makebox(0,0){$\scriptstyle(2,0)$}}
\put( 32,  2){\makebox(0,0){$\scriptstyle(3,0)$}}
\put( 42,  2){\makebox(0,0){$\scriptstyle(4,0)$}}
\put( 52,  2){\makebox(0,0){$\scriptstyle(5,0)$}}
\put( 62,  2){\makebox(0,0){$\scriptstyle(6,0)$}}
\put( 72,  2){\makebox(0,0){$\scriptstyle(7,0)$}}
\put( 82,  2){\makebox(0,0){$\scriptstyle(8,0)$}}
\put( 92,  2){\makebox(0,0){$\scriptstyle(9,0)$}}
\put(102,  2){\makebox(0,0){$\scriptstyle(10,0)$}}
\put(112,  2){\makebox(0,0){$\scriptstyle(11,0)$}}
\put(122,  2){\makebox(0,0){$\scriptstyle(12,0)$}}
\put(132,  2){\makebox(0,0){$\scriptstyle(13,0)$}}
\put(142,  2){\makebox(0,0){$\scriptstyle(14,0)$}}
\put(152,  2){\makebox(0,0){$\scriptstyle(15,0)$}}
\put(  7, 18){\makebox(0,0){$\scriptstyle(1,1)$}}
\put( 17, 28){\makebox(0,0){$\scriptstyle(2,2)$}}
\put( 27, 38){\makebox(0,0){$\scriptstyle(3,3)$}}
\put( 37, 48){\makebox(0,0){$\scriptstyle(4,4)$}}
\put( 47, 58){\makebox(0,0){$\scriptstyle(5,5)$}}
\put( 57, 68){\makebox(0,0){$\scriptstyle(6,6)$}}
\put( 67, 78){\makebox(0,0){$\scriptstyle(7,7)$}}
\put( 77, 88){\makebox(0,0){$\scriptstyle(8,8)$}}
\put( 87, 98){\makebox(0,0){$\scriptstyle(9,9)$}}
\put( 40, 30){\makebox(0,0){\footnotesize\textcolor{red}{$T_3x_3$}}}
\put( 42, 20){\makebox(0,0){\footnotesize\textcolor{red}{$1$}}}
\put( 42, 10){\makebox(0,0){\footnotesize\textcolor{red}{$1$}}}
\put( 42, 40){\makebox(0,0){\footnotesize\textcolor{red}{$1$}}}
\put( 50, 30){\makebox(0,0){\footnotesize\textcolor{red}{$t_3x_3$}}}
\put( 60, 20){\makebox(0,0){\footnotesize\textcolor{red}{$t_2x_2$}}}
\put( 62, 10){\makebox(0,0){\footnotesize\textcolor{red}{$1$}}}
\put( 60, 50){\makebox(0,0){\footnotesize\textcolor{red}{$T_5x_5$}}}
\put( 62, 40){\makebox(0,0){\footnotesize\textcolor{red}{$1$}}}
\put( 70, 30){\makebox(0,0){\footnotesize\textcolor{red}{$t_3x_3$}}}
\put( 80, 20){\makebox(0,0){\footnotesize\textcolor{red}{$t_2x_2$}}}
\put( 90, 10){\makebox(0,0){\footnotesize\textcolor{red}{$t_1x_1$}}}
\put( 62, 60){\makebox(0,0){\footnotesize\textcolor{red}{$1$}}}
\put( 70, 50){\makebox(0,0){\footnotesize\textcolor{red}{$t_5x_5$}}}
\put( 72, 40){\makebox(0,0){\footnotesize\textcolor{red}{$1$}}}
\put( 80, 30){\makebox(0,0){\footnotesize\textcolor{red}{$t_3x_3$}}}
\put( 90, 20){\makebox(0,0){\footnotesize\textcolor{red}{$t_2x_2$}}}
\put(100, 10){\makebox(0,0){\footnotesize\textcolor{red}{$t_1x_1$}}}
\put( 80, 70){\makebox(0,0){\footnotesize\textcolor{red}{$T_7x_7$}}}
\put( 90, 60){\makebox(0,0){\footnotesize\textcolor{red}{$t_6x_6$}}}
\put(100, 50){\makebox(0,0){\footnotesize\textcolor{red}{$t_5x_5$}}}
\put(102, 40){\makebox(0,0){\footnotesize\textcolor{red}{$1$}}}
\put(102, 30){\makebox(0,0){\footnotesize\textcolor{red}{$1$}}}
\put(110, 20){\makebox(0,0){\footnotesize\textcolor{red}{$t_2x_2$}}}
\put(120, 10){\makebox(0,0){\footnotesize\textcolor{red}{$t_1x_1$}}}
\put( 90, 80){\makebox(0,0){\footnotesize\textcolor{red}{$T_8x_8$}}}
\put(100, 70){\makebox(0,0){\footnotesize\textcolor{red}{$t_7x_7$}}}
\put(102, 60){\makebox(0,0){\footnotesize\textcolor{red}{$1$}}}
\put(110, 50){\makebox(0,0){\footnotesize\textcolor{red}{$t_5x_5$}}}
\put(112, 40){\makebox(0,0){\footnotesize\textcolor{red}{$1$}}}
\put(112, 30){\makebox(0,0){\footnotesize\textcolor{red}{$1$}}}
\put(120, 20){\makebox(0,0){\footnotesize\textcolor{red}{$t_2x_2$}}}
\put(130, 10){\makebox(0,0){\footnotesize\textcolor{red}{$t_1x_1$}}}
\put( 92, 90){\makebox(0,0){\footnotesize\textcolor{red}{$1$}}}
\put(100, 80){\makebox(0,0){\footnotesize\textcolor{red}{$t_8x_8$}}}
\put(110, 70){\makebox(0,0){\footnotesize\textcolor{red}{$t_7x_7$}}}
\put(112, 60){\makebox(0,0){\footnotesize\textcolor{red}{$1$}}}
\put(120, 50){\makebox(0,0){\footnotesize\textcolor{red}{$t_5x_5$}}}
\put(122, 40){\makebox(0,0){\footnotesize\textcolor{red}{$1$}}}
\put(130, 30){\makebox(0,0){\footnotesize\textcolor{red}{$t_3x_3$}}}
\put(140, 20){\makebox(0,0){\footnotesize\textcolor{red}{$t_2x_2$}}}
\put(150, 10){\makebox(0,0){\footnotesize\textcolor{red}{$t_1x_1$}}}
\end{picture}
\caption{Lattice Paths ($n=7$, $m=3$, $\lambda'=(65^24^221)$, $T_i=\prod_{k=i}^{n+m}t_k$.)}\label{figure:Paths}
\end{center}
\end{figure}

The $j$th path $P_j$ from $u_j$ to $v_j$
corresponds to the $j$th column of $c$.
The entries in this column can be obtained by reading the second
coordinates of the horizontal steps of $P_j$ from left to right.
By \thetag{\ref{eq:stat}},
each horizontal step from $(i,j)$ to $(i+1,j-1)$ contributes $t_j$ to
$\pmb{t}^{\overline U(c)}$ if $j<i$,
and contributes $T_j$ to
$\pmb{t}^{\overline U(c)}$ if $j=i$,
and this is realized by the above asigned weights to each edge.
Note that the generating function $h(u_i,v_j)$ is an elementary symmetric function,
i.e.
\[
h(u_i,v_j)
=e^{(N-i)}_{\lambda_{j}-j+i}(t_{1}x_{1},\dots,t_{N-i-1}x_{N-i-1},T_{N-i}x_{N-i}).
\]
Thus we obtain the desired result \thetag{\ref{eq:GenFunc}} from \thetag{\ref{eq:LGV2}}.
This completes the proof.
\end{demo}
Now we are in the position to give a proof of Theorem~\ref{thm:gf}.

\begin{demo}{Proof of Theorem~\ref{thm:gf}}
If $(M_{ij})_{1\leq i,j\leq n}$ is any $n\times n$ matrix,
we have $\det(M_{ij})=\det(M_{n+1-i,n+1-j})$ in general.
Thus the determinant in the right-hand side of \thetag{\ref{eq:GenFunc}}
is equal to
\[
\det\left(e^{(m+i)}_{\lambda_{n-j}+j-i}(t_{1}x_{1},\dots,t_{m+i-1}x_{m+i-1},T_{m+i}x_{m+i})\right)_{0\leq i,j\leq n-1}.
\]
Note that the column indices are $\{\lambda_{n-j}+j|0\leq j\leq n-1\}=I_{n}(\lambda)$,
and the weight of this determinant is $\Pf\left(A^{I_{n}(\lambda)}_{I_{n}(\lambda)}\right)$.
If we take the sum over all partitions,
we obtain the desired identity \thetag{\ref{eq:General}}
from the minor summation formula \thetag{\ref{eq:msf}}.
This complete the proof.
\end{demo}
\begin{demo}{Proof of Corollary~\ref{cor:doubly_refined}}
To prove (i),
substitute $t_1=t$, $t_k=u$, $t_i=1$ ($i\neq 1,k$),
$x_i=1$ ($1\leq i\leq n+m$),
and $A=\bar S_{n+N}$ into \thetag{\ref{eq:General}}.
To prove (ii),
substitute $t_1=t$, $t_i=1$ ($i=2,\dots,n+m$),
$x_i=1$ ($1\leq i\leq n+m$),
and $A=\bar C_{n+N}(t)$ into \thetag{\ref{eq:General}}.
\end{demo}
\begin{demo}{Proof of Corollary~\ref{cor:doubly_refined_prop}}
We consider the right-hand side of \thetag{\ref{eq:gen_double_refined3}}
where $N$ should be taken large enough.
Here we assume the row and column indices run over $[2n+N]$.
For example,
if $n=3$, $m=0$ and $N=7$ 
then the right-hand side of \thetag{\ref{eq:gen_double_refined3}} looks like
\[
\Pf\left(
  \begin{array}{ccc|ccccccc}
   0 &  0 &  0 &  0 &  0 &  1 & 1+t&  t &  0 &  0 \\
   0 &  0 &  0 &  0 &  1 &  t &  0 &  0 &  0 &  0 \\
   0 &  0 &  0 &  1 &  0 &  0 &  0 &  0 &  0 &  0 \\\hline
   0 &  0 & -1 &  0 &  1 & -u & u^2&-u^3& u^4&-u^5\\
   0 & -1 &  0 & -1 &  0 &  1 & -u & u^2&-u^3& u^4\\
  -1 & -t &  0 &  u & -1 &  0 &  1 & -u & u^2&-u^3\\
 -1-t&  0 &  0 &-u^2&  u & -1 &  0 & 1  & -u & u^2\\
  -t &  0 &  0 & u^3&-u^2&  u & -1 & 0  &  1 & -u \\
   0 &  0 &  0 &-u^4& u^3&-u^2&  u & -1 &  0 &  1 \\
   0 &  0 &  0 & u^5&-u^4& u^3&-u^2&  u & -1 &  0 \\
  \end{array}
\right).
\]
Add $u$ times column $2n+N-1$ to column $2n+N$,
add $u$ times column $2n+N-2$ to column $2n+N-1$,
$\dots$,
and 
add $u$ times column $n+2$ to column $n+3$.
Then,
add $tu$ times column $n+1$ to column $n+2$
if $m=0$,
or
add $u$ times column $n+1$ to column $n+2$
otherwise.
Perform the same operation on the rows.
Thus we obtain the right-hand side of \thetag{\ref{eq:gen_double_refined2}}
where $m$ is replaced by $m+1$.
In the above example,
the resulting matrix looks like
\[
\Pf\left(
  \begin{array}{ccc|ccccccc}
    0 &  0 &  0 &  0 &  0 &  1 &1+t+u&t+u+tu & tu &  0 \\
    0 &  0 &  0 &  0 &  1 & t+u& tu &  0 &  0 &  0 \\
    0 &  0 &  0 &  1 & tu &  0 &  0 &  0 &  0 &  0 \\\hline
    0 &  0 & -1 &  0 &  1 &  0 &  0 &  0 &  0 &  0 \\
    0 & -1 &  0 & -1 &  0 &  1 &  0 &  0 &  0 &  0 \\
   -1 & -t &  0 &  0 & -1 &  0 &  1 &  0 &  0 &  0 \\
-1-t-u&  0 &  0 &  0 &  0 & -1 &  0 &  1 &  0 &  0 \\
-t-u-tu& 0 &  0 &  0 &  0 &  0 & -1 &  0 &  1 &  0 \\
  -tu &  0 &  0 &  0 &  0 &  0 &  0 & -1 &  0 &  1 \\
   0  &  0 &  0 &  0 &  0 &  0 &  0 &  0 & -1 &  0 \\
  \end{array}
\right).
\]
This proves the second equality.
To prove the first equality,
we perform similar operations
on the right-hand side of \thetag{\ref{eq:gen_double_refined}}.
\end{demo}
\begin{demo}{Proof of Corollary~\ref{cor:refined}}
To prove (i),
Substitute $t_i=1$ ($1\leq i\leq n+m$),
$x_i=1$ ($1\leq i\leq n+m$)
and $A=\bar S_{n+N}$ into \thetag{\ref{eq:General}}.
To prove (ii),
Substitute $t_i=1$ ($1\leq i\leq n+m$),
$x_i=1$ ($1\leq i\leq n+m$)
and $A=\bar R_{n+N}$ into \thetag{\ref{eq:General}}
since $\SHAPE{c}'$ should be even.
The other identities can be proven similarly
using \thetag{\ref{eq:General}} and Proposition~\ref{prop:skew}(iii).
\end{demo}
\begin{demo}{Proof of Corollary~\ref{cor:MT}}
As before,
we substitute $t_1=t$, $t_i=1$ ($2\leq i\leq n+m$),
$x_i=1$ ($1\leq i\leq n+m$)
and $A=\bar L_{n+N}^{(n,k)}(\varepsilon)$
into \thetag{\ref{eq:General}}.
Proposition~\ref{prop:skew}(iii) proves our claim.
\end{demo}
\begin{demo}{Proof of Corollary~\ref{cor:refined-1}}
To prove \thetag{\ref{eq:gen_refined-1}},
we substitute $x_{i}=q^i$ for $i=1,\dots,n+m$,
$t_1=t$, and $t_k=1$ for $k=2,\dots,n+m$
into \thetag{\ref{eq:General}}.
Since 
$e^{(n)}_{r}(x_1,\dots,x_n)
=e^{(n-1)}_{r}(x_2,\dots,x_n)
+x_{1}e^{(n-1)}_{r-1}(x_2,\dots,x_n)$
and
$e^{(n)}_{r}(q,q^2,\dots,q^n)=q^{\binom{r+1}2}\left[{n\atop r}\right]_{q}$,
we obtain
\[
e^{(n)}_{r}(tq,q^2,\dots,q^n)
=q^{\frac{r(r+1)}2}
\left\{
q^r\left[{{n-1}\atop{r}}\right]_{q}
+\left[{{n-1}\atop{r-1}}\right]_{q}t
\right\}.
\]
Thus,
if we put $q=-1$,
then the entries $b^{(m)}_{ij}(\pmb{t},\pmb{x})$ in \thetag{\ref{eq:entry_b}}
becomes
\[
=(-1)^{\frac{(j-i)(j-i+1)}2}
\left\{
(-1)^{j-i}\left[{{n-1}\atop{r}}\right]_{-1}
+\left[{{n-1}\atop{r-1}}\right]_{-1}t
\right\}.
\]
If we substitute $A=\bar S_{n+N}$
then we obtain the desired identity from
Proposition~\ref{prop:skew}.
The other two identities
\thetag{\ref{eq:refined_saturated-1}}
and
\thetag{\ref{eq:odd_column-1}}
can be shown similarly.
\end{demo}

\section{
Constant term identities
}\label{sec:cti}

In \cite{Z1},
D.~Zeilberger proved
the following constant term identity.
Let $D$ be the sum of all the $n\times n$ minors of 
the $n\times(2n+m-1)$ matrix $X$ given by
\[
X_{ij}=\binom{m+i}{j-i},
\qquad
0\leq i\leq n-1,
\quad
0\leq j\leq2n+m-2,
\]
and let $C$ be the constant term of
\[
\prod_{1\leq i\leq j\leq n}\left(1-\frac{z_{i}}{z_{j}}\right)
\prod_{i=1}^{n}\left(1+\frac1{z_{i}}\right)^{m+n-i}
\prod_{i=1}^{n}\frac1{1-z_{i}}
\prod_{1\leq i\leq j\leq n}\frac1{1-z_{i}z_{j}},
\]
then $D=C$ holds.
The aim of this section is to give a generalization of this constant term identity,
which gives the constant term identities for all conjectures we treat.
In that sense,
Theorem~\ref{th:CT} gives a generalization of Zeilberger's theorem,
and, as corollaries,
we obtain Corollary~\ref{cor:const_doubly_refined} for the doubly refined TSSCPP conjecture,
Corollary~\ref{cor:const_refined} for the refined TSSCPP conjecture
and Corollary~\ref{cor:const_mt} for Conjecture~\ref{conj:MT},
which are the main results of this section.

\bigbreak

Let $m$ and $n\geq1$ be non-negative integers.
Let $N$ be an even integer such that $N\geq n+m-1$.
Let $A$ be 
an $(n+N)\times(n+N)$ skew-symmetric matrix.
Let $B_{n,m}^{N}(t)=(b^{(m)}_{ij}(t))_{0\leq i\leq n-1,0\leq j\leq n+N-1}$
be the $n\times(n+N)$ matrix defined in \thetag{\ref{eq:entry_doubly_refined}}.
Let $D_{n,m}(A,t,u)$ be the sum
\begin{equation}
\sum_{I\in\binom{[n+N]}{n}}
(-1)^{s(\overline{I},I)}
\Pf\left(
A^{\overline{I}}_{\overline{I}}
\right)
\det\left( B_{n,m}^{N}(t,u)_I\right),
\end{equation}
and we also write $D_{n,m}(A,t)$ for $D_{n,m}(A,t,1)$,
and $D_{n,m}(A)$ for $D_{n,m}(A,1)$.
Let $T$ be the $n\times (n+N)$ matrix whose $(i,j)$th entry is
$z_i^{j-1}$ for $1\leq i\leq n$ and $1\leq j\leq n+N$.
If we put
\[
G_{A}(z_1,\cdots,z_n)=
\Pf\left(
\begin{array}{cc}
O_{n}&J_{n}T\\
-{}^t\!TJ_{n}&A
\end{array}\right),
\]
then $G_{A}(z_1,\cdots,z_n)$ is an anti-symmetric polynomial
in the variables $z_1,\dots,z_n$,
so that we can write 
$G_{A}(z_1,\cdots,z_n)=F_{A}(z_1,\cdots,z_n)\prod_{1\leq i<j\leq n}(z_j-z_i)$
with a symmetric polynomial $F_{A}(z_1,\cdots,z_n)$.
For example,
if we take $\bar S_{\infty}$, $\bar R_{\infty}$ and $\bar S_{\infty}$ for $A$,
we can easily see that
\begin{align}
&F_{\bar S_{\infty}}(z_1,\cdots,z_n)
=\sum_{\lambda}s^{(n)}_{\lambda}(\pmb{x})
=\prod_{i=1}^{n}\frac1{1-z_{i}}
\prod_{1\leq i<j\leq n}\frac1{1-z_{i}z_{j}},
\label{eq:SchurSum}\\
&F_{\bar R_{\infty}}(z_1,\cdots,z_n)
=\sum_{\lambda\text{ even}}s^{(n)}_{\lambda}(\pmb{x})
=\prod_{i=1}^{n}\frac1{1-z_{i}^2}
\prod_{1\leq i<j\leq n}\frac1{1-z_{i}z_{j}},
\label{eq:SchurSumEven}\\
&F_{\bar C_{\infty}}(z_1,\cdots,z_n)
=\sum_{\lambda'\text{ even}}s^{(n)}_{\lambda}(\pmb{x})
=\prod_{1\leq i<j\leq n}\frac1{1-z_{i}z_{j}},
\end{align}
where $s^{(n)}_{\lambda}(\pmb{x})$ denotes the Schur function
in the $n$ variables $\pmb{x}=(x_1,\dots,x_n)$
corresponding to the partition $\lambda$
(see \cite[I, 5, Ex.4, 5]{M}).
It can be also shown that
\begin{align}
&F_{\bar R_{\infty}(t)}(z_1,\cdots,z_n)
=\sum_{\lambda}t^{r(\lambda)}s^{(n)}_{\lambda}(\pmb{x})
=\prod_{i=1}^{n}\frac{1+tz_{i}}{1-z_{i}^2}
\prod_{1\leq i<j\leq n}\frac1{1-z_{i}z_{j}},\\
&F_{\bar C_{\infty}(t)}(z_1,\cdots,z_n)
=\sum_{\lambda}t^{c(\lambda)}s^{(n)}_{\lambda}(\pmb{x})
=\prod_{i=1}^{n}\frac1{1-tz_{i}}\prod_{1\leq i<j\leq n}\frac1{1-z_{i}z_{j}},
\end{align}
(see \cite[I, 5, Ex.7, 8]{M}).
To derive these identities are the original motivation 
of the minor summation formula (see \cite{I1,IOW,IW1}).
I.G.~Macdonald obtained the bounded version of \thetag{\ref{eq:SchurSum}}:
\begin{align}
&\lim_{\varepsilon\to0}\varepsilon^{-\lfloor\frac{n}2\rfloor}F_{\bar L_{\infty}^{(n,k)}(\varepsilon)}(z_1,\cdots,z_n)
=\sum_{{\lambda}\atop{\lambda_{1}\leq k}}s^{(n)}_{\lambda}(\pmb{x})
=\frac
{\det(z_{i}^{j-1}-z_{i}^{k+2n-j})_{1\leq i,j\leq n}}
{\prod_{i=1}^{n}(1-z_{i})\prod_{1\leq i<j\leq n}(z_{j}-z_{i})(1-z_{i}z_{j})},
\label{eq:SchurSumBounded}
\end{align}
(see \cite[I, 5, Ex.16]{M}).
In fact,
this identity \thetag{\ref{eq:SchurSumBounded}}
can be also derived from the minor summation formula \thetag{\ref{eq:msf_even}},
\thetag{\ref{eq:msf_odd}}
and Schur's Pfaffian \cite{IOTZ}
(also see \cite{I1}).

\bigbreak

If we write
$
h^{(m)}_{i}(z,t,u)=\sum_{j\geq0}
b^{(m)}_{ij}(t,u)z^{j-i}
$
where $b^{(m)}_{ij}(t)$ is as in \thetag{\ref{eq:entry_doubly_refined}},
then we have
\begin{equation}
h^{(m)}_{i}(z,t,u)
=\begin{cases}
(1+z)^{m+i}
&\text{ if $m+i=0$,}\\
(1+z)^{m+i-1}(1+tuz)
&\text{ if $m+i=1$,}\\
(1+z)^{m+i-2}(1+tz)(1+uz)
&\text{ otherwise.}
\end{cases}
\end{equation}
We also write $h^{(m)}_{i}(z,t)$ for $h^{(m)}_{i}(z,t,1)$,
and $h^{(m)}_{i}(z)=h^{(m)}_{i}(z,1)=(1+z)^{m+i}$.
\begin{theorem}
\label{th:CT}
Let $C_{n,m}(A,t,u)$ denote the constant term of
\begin{align}
\prod_{1\leq i<j\leq n}\left(1-\frac{z_j}{z_i}\right)
\prod_{k=1}^{n}h^{(m)}_{n-k}(z_{k}^{-1},t,u)
F_{A}(z_1,\cdots,z_n).
\end{align}
Then $D_{n,m}(A,t,u)$ is equal to $C_{n,m}(A,t,u)$.
\end{theorem}
The following corollary gives an constant term expression for the doubly refined
enumeration of TSSCPPs.
See Corollary~\ref{cor:doubly_refined_prop} for the relations between
these constant terms.
\begin{corollary}
\label{cor:const_doubly_refined}
Let $m$ and $n\geq1$ be non-negative integers.
\begin{enumerate}
\item[(i)]
If $r$ is an integer such that $2\leq r\leq n+m$,
then 
$\sum_{c\in\CSPP{n,m}}t^{\overline U_1(c)}u^{\overline U_r(c)}$
is equal to
\begin{align}
\CT{\pmb{z}}
\prod_{1\leq i<j\leq n}\left(1-\frac{z_j}{z_i}\right)
\prod_{k=1}^{n}h^{(m)}_{n-k}(z_{k}^{-1},t,u)
\prod_{i=1}^{n}\frac1{1-z_{i}}
\prod_{1\leq i<j\leq n}\frac1{1-z_{i}z_{j}}.
\label{eq:const_double_refined}
\end{align}
\item[(ii)]
If $r$ is an integer such that $2\leq r\leq n+m$,
then 
$\sum_{c\in\CCSPP{n,m}}t^{\overline U_1(c)}u^{\overline U_r(c)}$
is equal to
\begin{align}
\CT{\pmb{z}}
\prod_{1\leq i<j\leq n}\left(1-\frac{z_j}{z_i}\right)
\prod_{k=1}^{n}h^{(m)}_{n-k}(z_{k}^{-1},t,u)
\prod_{i=1}^{n}\frac1{1-z_{i}^2}
\prod_{1\leq i<j\leq n}\frac1{1-z_{i}z_{j}}.
\label{eq:const_double_refined_row}
\end{align}
\item[(iii)]
If $r$ is an integer such that $1\leq r\leq n+m$,
then 
$\sum_{c\in\CSPP{n,m}}t^{\overline U_1(c)}u^{V^\text{C}(c)}$
is equal to
\begin{align}
\CT{\pmb{z}}
\prod_{1\leq i<j\leq n}\left(1-\frac{z_j}{z_i}\right)
\prod_{k=1}^{n}h^{(m)}_{n-k}(z_{k}^{-1},t)
\prod_{i=1}^{n}\frac1{1-uz_{i}}
\prod_{1\leq i<j\leq n}\frac1{1-z_{i}z_{j}}.
\label{eq:const_double_refined_column}
\end{align}
\end{enumerate}
\end{corollary}
The following corollary gives an constant term expression for the refined
enumeration of TSSCPPs.
See Corollary~\ref{cor:column_even} for the relations between
these constant terms.
\begin{corollary}
\label{cor:const_refined}
Let $m$ and $n\geq1$ be non-negative integers.
\begin{enumerate}
\item[(i)]
If $r$ is an integer such that $1\leq r\leq n+m$,
then 
$\sum_{c\in\CSPP{n,m}}t^{\overline U_r(c)}$
is equal to
\begin{align}
\CT{\pmb{z}}
\prod_{1\leq i<j\leq n}\left(1-\frac{z_j}{z_i}\right)
\prod_{k=1}^{n}h^{(m)}_{n-k}(z_{k}^{-1},t)
\prod_{i=1}^{n}\frac1{1-z_{i}}
\prod_{1\leq i<j\leq n}\frac1{1-z_{i}z_{j}}.
\label{eq:const_refined}
\end{align}
\item[(ii)]
If $r$ is an integer such that $1\leq r\leq n+m$,
then 
$\sum_{c\in\CCSPP{n,m}}t^{\overline U_r(c)}$
is equal to
\begin{align}
\CT{\pmb{z}}
\prod_{1\leq i<j\leq n}\left(1-\frac{z_j}{z_i}\right)
\prod_{k=1}^{n}h^{(m)}_{n-k}(z_{k}^{-1},t)
\prod_{i=1}^{n}\frac1{1-z_{i}^2}
\prod_{1\leq i<j\leq n}\frac1{1-z_{i}z_{j}}.
\label{eq:const_refined_column}
\end{align}
\item[(iii)]
If $r$ is an integer such that $1\leq r\leq n+m$,
then 
$\sum_{c\in\RCSPP{n,m}}t^{\overline U_r(c)}$
is equal to
\begin{align}
\CT{\pmb{z}}
\prod_{1\leq i<j\leq n}\left(1-\frac{z_j}{z_i}\right)
\prod_{k=1}^{n}h^{(m)}_{n-k}(z_{k}^{-1},t)
\prod_{1\leq i<j\leq n}\frac1{1-z_{i}z_{j}}.
\label{eq:const_refined_row}
\end{align}
In particular if we put $t=1$ in this equation,
then we see that $\sharp\RCSPP{n,m}$ equals
\begin{align}
\CT{\pmb{z}}
\prod_{1\leq i<j\leq n}\left(1-\frac{z_j}{z_i}\right)
\prod_{k=1}^{n}\left(1+\frac1{z_k}\right)^{n+m-k}
\prod_{i=1}^{n}\frac1{1-z_{i}^2}
\prod_{1\leq i<j\leq n}\frac1{1-z_{i}z_{j}}.
\label{eq:const_row}
\end{align}
\item[(iv)]
The generating function $\sum_{c\in\CSPP{n,m}}t^{V^\text{C}(c)}$
is equal to
\begin{align}
\CT{\pmb{z}}
\prod_{1\leq i<j\leq n}\left(1-\frac{z_j}{z_i}\right)
\prod_{k=1}^{n}\left(1+\frac1{z_k}\right)^{n+m-k}
\prod_{i=1}^{n}\frac1{1-tz_{i}}
\prod_{1\leq i<j\leq n}\frac1{1-z_{i}z_{j}}.
\label{eq:const_refined_VC}
\end{align}
\item[(v)]
The generating function $\sum_{c\in\CSPP{n,m}}t^{V^\text{R}(c)}$
is equal to
\begin{align}
\CT{\pmb{z}}
\prod_{1\leq i<j\leq n}\left(1-\frac{z_j}{z_i}\right)
\prod_{k=1}^{n}\left(1+\frac1{z_k}\right)^{n+m-k}
\prod_{i=1}^{n}\frac{1+tz_{i}}{1-z_{i}^2}
\prod_{1\leq i<j\leq n}\frac1{1-z_{i}z_{j}}.
\label{eq:const_refined_VR}
\end{align}
\end{enumerate}
\end{corollary}
The following corollary gives a constant term identity to answer
Conjecture~\ref{conj:refined_MT} and Conjecture~\ref{conj:MT}.
\begin{corollary}
\label{cor:const_mt}
Let $m$ and $n\geq1$ be non-negative integers.
If $r$ is an integer such that $1\leq r\leq n+m$,
then 
$\sum_{c\in\CSPP{n,m}^{k}}t^{\overline U_r(c)}$
is equal to
\begin{align}
\CT{\pmb{z}}
\prod_{1\leq i<j\leq n}\left(1-\frac{z_j}{z_i}\right)
\prod_{i=1}^{n}h^{(m)}_{n-i}(z_{i}^{-1},t)
\frac
{\det(z_{i}^{j-1}-z_{i}^{k+2n-j})_{1\leq i,j\leq n}}
{\prod_{i=1}^{n}(1-z_{i})\prod_{1\leq i<j\leq n}(z_{j}-z_{i})(1-z_{i}z_{j})}.
\label{eq:const_mt_refined}
\end{align}
Especially,
when $t=1$,
the number of elements of $\CSPP{n,m}^{k}$ is equal to
\begin{align}
\CT{\pmb{z}}
\prod_{1\leq i<j\leq n}\left(1-\frac{z_j}{z_i}\right)
\prod_{i=1}^{n}\left(1+\frac1{z_i}\right)^{n+m-k}
\frac
{\det(z_{i}^{j-1}-z_{i}^{k+2n-j})_{1\leq i,j\leq n}}
{\prod_{i=1}^{n}(1-z_{i})\prod_{1\leq i<j\leq n}(z_{j}-z_{i})(1-z_{i}z_{j})}.
\label{eq:const_mt}
\end{align}
\end{corollary}
Christian Krattenthaler has obtained an equivalent result to \thetag{\ref{eq:const_mt}} in \cite{Kr3}
concerning Conjecture~\ref{conj:MT}
(i.e. Conjecture~7 of \cite{MRR2}).
\begin{demo}{Proof of Theorem~\ref{th:CT}}
We use the notation $\pmb{z}=(z_1,\dots,z_{n})$
and let $\CT{\pmb{z}}$ denote the constant term in $\pmb{z}$.
We also write 
$\alpha(I)=(-1)^{s(\overline I,I)}\Pf\left(A^{\overline I}_{\overline I}\right)$
for brevity.
Then $D_{n,m}(A,t,u)$ is equal to
\begin{align*}
&\sum_{I=\{j_1,\dots,j_n\}\in\binom{[n+N]}{n}}
\alpha(I)
\det\left(
{b^{(m)}_{i-1,j_k-1}}(t,u)
\right)_{1\leq i,k\leq n}\\
&=\sum_{I=\{j_1,\dots,j_n\}\in\binom{[n+N]}{n}}
\alpha(I)
\sum_{\pi\in S_{n}}
\sgn\pi\,\prod_{k=1}^{n}
b^{(m)}_{\pi(k)-1,j_k-1}(t,u).
\end{align*}
This sum equals
\begin{align*}
&\sum_{I=\{j_1,\dots,j_n\}\in\binom{[n+N]}{n}}
\alpha(I)
\CT{\pmb{z}}\sum_{\pi\in S_{n}}
\sgn\pi\,
\prod_{k=1}^{n}
\frac{ h^{(m)}_{\pi({k})-1}(z_{\pi(k)},t,u)}{z_{k}^{j_{k}-\pi({k})}}
\\
=&\sum_{I=\{j_1,\dots,j_n\}\in\binom{[n+N]}{n}}
\alpha(I)
\CT{\pmb{z}}\prod_{k=1}^{n}
\frac{h^{(m)}_{0}(z_{k},t,u)}
{z_l^{j_l-1}}
\sum_{\pi\in S_{n}}\sgn\pi\,
\prod_{k=1}^{n}
\left\{(1+z_{k})z_{k}\right\}^{\pi(k)-1}.
\end{align*}
Using the Vandermonde determinant 
$
\sum_{\pi\in S_{n}}\sgn\pi\,
\prod_{k=1}^{n}y_k^{\pi(k)-1}=\prod_{i<j}(y_j-y_i),
$
we obtain this sum becomes
\begin{align*}
&\CT{\pmb{z}}
\prod_{k=1}^{n}h^{(m)}_{0}(z_{k},t,u)
\prod_{1\leq i<j\leq n}\left\{z_j(1+z_j)-z_i(1+z_i)\right\}
\sum_{I=\{j_1,\dots,j_n\}}
\alpha(I)
\prod_{k=1}^{n}z_k^{-j_k+1}
\\
=&\CT{\pmb{z}}
\prod_{k=1}^{n}h^{(m)}_{0}(z_{k}^{-1},t,u)
\prod_{1\leq i<j\leq n}\left\{z_j^{-1}(1+z_j^{-1})-z_i^{-1}(1+z_i^{-1})\right\}
\sum_{I=\{j_1,\dots,j_n\}}
\alpha(I)
\prod_{k=1}^{n}z_k^{j_k-1}.
\end{align*}
This identity follows since the constant term is not changed 
by the transformation $z_{l}\rightarrow z_{l}^{-1}$.
If we use the fact that 
$\CT{\pmb{z}}g(z_1,\cdots,z_n)
=\frac{1}{n!}\CT{\pmb{z}}\sum_{\sigma\in S_n}
g(z_{\sigma(1)},\cdots,z_{\sigma(n)})$
for any polynomial $g\in\Bbb{C}[z_1,\dots,z_n]$,
then the above sum becomes
\begin{align*}
&\frac{1}{n!}\CT{\pmb{z}}
\prod_{k=1}^{n}h^{(m)}_{0}(z_{k}^{-1},t,u)
\prod_{1\leq i<j\leq n}\left\{z_j^{-1}(1+z_j^{-1})-z_i^{-1}(1+z_i^{-1})\right\}\\
&\qquad\qquad\times
\sum_{\sigma\in S_n}\sgn\sigma\,
\sum_{I=\{j_1,\dots,j_n\}}
\alpha(I)
\prod_{k=1}^{n}z_{\sigma(k)}^{j_k-1}.
\end{align*}
Now use the minor summation formula \thetag{\ref{eq:msf}}
to obtain
\begin{align}
\sum_{I=\{j_1,\dots,j_n\}\in\binom{[n+N]}{n}}
(-1)^{s(\overline I,I)}\Pf\left(A^{\overline I}_{\overline I}\right)
\sum_{\sigma\in S_n}\sgn(\sigma)\prod_{k=1}^{n}z_{\sigma(k)}^{j_k-1}
=G_{A}(z_1,\dots,z_n).
\label{eq:msf_app}
\end{align}
Substituting \thetag{\ref{eq:msf_app}} into the above identity,
we see that $D_{n,m}(A,t,u)$ is equal to
\begin{align*}
&\frac{1}{n!}
\CT{\pmb{z}}
\prod_{k=1}^{n}h^{(m)}_{0}(z_{k}^{-1},t,u)
\prod_{1\leq i<j\leq n}\left\{z_j^{-1}(1+z_j^{-1})-z_i^{-1}(1+z_i^{-1})\right\}
G_{A}(z_1,\cdots,z_n)
\end{align*}
which is equal to
\begin{align*}
&\frac{1}{n!}
\CT{\pmb{z}}
\prod_{k=1}^{n}h^{(m)}_{0}(z_{k}^{-1},t,u)
G_{A}(z_1,\cdots,z_n)
\sum_{\pi\in S_n}\sgn(\pi)
\prod_{k=1}^{n}\left\{z_k^{-1}(1+z_k^{-1})\right\}^{\pi(k)-1}
\\
=&\CT{\pmb{z}}
\prod_{k=1}^{n}h^{(m)}_{0}(z_{k}^{-1},t,u)
\prod_{k=1}^{n}\left\{z_k^{-1}(1+z_k^{-1})\right\}^{k-1}
G_{A}(z_1,\cdots,z_n).
\end{align*}
Now by changing all the index $k\rightarrow n+1-k$,
we obtain
\begin{align*}
&\CT{\pmb{z}}
\prod_{k=1}^{n}h^{(m)}_{0}(t,z_{k}^{-1})
\prod_{k=1}^{n}\left\{z_k^{-1}(1+z_k^{-1})\right\}^{n-k}
F_{A}(z_1,\cdots,z_n)\prod_{1\leq i<j\leq n}(z_i-z_j)
\\
=&\CT{\pmb{z}}
\prod_{k=1}^{n}h^{(m)}_{0}(t,z_{k}^{-1})
\prod_{k=1}^{n}(1+z_k^{-1})^{n-k}
\prod_{1\leq i<j\leq n}\left(1-\frac{z_j}{z_i}\right)
\,F_{A}(z_1,\cdots,z_n).
\end{align*}
This complete the proof.
\end{demo}


\section{
Concluding remarks
}\label{sec:remark}

First of all,
we should note that the evaluations of the Pfaffians
appearing in Section~\ref{sec:gf} are still open.
Mills, Robbins and Rumsey had the simple forms of their conjectures in \cite{MRR2}
when $m=0$.
In this paper we gave the Pfaffian forms and constant term expressions of their conjectures
in Section~\ref{sec:gf} and Section~\ref{sec:cti}.
In particular the Pfaffians in Section~\ref{sec:gf} are approximately
of size $(2n+m-1)$.
But it is also possible to make them into Pfaffians of size $n$ or of size $n+1$
using \thetag{\ref{eq:msf_even}} or \thetag{\ref{eq:msf_odd}}
(cf. \cite{A2,Kr1,Ste1}).
Here we adopt the above Pfaffians since they do not depend on whether
$n$ is even or odd.
To evaluate these Pfaffians or the constant terms,
maybe one needs the other tools.
In the forth coming paper \cite{I2} which will appear soon,
we will study the other two conjectures by Mills, Robbins and Rumsey,
(i.e. \cite[Conjecture~4, Conjecture~6]{MRR2}).


\bigbreak
\noindent
{\bf Acknowledgment:}
The author would like to express his deep gratitude 
to Prof. Soichi Okada
for his valuable comments and suggestions.



\begin{thebibliography}{99}


\bibitem{A1}
G.E.~Andrews,
``Pfaff's method (I):
the Mills-Robbins-Rumsey determinant'',
{\em Discrete Math.}
{\bf  193} (1998),  43--60.

\bibitem{A2}
G.E.~Andrews,
``Plane partitions V:
the TSSCPP conjecture'',
{\em J. Combin. Theory Ser. A}
{\bf  66} (1994), 28--39.



\bibitem{AB}
G.E.~Andrews and W.H.~Burge,
``Determinant identities'',
{\em Pacific J. Math.}
{\bf  158} (1993), 1--14.



\bibitem{B}
D.M.~Bressound,
{\em Proofs and Confirmations},
Cambridge U.P.



\bibitem{F1}
P.~Di~Francesco,
``A refined Razumnov-Stroganov conjecture'',
{\tt arXiv:cond-mat/0407477}.


\bibitem{FZ1}
P.~Di~Francesco and P.~Zinn-Justin,
``Around the Razumov-Stroganov conjecture: proof of a multi-parameter sum rule'',
{\tt arXiv:math-ph/0410061}, Electron. J. Combin. {\bf  12} (2005), R6.


\bibitem{E}
T.~Eisenk\"olbl,
``(-1)-enumeration of plane partitions with complementation symmetry''
Adv. in Appl. Math.
{\bf  30} (2003), 53--95.


\bibitem{GV}
I.~Gessel and G.~Viennot,
{\em Determinants, Paths, and Plane Partitions},
preprint (1989).

\bibitem{Ha}
T.~Hashimoto,
``A central element in the universal enveloping algebra
of type $\mathsf{D}_n$ via minor summation formula of Pfaffians'',
{\tt arXiv:math.RT/0602055}.




\bibitem{I1}
M.~Ishikawa,
``Minor summation formula and a proof of Stanley's open problem'',
{\tt arXiv:math.CO/0408204}.

\bibitem{I2}
M.~Ishikawa,
``On refined enumerations of totally symmetric self-complementary
plane partitions II'',
in preparation.



\bibitem{IOTZ}
M.~Ishikawa, H.~Tagawa, S.~Okada and J.~Zeng,
``Generalizations of Cauchy's determinant and Schur's Pfaffian'',
{\tt arXiv:math.CO/0411280},
to appear in {\em Adv. in Appl. Math.}


\bibitem{IOW}
M.~Ishikawa, S.~Okada and M.~Wakayama,
``Minor summation of Pfaffians and generalized Littlewood type formulas'',
{\em J. Alg.} {\bf  183}  (1996), 193-216.


\bibitem{IW1}
M.~Ishikawa and M.~Wakayama,
``Minor summation formula of Pfaffians'',
{\em Linear and Multilinear algebra}
{\bf  39} (1995), 285--305.


\bibitem{IW2}
M.~Ishikawa and M.~Wakayama,
``Applications of the minor summation formula III:
Pl\"ucker relations, lattice pathes and Pfaffians'',
{\tt arXiv:math.CO/0312358},
{\em J. Combin. Theory Ser. A}
{\bf 113} (2006) 113-155.



\bibitem{Kn}
D.E.~Knuth,
``Overlapping Pfaffians'',
{\em Electron. J. Combin.}
{\bf  3}, 151--163.

\bibitem{Kr1}
C.~Krattenthaler,
``Determinant identities and a generalization of the number of totally symmetric self-complementary plane partitions'',
{\em Electron. J. Combin.} {\bf  4(1)} (1997), \#R27.

\bibitem{Kr2}
C.~Krattenthaler,
``Advanced determinant calculus'',
{\em Sem. Lothar. Combin.} {\bf 42} ("The Andrews Festschrift") (1999), Article B42q.

\bibitem{Kr3}
C.~Krattenthaler,
a manuscript on the magog trapezoids (private communication).

\bibitem{Ku1}
G.~Kuperberg,
``Another proof of the alternating-aign matrix conjecture'',
{\em Int. Math. Res. Not.} {\bf  3} (1996), 139--150.
{\tt arXiv:math.CO/9810091}.

\bibitem{Ku2}
G.~Kuperberg,
``An exploration of the permanent-determinant method'',
{\em Electron. J. Combin.} {\bf  5} (1998), \#R64,
{\tt arXiv:math.CO/9810091}.

\bibitem{Ku3}
G.~Kuperberg,
``Symmetry classes of alternating-sign matrices under one roof'',
{\em Ann. of Math.} (2) {\bf  156} (2002), 835-866,
{\tt arXiv:math.CO/0008184}.



\bibitem{M}
I.~G.~Macdonald,
{\em Symmetric Functions and Hall Polynomials (2nd ed.)},
Oxford Univ. Press, (1995).

\bibitem{MRR1}
W.H.~Mills, D.P.~Robbins and H.~Rumsey,
``Alternating sign matrices and descending plane partitions'',
{\em J. Combin. Theory Ser. A}
{\bf  34}, (1983), 340--359.

\bibitem{MRR2}
W.H.~Mills, D.P.~Robbins and H.~Rumsey,
``Self-complementary totally symmetric plane partitions'',
{\em J. Combin. Theory Ser. A}
{\bf  42}, (1986), 277--292.



\bibitem{O2}
S.~Okada,
``Enumeration of symmetry classes of alternating sign matrices and characters of classical groups'',
{\tt arXiv:math.CO/0308234},
to appear.



\bibitem{Ro}
D.P. Robbins, 
``Symmetry classes of alternating sign matrices'',
{\tt arXiv:math.CO/0008045}.


\bibitem{RS1}
A.V.~Razumov and Yu. G.~Stroganov,
``On refined enumerations of some symmetry classes of ASMs'',
{\tt arXiv:math-ph/0312071}.




\bibitem{Sta1}
R.P.~Stanley,
``Symmetries of plane partitions'',
{\em J. Combin. Theory Ser. A}
{\bf  43}, (1986), 103--113.

\bibitem{Sta2}
R.P.~Stanley,
{\em Enumerative combinatorics, Volume II},
Cambridge University Press,
(1999).


\bibitem{Ste1}
J.R.~Stembridge,
``Nonintersecting paths, Pfaffians, and plane partitions''
{\em Adv. math.},
{\bf 83} (1990), 96--131.

\bibitem{Ste2}
J.R.~Stembridge,
``Strange Enumerations of CSPP's and TSPP's'',
preprint.


\bibitem{Str1}
Yu.G .~Stroganov,
``A new way to deal with Izergin-Korepin determinant at root of unity''
{\tt arXiv:math-ph/0204042}.


\bibitem{Z1}
D.~Zeilberger,
``A constant term identity
featuring the ubiquitous (and mysterious)
Andrews-Mills-Robbins-Rumsey numbers'',
{\em J. Combin. Theory Ser. A}
{\bf  66} (1994), 17--27.

\bibitem{Z2}
D.~Zeilberger,
``Proof of the refined alternating sign matrix conjecture'',
{\em  New York J. Math.}
{\bf  2} (1996), 59--68.



\end{thebibliography}
\end{document}